\documentclass[11pt]{article}

\usepackage{amsfonts,amssymb}
\usepackage{a4wide}

\title{On the McCool group $M_3$ and its associated Lie algebra}
\author{V. Metaftsis and A.I. Papistas}
\date{}

\newcommand{\aut}{{\rm Aut}}
\newcommand{\x}{\chi}
\newcommand{\g}{\gamma}

\newcommand{\Z}{\mathbb Z}
\newcommand{\qed}{$\Box$}

\renewcommand{\L}{{\rm gr}}

\newcommand{\ia}[1]{{\rm I}_{#1}{\rm A}}

\newcommand{\pf}{{\noindent \it Proof.\ }}

\newtheorem{theorem}{Theorem}
\newtheorem{corollary}{Corollary}
\newtheorem{lemma}{Lemma}
\newtheorem{remark}{Remark}
\newtheorem{proposition}{Proposition}

\begin{document}

\maketitle

\begin{abstract}
We prove that the Lie algebra of the McCool group $M_3$ is torsion
free. As a result we are able to give a presentation for the Lie
algebra of $M_3$. Furthermore, $M_{3}$ is a Magnus group.
\end{abstract}


\section{Introduction}

Throughout this paper, by \lq \lq Lie algebra\rq \rq,  we mean Lie
algebra over the ring of integers $\mathbb{Z}$. Let $G$ be a
group. We denote by $(a,b)$ the commutator $(a,b)=a^{-1}b^{-1}ab$.
For a positive integer $c$, let $\g_c(G)$ be the $c$-th term of
the lower central series of $G$. The (restricted) direct sum of
the quotients $\g_c(G)/\g_{c+1}(G)$ is the {\it associated graded
Lie algebra} of $G$, $\L(G) = \bigoplus_{c \geq 1}
\g_c(G)/\g_{c+1}(G)$. The Lie bracket multiplication in $\L(G)$ is
defined as $[a \g_{c+1}(G), b\g_{d+1}(G)] = (a, b) \g_{c+d+1}(G)$,
with $a \in \g_{c}(G)$, $b \in \g_{d}(G)$ and $(a, b)
\in\g_{c+d}(G)$ and extends the multiplication linearly.

For a group $G$, we write ${\rm IA}(G)$ for the kernel of the
natural group homomorphism from ${\rm Aut}(G)$ into ${\rm
Aut}(G/G^{\prime})$ with $G^{\prime} = \gamma_{2}(G)$. For a
positive integer $c \geq 2$, the natural group epimorphism from
$G$ onto $G/\g_{c}(G)$ induces a group homomorphism $\pi_{c}$ from
the automorphism group ${\aut}(G)$ into the automorphism group
${\aut}(G/\g_{c}(G))$. Write ${\rm I}_{c}{\rm A}(G) = {\rm
Ker}\pi_{c}$. Note that $\ia{2}(G) = {\rm IA}(G)$. It is proved by
Andreadakis \cite[Theorem 1.2]{andreadakis}, that if $G$ is
residually nilpotent (that is, $\bigcap_{c\geq 1}\gamma_{c}(G) =
\{1\}$) then $\bigcap_{c \geq 2} \ia{c}(G) = \{1\}$. For a
positive integer $n$, with $n \geq 2$, we write $F_{n}$ for a free
group of rank $n$ with a free generating set $\{x_{1}, \ldots,
x_{n}\}$. It was shown by Magnus \cite{magnus}, using work of
Nielsen \cite{nielsen}, that IA$(F_{n})$ has a finite generating
set $\{\x_{ij}, \x_{ijk}: 1 \leq i, j, k \leq n; i \neq j, k; j <
k\}$, where $\x_{ij}$ maps $x_{i}\mapsto x_{i}(x_{i},x_{j})$ and
$\x_{ijk}$ maps $x_{i}\mapsto x_{i}(x^{-1}_{j},x^{-1}_{k})$, with
both $\x_{ij}$ and $\x_{ijk}$ fixing the remaining basis elements.
Let $M_{n}$ be the subgroup of IA$(F_{n})$ generated by the subset
$S=\{ \x_{ij} : 1\le i,j\le n;\ i\neq j\}$. Then $M_n$ is called
the {\it McCool group} or the {\it basis conjugating automorphisms
group}. It is easily verified that the following relations are
satisfied by the elements of $S$, provided that, in each case, the
subscripts $i, j, k, q$ occurring are distinct:
\begin{eqnarray}
(\chi_{ij}, \chi_{kj}) & = & 1 \nonumber\\
(\chi_{ij}, \chi_{kq}) & = & 1 \label{relations}\\
(\chi_{ij} \chi_{kj}, \chi_{ik}) & = & 1.\nonumber
\end{eqnarray}
It has been proved in \cite{mccool} that $M_n$ has a presentation
$\langle S\mid Z \rangle$, where $Z$ is the set of all possible
relations of the above forms. Since $\gamma_{c}(M_{n}) \subseteq
\gamma_{c}({\rm IA}(F_{n})) \subseteq {\rm I}_{c+1}{\rm A}(F_{n})$
for all $c \geq 1$, and since $F_{n}$ is residually nilpotent, we
have $\bigcap_{c \geq 1} \gamma_{c}(M_{n}) = \{1\}$ and so,
$M_{n}$ is residually nilpotent.

The study of $M_n$ is strongly connected to that of $B_n$, the braid group.
It is well known that $M_n$ and $B_n$, as subgroups of Aut$(F_n)$ intersect at $P_n$,
the pure braid group. The graded algebra gr$(P_n)$ has been studied extensively and
Kohno \cite{kohno} and Falk and Randell \cite{fr} provide an important description of
these Lie algebras. Using similar techniques, Cohen et. al. \cite{cpvw} show that the
graded Lie algebra of the upper triangular McCool group, gr$(M_n^+)$, is additively isomorphic
to the direct sum of free Lie subalgebras
$$\bigoplus L[\x_{k1},\ldots,\x_{k,k-1}]$$ with
$[\x_{kj},\x_{st}]=0$ if $\{i,j\}\cup\{s,t\}=\emptyset$,
$[\x_{kj},\x_{sj}]=0$ if $\{s,k\}\cap\{j\}=\emptyset$ and
$[\x_{ik},\x_{ij}+\x_{kj}]=0$ for $j<k<i$.

The description of the Lie algebra gr$(M_n)$, with $n \geq 2$,
which is stated as a problem in \cite{barmik}, is a non trivial
problem. Since $M_{2} = {\rm IA}(F_{2}) \cong F_{2}$, it is well
known that ${\rm gr}(M_{2})$ is a free Lie algebra of rank $2$.
For $n \geq 3$, it seems that the known techniques for gr$(P_n)$
are not able to prove analogous results. In the present paper, we
concentrate on $M_3$ and we perform a thorough analysis on its Lie
algebra gr$(M_3)$, in order to understand its structure. For that
we use the presentation of $M_3$ given by McCool in \cite{mccool}.
Our analysis implies the following.

\begin{theorem}\label{th1.1}
Let gr$(M_3)$ be the graded Lie algebra of $M_3$. Then
\begin{enumerate}
\item gr$(M_3)$ is torsion-free $\Z$-module. In particular, ${\rm gr}(M_{3}) \cong L/J$ as Lie algebras, where $L$ is a free Lie algebra of rank $6$ and $J$ is a free Lie algebra. 
\item $M_3$ is residually nilpotent and each $\g_c(M_3)/\g_{c+1}(M_3)$ is torsion free.
That is, $M_{3}$ is a Magnus group.
\end{enumerate}
\end{theorem}

In the next few lines, we briefly explain our approach to prove
Theorem \ref{th1.1}. We use Lazard elimination to decompose the
free Lie algebra $L$ on six generators $\{x_1,\ldots,x_6\}$ into
the free Lie algebras on the $\Z$-modules $V_i$ generated by
$\{x_{2i-1},x_{2i}\}$ (with $i = 1,2,3$) and the free Lie algebra
on some specific $\mathbb{Z}$-module $W$. As $W$ is graded,
$W=\bigoplus_{n\ge2}W_n$ and we decompose each $W_n$ according to
the number of generators of $V_i$ that appear in its generators.
So we describe in detail the generating sets of $W_n$. This allow
us to describe a generating set for the derived algebra $L'$.
Next, we use a Lie algebra automorphism of $L^{\prime}$ to change
the generating set of $L'$ in such a way that the relations of
$\L(M_3)$ induced by the presentation of $M_3$ become part of the
new generating set of $L'$. Again, we decompose $L'$ using the new
generating set. This decomposition, although long and tedious, is
necessary in order to help us understand $J$, the ideal of $L$
generated by the relations of $M_{3}$ viewed in $L$. We describe
the components of $J$ studying their intersection with the
components of $L'$. In Theorem \ref{th1}, our key result, we show
that the components of $J$ are direct summands of the components
of $L$. For this we need to study the homogeneous components of
$J$ by means of Lyndon polynomials, a filtration of
tensor powers and symmetric powers. As a consequence, $L/J$ is a
torsion-free $\mathbb{Z}$-module. By a result of Witt (see \cite[Theorem 2.4.2.5]{baht}), we have $J$ is a free Lie algebra. Finally, we prove that $L/J$ is
isomorphic to the Lie algebra of $M_3$. We should point out that the above method does not work for $n>3$, due to the complexity
of the calculations.

\section{Preliminary results}
In this section we give some preliminary results that will help us
decompose $\L(M_3)$.

\subsection{Some notation}

Given a free $\mathbb{Z}$-module $A$, we write $L(A)$ for the free
Lie algebra on $A$, that is the free Lie algebra on $\cal A$ where
$\cal A$ is an arbitrary $\mathbb{Z}$-basis of $A$. Thus we may
write $L(A) = L({\cal A})$. For a positive integer $c$, let
$L^{c}(A)$ denote the $c$-th homogeneous component of $L(A)$. It
is well-known that
$$
L(A) = \bigoplus_{c \geq 1}L^{c}(A).
$$
Throughout this paper, we write $L = L({\cal X})$ for the free Lie
algebra of rank $6$ with a free generating set ${\cal X} =
\{x_{1}, \ldots, x_{6}\}$. The elements of ${\cal X}$ are ordered
as $x_{1} < x_{2} < \cdots < x_{6}$. For a positive integer $c$,
we write $L^{c} = L^{c}({\cal X})$. From now on, we write
$$
\begin{array}{llllr}
y_{1} = [x_2,x_1], & y_{2} =[x_4,x_3], & y_{3} = [x_6,x_5], &  \\
y_{4} = [x_3,x_1], & y_{5} = [x_4,x_1], & y_{6} = [x_3,x_2], & y_{7} = [x_4,x_2], & ~~~~~(D)\\
y_{8} = [x_5,x_1], & y_{9} = [x_6,x_1], & y_{10} = [x_5,x_2], & y_{11} = [x_6,x_2], &\\
y_{12} = [x_5,x_3], & y_{13} = [x_6,x_3], & y_{14} = [x_5,x_4], &
y_{15} = [x_6,x_4]. &
\end{array}
$$
Let $J$ be the ideal of $L$ generated by the set
$$
{\cal V} = \{y_{1}, y_{2}, y_{3}, y_{6}+y_{7}, y_{7}+y_{5},
y_{9}+y_{8}, y_{10}+y_{8}, y_{12}+y_{13}, y_{15}+y_{13}\}.
$$

Let $F$ be a free group of rank $6$ with a free generating set
${\cal Y}_{F} = \{a_{1}, \ldots, a_{6}\}$. We order the elements
of ${\cal Y}_{F}$ as $a_{1} < a_{2} < \cdots < a_{6}$. It is well
known that $\L(F)$ is a free Lie algebra of rank $6$; freely
generated by the set $\{a_{i}F^{\prime}: i = 1, \ldots, 6\}$. The
free Lie algebras $L$ and $\L(F)$ are isomorphic to each other by a
natural isomorphism $\chi$ subject to $\chi(x_{i}) =
a_{i}F^{\prime}$, $i = 1, \ldots, 6$. From now on, we identify $L$
with $\L(F)$, and write $x_{i} = a_{i}F^{\prime}$, $i = 1, \ldots,
6$, and for $c \geq 2$,
$$
[x_{i_{1}}, \ldots, x_{i_{c}}] = (a_{i_{1}}, \ldots, a_{i_{c}})
\gamma_{c+1}(F)
$$
for all $i_{1}, \ldots, i_{c} \in \{1, \ldots, 6\}$. Furthermore,
for each $c \geq 1$, $L^{c} = \gamma_{c}(F)/\gamma_{c+1}(F)$.
Define
$$
\begin{array}{lll}
r_{1} = (a_{2}, a_{1}), r_{2} = (a_{4}, a_{3}), & & r_{3} = (a_{6}, a_{5}) \\
r_{4} = (a_{1}a_{2},a_{5}), & & r_{5} = (a_{3}a_{4}, a_{6}) \\
r_{6} = (a_{2}a_{1}, a_{4}), & & r_{7} = (a_{4}a_{3}, a_{2})  \\
r_{8} = (a_{6}a_{5}, a_{3}), & & r_{9}
= (a_{5}a_{6}, a_{1}),
\end{array}
$$
and ${\cal R} = \{r_{1}, \ldots, r_{9}\}$. Under the above
identification, we have
$$
\begin{array}{lll}
y_{1} & = & r_{1} \gamma_{3}(F), \\
y_{2} & = & r_{2} \gamma_{3}(F), \\
y_{3} & = & r_{3} \gamma_{3}(F), \\
y_{6} + y_{7} & = & r_{6} \gamma_{3}(F), \\
y_{7} + y_{5} & = & r_{7} \gamma_{3}(F), \\
y_{9} + y_{8} & = & r_{9} \gamma_{3}(F), \\
y_{10} + y_{8} & = & r_{4} \gamma_{3}(F), \\
y_{12} + y_{13} & = & r_{8} \gamma_{3}(F), \\
y_{15} + y_{13} & = & r_{5} \gamma_{3}(F).
\end{array}
$$
Let $N = {\cal R}^{F}$ be the normal closure of $\cal R$ in $F$. Thus $N$ is generated as a group
by the set $\{r^{g} = g^{-1}rg: r \in {\cal R}, g \in
F\}$. By using the presentation of $M_{3}$, we may
show that $M_{3} = F/N$. Since $r \in F^{\prime} \setminus
\gamma_{3}(F)$ for all $r \in {\cal R}$, we have $N
\subseteq F^{\prime}$ and so, $M_{3}/M_{3}^{\prime} \cong
F/F^{\prime}$.

\subsection{Lazard elimination}

For $\mathbb{Z}$-submodules $A$ and $B$ of any Lie algebra over
$\mathbb{Z}$, let $[A,B]$ be the $\mathbb{Z}$-submodule spanned by
$[a,b]$ where $a \in A$ and $b \in B$. Furthermore, $B \wr A$
denotes the $\mathbb{Z}$-submodule defined by
$$
B \wr A = B + [B,A] + [B,A,A] + \cdots.
$$
We use the left-normed convention for Lie commutators. One of our
main tools in this paper is Lazard elimination. The following
result is a version of Lazard's "Elimination Theorem" (see
\cite[Chapter 2, Section 2.9, Proposition 10]{bour}). In the form
written here it is a special case of (\cite[Lemma 2.2]{bks1}) or
(\cite[Lemma 2]{brmic}).

\begin{lemma}\label{le4}
Let $U$ and $V$ be free $\mathbb{Z}$-modules, and consider the
free Lie algebra $L(U \oplus V)$. Then $U$ and $V \wr U$ freely
generate Lie subalgebras $L(U)$ and $L(V \wr U)$, and there is a
$\mathbb{Z}$-module decomposition $L(U \oplus V) = L(U) \oplus L(V
\wr U)$. Furthermore,
$$
V \wr U = V \oplus [V,U] \oplus [V, U, U] \oplus \cdots
$$
and, for each $n \geq 0$, there is a $\mathbb{Z}$-module
isomorphism
$$
\zeta_{n}: [V, \underbrace{U, \ldots, U}_{n}] \longrightarrow V
\otimes \underbrace{U \otimes \cdots \otimes U}_{n}
$$
such that $\zeta_{n}([v, u_{1}, \ldots, u_{n}]) = v \otimes u_{1}
\otimes \cdots \otimes u_{n}$ for all $v \in V$ and $u_{1},
\ldots, v_{n} \in U$.
\end{lemma}

\vskip .120 in

As a consequence of Lemma \ref{le4} we have the following result.
This is a special case of \cite[Proof of Lemma 3]{brmic}.

\vskip .120 in

\begin{corollary}\label{c2}
For free $\mathbb{Z}$-modules $V_{1}, \ldots, V_{n}$, with $n \geq
2$, we write $L(V_{1} \oplus \cdots \oplus V_{n})$ for the free
Lie algebra on $V_{1} \oplus \cdots \oplus V_{n}$. Then there is a
$\mathbb{Z}$-module decomposition $L(V_{1} \oplus \cdots \oplus
V_{n}) = L(V_{1}) \oplus \cdots \oplus L(V_{n}) \oplus L(W)$,
where $W = W_{2} \oplus W_{3} \oplus \cdots$ such that, for all $m
\geq 2$, $W_{m}$ is the direct sum of submodules $[V_{i_{1}},
V_{i_{2}}, V_{i_{3}}, \ldots, V_{i_{m}}]$ $(i_{1} > i_{2} \leq
i_{3} \leq \cdots \leq i_{m})$. Each $[V_{i_{1}}, V_{i_{2}},
V_{i_{3}}, \ldots, V_{i_{m}}]$ is isomorphic to $V_{i_{1}} \otimes
V_{i_{2}} \otimes \cdots \otimes V_{i_{m}}$ as
$\mathbb{Z}$-module. Furthermore, $L(W)$ is the ideal of $L(V_{1}
\oplus \cdots \oplus V_{n})$ generated by the submodules
$[V_{i},V_{j}]$ with $i \neq j$.
\end{corollary}

For $i \in \{1, 2, 3\}$, let $V_{i}$ be the $\mathbb{Z}$-module
spanned by the set ${\cal V}_{i} = \{x_{2i-1}, x_{2i}\}$. Note
that ${\cal X} = {\cal V}_{1} \cup {\cal V}_{2} \cup {\cal
V}_{3}$. Thus $L = L(V_{1} \oplus V_{2} \oplus V_{3})$. By
Corollary \ref{c2}, we have
$$
L = (\bigoplus_{i=1}^{3} L(V_{i})) \oplus L(W),
$$
where $W = \bigoplus_{n \geq 2}W_{n}$ and
$$
W_{n} = \bigoplus_{i_{1} > i_{2} \leq i_{3} \leq \cdots \leq i_{n}
\atop i_{1}, \ldots, i_{n} \in \{1, 2, 3\}} [V_{i_{1}}, V_{i_{2}},
\ldots, V_{i_{n}}].
$$
Note that $L(W)$ is an ideal of $L$. For $n \geq
2$, we write
$$
\begin{array}{lll}
W_{n} & = & (\bigoplus_{\alpha + \beta + \gamma = n-2 \atop
\alpha, \beta, \gamma \geq 0} [V_{2}, V_{1}, ~_{\alpha}V_{1},
~_{\beta}V_{2},
~_{\gamma}V_{3}]) ~~\oplus \\
& & \\
&  & (\bigoplus_{\alpha_{1} + \beta_{1} + \gamma_{1} = n-2 \atop
\alpha_{1}, \beta_{1}, \gamma_{1} \geq 0} [V_{3}, V_{1},
~_{\alpha_{1}}V_{1}, ~_{\beta_{1}}V_{2},
~_{\gamma_{1}}V_{3}]) ~~\oplus \\
& & \\
& & (\bigoplus_{\delta + \epsilon = n-2 \atop \delta, \epsilon
\geq 0} [V_{3}, V_{2}, ~_{\delta}V_{2}, ~_{\epsilon}V_{3}]).
\end{array}
$$
Write
$$
[V_{2},V_{1}] = V_{211} \oplus V_{212}, ~~[V_{3},V_{1}] = V_{311}
\oplus V_{312}~~{\rm and}~~[V_{3},V_{2}] = V_{321} \oplus V_{322},
$$
where
$V_{211} = \langle y_{4}, y_{5} \rangle$, $V_{212} = \langle y_{6}, y_{7} \rangle$,
$V_{311} = \langle y_{8}, y_{11} \rangle$, $V_{312} = \langle y_{9}, y_{10} \rangle$, $V_{321} = \langle y_{13}, y_{14} \rangle$ and $V_{322} = \langle y_{12}, y_{15} \rangle$. Furthermore, we denote
$$
\begin{array}{ll}
W^{(2,1)}_{n} = & (\bigoplus_{\alpha + \beta + \gamma = n-2 \atop
\alpha \geq 1}([V_{211},~ _{\alpha}V_{1}, ~_{\beta}V_{2},
~_{\gamma}V_{3}] \oplus [V_{212},~ _{\alpha}V_{1}, ~_{\beta}V_{2},
~_{\gamma}V_{3}]))~ \oplus  \\
& \\
& (\bigoplus_{\beta_{1} + \gamma_{1} = n-2 \atop \beta_{1} \geq
1}([V_{211}, ~_{\beta_{1}}V_{2}, ~_{\gamma_{1}}V_{3}] \oplus
[V_{212},
~_{\beta_{1}}V_{2}, ~_{\gamma_{1}}V_{3}])) ~\oplus \\
& \\
& ([V_{211}, ~_{(n-2)}V_{3}] \oplus [V_{212}, ~_{(n-2)}V_{3}]),
\end{array}
$$
$$
\begin{array}{ll}
W^{(3,1)}_{n} = & (\bigoplus_{\alpha + \beta + \gamma = n-2 \atop
\alpha \geq 1}([V_{311},~ _{\alpha}V_{1}, ~_{\beta}V_{2},
~_{\gamma}V_{3}] \oplus [V_{312},~ _{\alpha}V_{1}, ~_{\beta}V_{2},
~_{\gamma}V_{3}]))~ \oplus \\
& \\
& (\bigoplus_{\beta_{1} + \gamma_{1} = n-2 \atop \beta_{1} \geq
1}([V_{311}, ~_{\beta_{1}}V_{2}, ~_{\gamma_{1}}V_{3}] \oplus
[V_{312},
~_{\beta_{1}}V_{2}, ~_{\gamma_{1}}V_{3}])) ~\oplus \\
& \\
& ([V_{311}, ~_{(n-2)}V_{3}] \oplus [V_{312}, ~_{(n-2)}V_{3}]),
\end{array}
$$
and
$$
\begin{array}{ll}
W^{(3,2)}_{n} = & (\bigoplus_{\beta + \gamma = n-2 \atop \beta
\geq 1}([V_{323}, ~_{\beta}V_{2}, ~_{\gamma}V_{3}] \oplus
[V_{324}, ~_{\beta}V_{2}, ~_{\gamma}V_{3}]))~ \oplus \\
& \\
& ([V_{323}, ~_{(n-2)}V_{3}] \oplus [V_{324}, ~_{(n-2)}V_{3}]).
\end{array}
$$
Thus, for any $n \geq 2$,
$$
W_{n} = W^{(2,1)}_{n} \oplus W^{(3,1)}_{n} \oplus W^{(3,2)}_{n}.
$$
For $(j,i) \in \{(2,1), (3,1), (3,2)\}$, let ${\cal
W}^{(j,i)}_{n}$ denote the natural $\mathbb{Z}$-basis of
$W^{(j,i)}_{n}$. Thus
$$
{\cal W}_{n} = {\cal W}^{(2,1)}_{n} \cup {\cal W}^{(3,1)}_{n} \cup
{\cal W}^{(3,2)}_{n}
$$
is a $\mathbb{Z}$-basis of $W_{n}$. Furthermore, by ${\cal W}$, we
write for the disjoint union of all ${\cal W}_{n}$ with $n \geq 2$
$$
{\cal W} = \bigcup_{n \geq 2} {\cal W}_{n}
$$
which is a $\mathbb{Z}$-basis of $W$. For the elements of ${\cal
W}_{n}$, $n \geq 2$, and so, for the elements of ${\cal W}$, we
introduce the following notation: Let $\alpha, \beta, \gamma$ be
non-negative integers with $\alpha + \beta + \gamma = n-2$, $n
\geq 2$. For $\mu \in \{3,4,5,6\}$, $\nu \in \{1,2\}$, let
$$
v^{(\mu,\nu)}_{n,(\alpha,\beta,\gamma)} =
[x_{\mu},x_{\nu},x_{1,1}, \ldots,x_{\alpha,1}, x_{1,2},
\ldots,x_{\beta,2} ,x_{1,3}, \ldots,x_{\gamma,3}]
$$
with $x_{1,1}, \ldots, x_{\alpha,1} \in {\cal V}_{1},$ $x_{1,2},
\ldots, x_{\beta,2} \in {\cal V}_{2},$ $x_{1,3}, \ldots,
x_{\gamma,3} \in {\cal V}_{3}$, and for $\lambda \in \{5,6\}$,
$\tau \in \{3,4\}$, let
$$
u^{(\lambda,\tau)}_{n,(\delta,\epsilon)} =
[x_{\lambda},x_{\tau},x_{1,2}, \ldots,x_{\delta,2},x_{1,3},
\ldots,x_{\epsilon,3}]
$$
with $x_{1,2}, \ldots, x_{\delta,2} \in {\cal V}_{2},$ $x_{1,3},
\ldots, x_{\epsilon,3} \in {\cal V}_{3}$ and $\delta + \epsilon =
n-2$. Therefore, for $n \geq 2$, we may write
\begin{flushleft}
${\cal W}_{n} = \{v^{(\mu,\nu)}_{n,(\alpha,\beta,\gamma)},
u^{(\lambda,\tau)}_{n,(\delta,\epsilon)}: \alpha, \beta, \gamma,
\delta, \epsilon \geq 0, \alpha+\beta+\gamma = \delta+\epsilon =
n-2,$
\end{flushleft}
\begin{flushright}
$\mu \in \{3, 4, 5, 6\}, \nu \in \{1,2\}, \lambda \in \{5,6\},
\tau \in \{3,4\}\}.$
\end{flushright}

\subsection{A generating set for the derived algebra $L^{\prime}$}

Since $L^{\prime} = \bigoplus_{m \geq 2}L^{m}$ and $L(W) \subseteq
L^{\prime}$, we have by the modular law
$$
L^{\prime} = (\bigoplus_{i=1}^{3}L^{\prime}(V_{i})) \oplus L(W).
$$
For a positive integer $n$, with $n \geq 2$, we write
$$
L^{n}_{\rm grad}(W) = L^{n} \cap L(W).
$$
That is, $L^{n}_{\rm grad}(W)$ is the $\mathbb{Z}$-submodule of
$L^{n}$ spanned by all Lie commutators of the form $[v_{1},
\ldots, v_{\kappa}]$ with $\kappa \geq 1$, $v_{i} \in W_{n(i)}$
and $n(1) + \cdots + n(\kappa) = n$. Note that $L^{2}_{\rm
grad}(W) = W_{2}$ and $L^{3}_{\rm grad}(W) = W_{3}$. It is easily
verified that, for $n \geq 2$,
$$
L^{n} = (\bigoplus_{i=1}^{3}L^{n}(V_{i})) \oplus L^{n}_{\rm
grad}(W).
$$
Furthermore, since $W$ is spanned by homogeneous elements, we have
$$
L(W) = \bigoplus_{n \geq 2} L^{n}_{\rm grad}(W).
$$

The following result is well known (see, for example, \cite{baht},
\cite{shme}). It gives us a way of constructing a
$\mathbb{Z}$-basis of a free Lie algebra.

\begin{lemma}\label{le6}
Let $L({\cal A})$ be a free Lie algebra on an ordered set ${\cal
A}$. Then $L^{\prime}({\cal A})$ is a free Lie algebra with a free
generating set ${\cal A}^{(1)}$
$$
{\cal A}^{(1)} = \{[a_{i_{1}}, \ldots, a_{i_{\kappa}}]: \kappa
\geq 2, a_{i_{1}} > a_{i_{2}} \leq a_{i_{3}} \leq \cdots \leq
a_{i_{\kappa}}, a_{i_{1}}, \ldots, a_{i_{\kappa}} \in {\cal A}\}.
$$
For a positive integer $d$, let $L^{(d)}({\cal A}) =
(L^{(d-1)}({\cal A}))^{\prime}$ with $L^{(0)}({\cal A}) = L({\cal
A})$ and $L^{(1)}({\cal A}) = L^{\prime}({\cal A})$. If ${\cal
A}^{(d)}$ is an ordered free generating set for $L^{(d)}({\cal
A})$, then ${\cal A}^{(d+1)}$ is a free generating set for
$L^{(d+1)}({\cal A})$ where
$$
{\cal A}^{(d+1)} = \{[a^{(d)}_{i_{1}}, \ldots,
a^{(d)}_{i_{\kappa}}]: \kappa \geq 2, a^{(d)}_{i_{1}} >
a^{(d)}_{i_{2}} \leq a^{(d)}_{i_{3}} \leq \cdots \leq
a^{(d)}_{i_{\kappa}}, a^{(d)}_{i_{1}}, \ldots,
a^{(d)}_{i_{\kappa}} \in {\cal A}^{(d)}\}.
$$
\end{lemma}

By Lemma \ref{le6} (for ${\cal A} = {\cal X}$), $L^{\prime} =
L^{\prime}({\cal X})$ is a free Lie algebra with a free generating
set
$$
{\cal X}^{(1)} = \{[x_{i_{1}}, \ldots, x_{i_{r}}]: r \geq 2, i_{1}
> i_{2} \leq i_{3} \leq \cdots \leq i_{r}, i_{1}, \ldots, i_{r}
\in \{1, \ldots, 6\}\}.
$$
The set ${\cal X}^{(1)}$ is called \emph{the standard free
generating set} of $L^{\prime}$. For a positive integer $n$, with
$n \geq 2$, let
$$
{\cal X}_{n} = {\cal X}^{(1)} \cap L^{n}.
$$
Note that ${\cal X}^{(1)}$ decomposes into disjoint finite subsets
${\cal X}_{n}$
$$
{\cal X}^{(1)} = \bigcup_{n \geq 2} {\cal X}_{n}.
$$
For $n \geq 2$ and $i = 1, 2, 3$, let
$$
\begin{array}{lll}
{\cal X}_{n,i} & = & {\cal X}_{n} \cap L^{n}(V_{i}) = {\cal
X}^{(1)} \cap L^{n}(V_{i}) \\
{\rm and} & & \\
{\cal X}_{n,W} & = & {\cal X}_{n} \cap L^{n}_{\rm grad}(W) = {\cal
X}^{(1)} \cap L^{n}_{\rm grad}(W).
\end{array}
$$
It is easily verified that
$$
{\cal X}_{2,i} = \{[x_{2i}, x_{2i-1}]\}, ~~{\cal X}_{2,W} =
\{y_{4}, \ldots, y_{15}\}
$$
and, for $n \geq 3$, ${\cal X}_{n,i}$ consists of all Lie
commutators of the form
$$
[x_{2i},x_{2i-1},~_{\alpha}x_{2i-1}, ~_{\beta}x_{2i}]
$$
with non-negative integers $\alpha$, $\beta$ and $\alpha + \beta =
n-2$. The proof of the following result is straightforward.

\begin{lemma}\label{le7}
For $n \geq 3$, ${\cal X}_{n,W}$ consists of all Lie commutators
of the form $[x_{i_{1}}, x_{i_{2}}, \ldots,$ \linebreak
$x_{i_{n}}]$, $i_{1}
> i_{2} \leq i_{3} \leq \cdots \leq i_{n}$, $i_{1}, \ldots, i_{n}
\in \{1, \ldots, 6\}$ and either $(i_{1}, i_{2}) \in \{(2,1),
(4,3)\}$ and at least one of $i_{3}, \ldots, i_{n}$ not in
$\{i_{1}, i_{2}\}$ or $(i_{1}, i_{2}) \in \{(j,i): 1 \leq i < j
\leq 6\} \setminus \{(2,1), (4,3), (6,5)\}$.
\end{lemma}

By the $\mathbb{Z}$-module decomposition of $L^{n}$, with $n \geq
2$, ${\cal X}_{n}$ decomposes into disjoint subsets ${\cal
X}_{n,i}$, $i = 1, 2, 3$, and ${\cal X}_{n,W}$
$$
{\cal X}_{n} = (\bigcup_{i=1}^{3}{\cal X}_{n,i}) \cup {\cal
X}_{n,W}.
$$
The elements of ${\cal X}_{n,i}$ $(i = 1, 2, 3)$ are arbitrarily
ordered, and this order can be extended  to ${\cal X}_{n}$ subject
to $u < v~~{\rm if}~~u \in {\cal X}_{n,i}, v \in {\cal
X}_{n,j}~~{\rm with}~~i < j$ and $u < v$ for $u \in
\bigcup_{i=1}^{3}{\cal X}_{n,i}$ and $v \in {\cal X}_{n,W}$. Thus
$$
{\cal X}^{(1)} = \bigcup_{n \geq 2} ((\bigcup_{i=1}^{3}{\cal
X}_{n,i}) \cup {\cal X}_{n,W})
$$
is an ordered (free) generating set for $L^{\prime}$.

\section{Changing generating sets}

\subsection{A new generating set for $L^{\prime}$}\label{c3}

It is essential, in this paper, to replace a given free generating
set of a free Lie algebra by another free generating set. This can
be done by using a Lie algebra automorphism. The following result
was proved in \cite[Lemma 2.1]{bks2}.

\begin{lemma}\label{le9}
Let $\cal Z$ be a countable set, and we assume that it is
decomposed into a disjoint union ${\cal Z} = {\cal Z}_{1} \cup
{\cal Z}_{2} \cup {\cal Z}_{3} \cup \cdots$ of finite subsets
${\cal Z}_{i} = \{z_{i,1}, \ldots, z_{i,k_{i}}\}$ in such a way
that $z_{i,1} < \cdots < z_{i,k_{i}}$ and $z_{i,k_{i}} <
z_{i+1,1}$ for all $i$. Let $L({\cal Z})$ denote the free Lie
algebra on $\cal Z$. Given $z \in {\cal Z}$, we let $L_{\cal
Z}(<z)$ denote the free Lie subalgebra of $L({\cal Z})$ that is 
generated by all $x \in {\cal Z}$ with $x < z$. For each $i$, let
$\phi_{i}$ be any automorphism of the free $\mathbb{Z}$-module
spanned by ${\cal Z}_{i}$. Let $\phi: {\cal Z} \longrightarrow
L({\cal Z})$ be the map given by $\phi(z_{i,j}) =
\phi_{i}(z_{i,j}) + u_{i,j}$, where $u_{i,j} \in L_{\cal Z}(<
z_{i,1})$. Then $\phi({\cal Z})$ is a free generating set of
$L({\cal Z})$.
\end{lemma}

In this section, we apply Lemma \ref{le9} to construct a different
free generating set of $L^{\prime}$. Namely, we have the following
result. Let $L$ be the free Lie algebra over $\mathbb{Z}$ of rank
$6$, freely generated by the ordered set ${\cal X} = \{x_{1},
\ldots, x_{6}\}$ with $x_{1} < x_{2} < \cdots < x_{6}$. For $n
\geq 2$, let the map $\psi_{n}: {\cal X}_{n} \longrightarrow
L^{\prime}$ be defined as follows: For $n = 2$,
$$
\begin{array}{lll}
\psi_{2}(y_{i}) & = & y_{i}, i = 1, \ldots, 5, 8, 11, 13, 14, \\
\psi_{2}(y_{6}) & = & y_{6} + y_{7}, \\
\psi_{2}(y_{7}) & = & y_{7} + y_{5}, \\
\psi_{2}(y_{9}) & = & y_{9} + y_{8}, \\
\psi_{2}(y_{10}) & = & y_{10} + y_{8}, \\
\psi_{2}(y_{12}) & = & y_{12} + y_{13}, \\
\psi_{2}(y_{15}) & = & y_{15} + y_{13},
\end{array}
$$
and, for $n \geq 3$,
$$
\begin{array}{lll}
\psi_{n}([x_{i_{1}}, \ldots, x_{i_{n}}]) & = &
[\psi_{2}([x_{i_{1}}, x_{i_{2}}]), x_{i_{3}}, \ldots, x_{i_{n}}]
\end{array}
$$
with $i_{1} > i_{2} \leq i_{3} \leq \cdots \leq i_{n}$, $i_{1},
\ldots, i_{n} \in \{1, \ldots, 6\}$. Then $\psi_{2}$ extends to an
automorphism of $L^{2}$ and so, $\psi_{2}({\cal X}_{2})$ is a
$\mathbb{Z}$-basis of $L^{2}$. By the definition of $\psi_{n}$ and
since ${\cal X}^{(1)}$ is a free generating set for $L^{\prime}$,
there exists an automorphism $\phi_{n}$ of the free
$\mathbb{Z}$-module spanned by ${\cal X}_{n}$ such that
$$
\psi_{n}([x_{i_{1}}, x_{i_{2}}, \ldots, x_{i_{n}}]) =
[\psi_{2}([x_{i_{1}}, x_{i_{2}}]), x_{i_{3}}, \ldots, x_{i_{n}}] =
\phi_{n}([x_{i_{1}}, \ldots, x_{i_{n}}]) + u_{i_{1},
\ldots,i_{n}},
$$
where $u_{i_{1}, \ldots,i_{n}} \in L_{{\cal X}^{(1)}}(<[x_{i_{1}},
\ldots, x_{i_{n}}])$. We define the map $\Psi: {\cal X}^{(1)}
\longrightarrow L^{\prime}$ by
$$
\Psi([x_{i_{1}}, \ldots, x_{i_{n}}]) = \psi_{n}([x_{i_{1}},
\ldots, x_{i_{n}}])
$$
with $i_{1} > i_{2} \leq i_{3} \leq \cdots \leq i_{n}$ and $i_{1},
\ldots, i_{n} \in \{1, \ldots, 6\}$. The map $\Psi$ satisfies the
conditions of Lemma \ref{le9} and so, $\Psi$ is an automorphism of
$L^{\prime}$, and $\Psi({\cal X}^{(1)})$ is a free generating set
of $L^{\prime}$.

\subsection{A description of $\Psi(L(W))$}

In this section, we give a suitable, for our purposes, description
of the free Lie algebra $\Psi(L(W))$. For $n \geq 2$, we let
${\cal X}_{n,\Psi} = \Psi({\cal X}_{n})$, that is,
$$
{\cal X}_{n,\Psi} = \Psi({\cal X}_{n}) = (\bigcup_{i=1}^{3}{\cal
X}_{n,i}) \cup \Psi({\cal X}_{n,W}).
$$
By Lemma \ref{le7} and the definition of $\Psi$, we have
$\Psi({\cal X}_{n,W})( = {\cal X}_{n,\Psi,W})$ consists of all Lie
commutators of the form $[\psi_{2}([x_{i_{1}},
x_{i_{2}}]),x_{i_{3}}, \ldots, x_{i_{n}}]$ for all
$[x_{i_1},x_{i_2},\ldots,x_{i_n}]\in{\cal X}_{n,W}.$ Thus, for $n
\geq 2$, ${\cal X}_{n,\Psi}$ decomposes into disjoint subsets
${\cal X}_{n,i}$, $i = 1, 2, 3$, and ${\cal X}_{n,\Psi,W}$
$$
{\cal X}_{n,\Psi} = (\bigcup_{i=1}^{3}{\cal X}_{n,i}) \cup {\cal
X}_{n,\Psi,W}.
$$
The elements of ${\cal X}_{n,i}$ $(i = 1, 2, 3)$ and ${\cal
X}_{n,\Psi,W}$ are arbitrarily ordered, and extend it to ${\cal
X}_{n,\Psi}$ subject to $u < v~~{\rm if}~~u \in {\cal X}_{n,r}, v
\in {\cal X}_{n,s}~~{\rm with}~~r < s$ and for $u \in
\bigcup_{i=1}^{3}{\cal X}_{n,i}$ and $v \in {\cal X}_{n,\Psi,W}$
we have $u < v$. For $\kappa \in \{1,2,3\}$, we write
$$
{\cal X}_{(\kappa,\Psi)} = \bigcup_{n \geq 2}{\cal
X}_{n,\kappa}~~~{\rm and}~~~{\cal X}_{(4,\Psi)} = \bigcup_{n \geq
2}{\cal X}_{n,\Psi,W}.
$$
Note that $\Psi({\cal X}^{(1)}) = \bigcup_{\kappa=1}^{4}{\cal
X}_{(\kappa,\Psi)}$.

For $\kappa = 1, \ldots, 4$, we let $U_{(\kappa,\Psi)}$ be the
${\mathbb{Z}}$-module spanned by the set ${\cal
X}_{(\kappa,\Psi)}$. Since $L^{\prime}$ is free on $\Psi({\cal
X}^{(1)})$, and $\Psi({\cal X}^{(1)})$ is a $\mathbb{Z}$-basis of
$U_{(1,\Psi)} \oplus \cdots \oplus U_{(4,\Psi)}$, we have
$L^{\prime}$ is free on $U_{(1,\Psi)} \oplus \cdots \oplus
U_{(4,\Psi)}$ and so, by Corollary \ref{c2},
$$
\begin{array}{ll}
L^{\prime} &=
 L(U_{(1,\Psi)}) \oplus \cdots \oplus L(U_{(4,\Psi)}) \oplus L(W^{(1,
\ldots,4,\Psi)}),
\end{array}
$$
where $W^{(1,\ldots,4,\Psi)} = \bigoplus_{n \geq
2}W^{(1,\ldots,4,\Psi)}_{n}$ and, for $n \geq 2$,
$$
W^{(1,\ldots,4,\Psi)}_{n} = \bigoplus_{i_{1} > i_{2} \leq i_{3}
\leq \cdots \leq i_{n} \atop i_{1}, \ldots,i_{n} \in \{1,
2,3,4\}}[U_{(i_{1},\Psi)},U_{(i_{2},\Psi)},
\ldots,U_{(i_{n},\Psi)}].
$$
Furthermore, by Corollary \ref{c2}, $L(W^{(1,\ldots,4,\Psi)})$ is
an ideal of $L^{\prime}$. Since $L(U_{(i,\Psi)}) =
L^{\prime}(V_{i})$, $i = 1, 2, 3$, we have
$$
L^{\prime} = L^{\prime}(V_{1}) \oplus L^{\prime}(V_{2}) \oplus
L^{\prime}(V_{3}) \oplus L(U_{(4,\Psi)}) \oplus L(W^{(1,
\ldots,4,\Psi)}).
$$
It is clear that $L(U_{(4,\Psi)}) \oplus L(W^{(1,
\ldots,4,\Psi)})$ is a Lie subalgebra of $L^{\prime}$. Our aim in
this subsection is to show the following result.

\begin{proposition}\label{p3}
Let $L$ be a free Lie algebra over $\mathbb{Z}$ of rank $6$,
freely generated by the ordered set ${\cal X} = \{x_{1}, \ldots,
x_{6}\}$ with $x_{1} < x_{2} < \cdots < x_{6}$. Let $\Psi$ be the
automorphism of $L^{\prime}$ defined in subsection \ref{c3}. Then
$$
\Psi(L(W)) = L(\Psi(W)) = L(U_{(4,\Psi)}) \oplus L(W^{(1, \ldots,
4, \Psi)})
$$
as $\mathbb{Z}$-modules.
\end{proposition}

For the proof of Proposition \ref{p3}, we need some extra notation
and technical results. Write
$$
\begin{array}{lll}
\Psi([V_{2},V_{1}]) & = & [V_{2},V_{1}]^{(1)} \oplus [V_{2}, V_{1}]^{(2)}, \\
\Psi([V_{3},V_{1}]) & = & [V_{3},V_{1}]^{(1)} \oplus [V_{3}, V_{1}]^{(2)}, \\
\Psi([V_{3},V_{2}]) & = & [V_{3},V_{2}]^{(1)} \oplus [V_{3}, V_{2}]^{(2)},
\end{array}
$$
where
$$
\begin{array}{ll}
[V_{2},V_{1}]^{(1)} = \langle y_{4}, y_{5} \rangle, &
[V_{2},V_{1}]^{(2)} = \langle \psi_{2}(y_{6}), \psi_{2}(y_{7})
\rangle,
\end{array}
$$
$$
\begin{array}{ll}
[V_{3},V_{1}]^{(1)} = \langle y_{8}, y_{11} \rangle, &
[V_{3},V_{1}]^{(2)} = \langle \psi_{2}(y_{9}), \psi_{2}(y_{10})
\rangle,
\end{array}
$$
and
$$
\begin{array}{ll}
[V_{3},V_{2}]^{(1)} = \langle y_{13}, y_{14} \rangle, & [V_{3},V_{2}]^{(2)} = \langle \psi_{2}(y_{12}), \psi_{2}(y_{15}) \rangle.
\end{array}
$$
Thus
$$
\Psi(W_{2}) (= W_{2,\Psi}) = W^{(1)}_{2,\Psi} \oplus
W^{(2)}_{2,\Psi},
$$
where, for $i = 1, 2$, $W^{(i)}_{2,\Psi} = [V_{2},V_{1}]^{(i)} \oplus
[V_{3},V_{1}]^{(i)} \oplus [V_{3},V_{2}]^{(i)}$. For a positive integer $n$, with $n
\geq 3$, non-negative integers $\alpha, \beta, \gamma$ with
$\alpha+\beta+\gamma = n-2$, $(j,i) \in  \{(2,1), (3,1)\}$ and $\mu = 1, 2$, we
put
$$
W^{(j,i,\mu)}_{n,\Psi,(\alpha,\beta,\gamma)} = [[V_{j},V_{i}]^{(\mu)},
~_{\alpha}V_{1}, ~_{\beta}V_{2}, ~_{\gamma}V_{3}]
$$
and for non-negative integers $\delta, \epsilon$ with $\delta +
\epsilon = n-2$
$$
W^{(3,2,\mu)}_{n,\Psi,(\delta,\epsilon)} = [[V_{3},V_{2}]^{(\mu)},
~_{\delta}V_{2}, ~_{\epsilon}V_{3}].
$$
For $(j,i) \in \{(2,1), (3,1)\}$, let
$$
W^{(j,i)}_{n,\Psi,(\alpha,\beta,\gamma)} = W^{(j,i,1)}_{n,\Psi,(\alpha,\beta,\gamma)} + W^{(j,i,2)}_{n,\Psi,(\alpha,\beta,\gamma)}
$$
and
$$
W^{(3,2)}_{n,\Psi,(\delta,\epsilon)} = W^{(3,2,1)}_{n,\Psi,(\delta,\epsilon)} + W^{(3,2,2)}_{n,\Psi,(\delta,\epsilon)}.
$$
Since ${\cal W} = \bigcup_{n \geq 2}{\cal W}_{n}$ is a $\mathbb{Z}$-basis of $W$ and since $\Psi$ is an automorphism of $L^{\prime}$, we have the above sums are direct, and further, the sum of $W^{(j,i)}_{n,\Psi,(\alpha,\beta,\gamma)}$ over all $3$-tuples $(\alpha, \beta, \gamma)$ with $\alpha + \beta + \gamma = n-2$ is direct denoted by $W^{(j,i)}_{n,\Psi}$. Our next technical result shows that
$$
\sum_{\delta + \epsilon = n-2 \atop \delta,~ \epsilon \geq
0}W^{(3,2)}_{n,\Psi,(\delta,\epsilon)}
$$
is in fact direct. By the action of $\psi_{2}$ on ${\cal
X}_{2}$, ${\cal W}_{n}$ is a part of a free generating set and
since ${\cal X}^{(1)}$ is a free generating of $L^{\prime}$, we have the following result.

\begin{lemma}\label{le11}
Let $n$ be a positive integer, with $n \geq 3$. Then, the sum of
$W^{(3,2)}_{n,\Psi,(\delta,\epsilon)}$ over all $2$-tuples
$(\delta,\epsilon)$ of non-negative integers $\delta, \epsilon$
with $\delta+\epsilon = n-2$ is direct, denoted by
$W^{(3,2)}_{n,\Psi}$.
\end{lemma}

By the definition of $\Psi$, we have $W^{(2,1)}_{n,\Psi} +
W^{(3,1)}_{n,\Psi} + W^{(3,2)}_{n,\Psi}$, $n \geq 2$, is direct,
denoted by $W_{n,\Psi}$. For $n \geq 2$ and $i = 1, 2$, we write
$$
W^{(i)}_{n,\Psi} = W^{(2,1,i)}_{n,\Psi} \oplus
W^{(3,1,i)}_{n,\Psi} \oplus W^{(3,2,i)}_{n,\Psi}.
$$
Since $W^{(i)}_{\Psi} = \sum_{n \geq 2} W^{(i)}_{n,\Psi}$, $i = 1,
2$, is direct, and $W^{(1)}_{\Psi} \cap W^{(2)}_{\Psi} = \{0\}$,
we denote
$$
W_{\Psi} =  W^{(1)}_{\Psi} \oplus W^{(2)}_{\Psi}.
$$
For $n \geq 2$ and $i = 1, 2$, let ${\cal W}^{(i)}_{n,\Psi}$
denote the natural ${\mathbb{Z}}$-basis of $W^{(i)}_{n,\Psi}$.
More precisely, let $\alpha, \beta, \gamma$ be nonnegative
integers with $\alpha + \beta + \gamma = n-2$, $n \geq 2$. For
$\mu \in \{3,4,5,6\}$, $\nu \in \{1,2\}$, we define
$$
v^{(\mu,\nu,\Psi)}_{n,(\alpha,\beta,\gamma)} =
[\psi_{2}([x_{\mu},x_{\nu}]),x_{1,1}, \ldots,x_{\alpha,1},
x_{1,2}, \ldots,x_{\beta,2},x_{1,3}, \ldots,x_{\gamma,3}]
$$
with $x_{1,1}, \ldots, x_{\alpha,1} \in {\cal V}_{1},$ $x_{1,2},
\ldots,x_{\beta,2} \in {\cal V}_{2}$, $x_{1,3}, \ldots,
x_{\gamma,3} \in {\cal V}_{3}$, and for $\lambda \in \{5,6\}$,
$\tau \in \{3,4\}$,
$$
u^{(\lambda,\tau,\Psi)}_{n,(\delta,\epsilon)} =
[\psi_{2}([x_{\lambda},x_{\tau}]),x_{1,2},
\ldots,x_{\delta,2},x_{1,3}, \ldots,x_{\epsilon,3}]
$$
with $x_{1,2}, \ldots, x_{\delta,2} \in {\cal V}_{2},$ $x_{1,3},
\ldots, x_{\epsilon,3} \in {\cal V}_{3}$ and $\delta + \epsilon =
n-2$. Thus, for $n \geq 2$,
\begin{flushleft}
${\cal W}^{(1)}_{n,\Psi} =
\{v^{(3,1,\Psi)}_{n,(\alpha,\beta,\gamma)},
v^{(4,1,\Psi)}_{n,(\alpha,\beta,\gamma)},
v^{(5,1,\Psi)}_{n,(\alpha,\beta,\gamma)},
v^{(6,2,\Psi)}_{n,(\alpha,\beta,\gamma)},$
\end{flushleft}
\begin{flushright}
$u^{(6,3,\Psi)}_{n,(\delta,\epsilon)}, u^{(5,4,\Psi)}_{n,(\delta,\epsilon)}: \alpha, \beta, \gamma,
\delta, \epsilon \geq 0, \alpha + \beta + \gamma = \delta +
\epsilon = n-2\}$
\end{flushright}
and
\begin{flushleft}
${\cal W}^{(2)}_{n,\Psi} =
\{v^{(3,2,\Psi)}_{n,(\alpha,\beta,\gamma)},
v^{(4,2,\Psi)}_{n,(\alpha,\beta,\gamma)},
v^{(6,1,\Psi)}_{n,(\alpha,\beta,\gamma)},
v^{(5,2,\Psi)}_{n,(\alpha, \beta, \gamma)},$
\end{flushleft}
\begin{flushright}
$u^{(5,3,\Psi)}_{n,(\delta,\epsilon)}, u^{(6,4,\Psi)}_{n,(\delta,\epsilon)}: \alpha, \beta, \gamma,
\delta, \epsilon \geq 0, \alpha + \beta + \gamma = \delta +
\epsilon = n-2\}.$
\end{flushright}
Furthermore, we write ${\cal W}^{(i)}_{\Psi} = \bigcup_{n \geq 2}
{\cal W}^{(i)}_{n,\Psi}$ ($i = 1, 2$), which is the natural
$\mathbb{Z}$-basis of $W^{(i)}_{\Psi}$ and so, ${\cal W}_{\Psi} =
{\cal W}^{(1)}_{\Psi} \cup {\cal W}^{(2)}_{\Psi}$ is a
$\mathbb{Z}$-basis of $W_{\Psi}$. We note that $W^{(2)}_{\Psi}
\subseteq J$. For $n \geq 2$, we write ${\cal W}_{n,\Psi} = {\cal
W}^{(1)}_{n,\Psi} \cup {\cal W}^{(2)}_{n,\Psi}$ and so, ${\cal
W}_{\Psi} = \bigcup_{n \geq 2} {\cal W}_{n,\Psi}$.

\begin{lemma}\label{l9}
Let $L({\cal W}_{\Psi})$ be the Lie subalgebra of $L^{\prime}$
generated by ${\cal W}_{\Psi}$. Then $L({\cal W}_{\Psi}) =
L(U_{(4,\Psi)}) \oplus L(W^{(1, \ldots,4,\Psi)}).$
\end{lemma}

\pf Recall that
$$
L^{\prime} = L^{\prime}(V_{1}) \oplus L^{\prime}(V_{2}) \oplus
L^{\prime}(V_{3}) \oplus L(U_{(4,\Psi)}) \oplus L(W^{(1,
\ldots,4,\Psi)}).
$$
We shall prove our claim into two steps.

\vskip .120 in

\emph{Step 1}. $L({\cal W}_{\Psi}) \subseteq L(U_{(4,\Psi)})
\oplus L(W^{(1, \ldots,4,{\Psi})})$

\vskip .120 in

Since $L(U_{(4,\Psi)}) \oplus L(W^{(1, \ldots,4,\Psi)})$ is a Lie
subalgebra of $L^{\prime}$, it is enough to show that
$$
{\cal W}_{\Psi} \subseteq L(U_{(4,\Psi)}) \oplus L(W^{(1,
\ldots,4,\Psi)}).
$$
Since ${\cal W}_{\Psi} = \bigcup_{n \geq 2}{\cal W}_{n,\Psi}$, it
is enough to show that, for any $n \geq 2$,
$$
{\cal W}_{n,\Psi} \subseteq L(U^{(4,\Psi)}) \oplus L(W^{(1,
\ldots,4,\Psi)}).
$$
For $n = 2,3$, it is clear that ${\cal W}_{n,\Psi} \subseteq
L(U^{(4,\Psi)}) \oplus L(W^{(1, \ldots,4,\Psi)})$. Thus we assume
that $n \geq 4$. Since $L^{\prime} = L^{\prime}(V_{1}) \oplus
L^{\prime}(V_{2}) \oplus L^{\prime}(V_{3}) \oplus L(U_{(4,\Psi)})
\oplus L(W^{(1, \ldots,4,\Psi)})$, by the definition of $\Psi$,
and since $L(W^{(1, \ldots,4,\Psi)})$ is an ideal of $L^{\prime}$,
we have each $v^{(\mu,\nu,\Psi)}_{n,(\alpha,\beta,\gamma)},
u^{(\lambda,\tau,\Psi)}_{n,(\delta,\epsilon)}$ belongs to
$L(U_{(4,\Psi)}) \oplus L(W^{(1, \ldots,4,\Psi)})$ and so, we
obtain the desired result.

\vskip .120 in

\emph{Step 2}. $L^{\prime} = L^{\prime}(V_{1}) \oplus
L^{\prime}(V_{2}) \oplus L^{\prime}(V_{3}) \oplus L({\cal
W}_{\Psi})$

\vskip .120 in

By Step 1,
$$
(\bigoplus_{i=1}^{3}L^{\prime}(V_{i})) \bigcap L({\cal W}_{\Psi})
= \{0\}.
$$
Put
$$
L^{\prime}_{\Psi} = (\bigoplus_{i=1}^{3}L^{\prime}(V_{i})) \oplus
L({\cal W}_{\Psi}).
$$
We claim that $L^{\prime} = L^{\prime}_{\Psi}$. It is enough to
show that $L^{\prime} \subseteq L^{\prime}_{\Psi}$. By Corollary
\ref{c2},
$$
L = (\bigoplus_{i=1}^{3}L(V_{i})) \oplus L(W),
$$
and since $L^{\prime}(V_{1}) \oplus L^{\prime}(V_{2}) \oplus
L^{\prime}(V_{3}) \oplus L(W) \subseteq L^{\prime}$ and $(V_{1}
\oplus V_{2} \oplus V_{3}) \cap L^{\prime} = \{0\}$, we have, by
using the modular law,
$$
L^{\prime} = (\bigoplus_{i=1}^{3}L^{\prime}(V_{i})) \oplus L(W).
$$
We shall show that $W \subseteq L({\cal W}_{\Psi})$. Since
$L({\cal W}_{\Psi})$ is a $\mathbb{Z}$-module, and ${\cal W} =
\bigcup_{n \geq 2}{\cal W}_{n}$ is a $\mathbb{Z}$-basis of $W$, it
is enough to show that
$$
{\cal W}_{n} \subseteq L({\cal W}_{\Psi})
$$
for all $n \geq 2$. Furthermore, it is enough to show that
$$
v^{(\mu,\nu)}_{n,(\alpha,\beta,\gamma)},
u^{(\lambda,\tau)}_{n,(\delta,\epsilon)} \in L({\cal W}_{\Psi})
$$
for all non-negative integers $\alpha, \beta, \ldots, \epsilon
\geq 0, \alpha+\beta+\gamma = \delta+\epsilon = n-2,$ $\mu \in
\{3, 4, 5, 6\}, \nu \in \{1,2\}, \lambda \in \{5,6\}, \tau \in
\{3,4\}.$ Recall that
$$
[x_{3},x_{2}] = y_{6} =  \psi_{2}(y_{6}) - y_{7} =
\psi_{2}(y_{6})- \psi_{2}(y_{7}) + y_{5},
$$
$$
[x_{4},x_{2}] = y_{7}  =  \psi_{2}(y_{7}) - y_{5},
$$
$$
[x_{6},x_{1}] = y_{9} =  \psi_{2}(y_{9}) - y_{8},
$$
$$
[x_{5},x_{2}] = y_{10}  =  \psi_{2}(y_{10}) - y_{8},
$$
$$
[x_{5},x_{3}] = y_{12}  =  \psi_{2}(y_{12}) - y_{13}
$$
{\rm and}
$$
[x_{6},x_{4}] = y_{15}  =  \psi_{2}(y_{15}) - y_{13}.
$$
Write ${\cal H} = \{(\mu,\nu): \mu \in \{3,4,5,6\}, \nu \in
\{1,2\}\}$. If $(\mu,\nu) \in {\cal H} \setminus
\{(3,2),(4,2),(5,2),(6,1)\}$, then
$v^{(\mu,\nu)}_{n,(\alpha,\beta,\gamma)} \in W_{n,\Psi}$ for all
non-negative integers $\alpha, \beta, \gamma$ with $\alpha + \beta
+ \gamma = n-2$. Suppose that $(\mu, \nu) = (3,2)$. Then
$$
\begin{array}{lll}
v^{(3,2)}_{n,(\alpha,\beta,\gamma)} & = & [y_{6},
x_{1,1},\ldots,x_{\alpha,1}, x_{1,},\ldots,x_{\beta,2},
x_{1,3},\ldots,x_{\gamma,3}]
\\
& = & [\psi_{2}(y_{6}), x_{1,1},\ldots,x_{\alpha,1},
x_{1,2},\ldots,x_{\beta,2},
x_{1,3},\ldots,x_{\gamma,3}]~-\\
& & [\psi_{2}(y_{6}), x_{1,1},\ldots,x_{\alpha,1},
x_{1,2},\ldots,x_{\beta,2}, x_{1,3},\ldots,x_{\gamma,3}]~+\\
& & [y_{5}, x_{1,1},\ldots,x_{\alpha,1},
x_{1,2},\ldots,x_{\beta,2}, x_{1,3},\ldots,x_{\gamma,3}]
\end{array}
$$
where $x_{1,1},\ldots,x_{\alpha,1} \in {\cal V}_{1}$,
$x_{1,2},\ldots,x_{\beta,2} \in {\cal V}_{2}$,
$x_{1,3},\ldots,x_{\gamma,3} \in {\cal V}_{3}$ and so,
$v^{(6,1)}_{n,(\alpha,\beta,\gamma)} \in W_{n,\Psi}$ for all
$\alpha, \beta, \gamma$. By using similar arguments as before, we
have $v^{(4,2)}_{n,(\alpha,\beta,\gamma)},
v^{(5,2)}_{n,(\alpha,\beta,\gamma)}$ and
$v^{(6,1)}_{n,(\alpha,\beta,\gamma)} \in W_{n,\Psi}$. Therefore,
for $(\mu,\nu) \in {\cal H}$,
$v^{(\mu,\nu)}_{n,(\alpha,\beta,\gamma)} \in W_{n,\Psi}$ for all
$\alpha, \beta, \gamma$. It is easily verified that, for
$(\lambda,\tau) \in \{(5,4), (6,3)\}$,
$u^{(\lambda,\tau)}_{n,(\delta,\epsilon)} \in W_{n,\Psi}$ for all
$\delta,\epsilon$. Thus, we concentrate on the cases $(5,3)$ and
$(6,4)$. Then
$$
\begin{array}{lll}
u^{(5,3)}_{n,(\delta,\epsilon)} & = & [y_{12}, x_{1,2},\ldots,
x_{\delta,2}, x_{1,3},\ldots,x_{\epsilon,3}] \\
& = & [\psi_{2}(y_{12}), x_{1,2},\ldots, x_{\delta,2},
x_{1,3},\ldots,x_{\epsilon,3}]~- \\
&  & [y_{13}, ,x_{1,2},\ldots, x_{\delta,2},
x_{1,3},\ldots,x_{\epsilon,3}],
\end{array}
$$
where $x_{1,2},\ldots, x_{\delta,2} \in {\cal V}_{2},
x_{1,3},\ldots,x_{\epsilon,3} \in {\cal V}_{3}$. Similarly,
$u^{(6,4)}_{n,(\delta,\epsilon)} \in W_{n,\Psi}$ for all $\delta$
and $\epsilon$. Therefore, for all $n \geq 2$, ${\cal W}_{n}
\subseteq L({\cal W}_{\Psi})$. Since $L(W)$ is generated by ${\cal
W}$, we have $L(W) \subseteq L^{\prime}_{\Psi}$ and so,
$L^{\prime} = L^{\prime}_{\Psi}$. Thus
$$
\begin{array}{lll}
({\rm by ~Step~1})~~L(W_{\Psi}) & = & L^{\prime}_{\Psi} \cap
(L(U_{(4,\Psi)}) \oplus
L(W^{(1, \ldots,4,{\Psi})})) \\
& = & L^{\prime} \cap (L(U_{(4,\Psi)}) \oplus
L(W^{(1, \ldots,4,{\Psi})})) \\
& = & L(U_{(4,\Psi)}) \oplus L(W^{(1, \ldots,4,{\Psi})})
\end{array}
$$
and so, we obtain the required result.

\begin{lemma}\label{l10}
Let $L({\cal W}_{\Psi})$ be the Lie subalgebra of $L^{\prime}$
generated by ${\cal W}_{\Psi}$. Then $\Psi(L(W)) = L(\Psi(W)) =
L({\cal W}_{\Psi})$. In particular, $L({\cal W}_{\Psi})$ is free
on ${\cal W}_{\Psi}$, and it is an ideal in $L^{\prime}$.
\end{lemma}

\pf We first show that $\Psi(L(W)) \subseteq L({\cal W}_{\Psi})$.
Let ${\cal W}$ be the natural $\mathbb{Z}$-basis of $W$. Recall
that ${\cal W} = \bigcup_{n \geq 2} {\cal W}_{n}$, where ${\cal
W}_{n}$ is the $\mathbb{Z}$-basis of $W_{n}$ and its elements are
denoted by $v^{(\mu,\nu)}_{n,(\alpha,\beta,\gamma)},
u^{(\lambda,\tau)}_{n,(\delta,\epsilon)}$. Since $\Psi({\cal W}) =
\bigcup_{n \geq 2} \Psi({\cal W}_{n})$, it is enough to show that
$\Psi(v^{(\mu,\nu)}_{n,(\alpha,\beta,\gamma)}),
\Psi(u^{(\lambda,\tau)}_{n,(\delta,\epsilon)})$ \linebreak $\in
L(W_{\Psi})$. By Lemma \ref{le6} (for ${\cal A} = {\cal X}$),
${\cal X}^{(1)}$ is a free generating set for $L^{\prime}$. Recall
that ${\cal X}^{(1)}$ is an ordered set. Furthermore, for a
positive integer $d$, if ${\cal X}^{(d)}$ is an ordered free
generating set for $L^{(d)}({\cal X})$, with $L^{(1)}({\cal X}) =
L^{\prime}$, then ${\cal X}^{(d+1)}$ is a free generating set for
$L^{(d+1)}({\cal X})$ where ${\cal X}^{(d+1)} =
\{[x^{(d)}_{i_{1}}, \ldots, x^{(d)}_{i_{\kappa}}]: \kappa \geq 2,
x^{(d)}_{i_{1}} > x^{(d)}_{i_{2}} \leq x^{(d)}_{i_{3}} \leq \cdots
\leq x^{(d)}_{i_{\kappa}}, x^{(d)}_{i_{1}}, \ldots,
x^{(d)}_{i_{\kappa}} \in {\cal X}^{(d)}\}$. Write $X^{(d)}$ for
the $\mathbb{Z}$-module spanned by ${\cal X}^{(d)}$. Note that
$L^{(d)}({\cal X})/L^{(d+1)}({\cal X})$ is isomorphic to $X^{(d)}$
as $\mathbb{Z}$-module. Furthermore, for a positive integer $n$,
with $n \geq 2$, there exists a (unique) positive integer $m(n)
\geq 2$ such that $L^{n} = \bigoplus_{j=1}^{m(n)}(L^{n} \bigcap
X^{(j)}).$ Thus, any non-zero simple Lie commutator $w =
[x_{i_{1}}, \ldots, x_{i_{n}}]$, with $x_{i_{1}}, \ldots,
x_{i_{n}} \in {\cal X}$, is written as $w =
\sum_{j=1}^{m(n)}u_{j}$, where $u_{j} \in L^{n} \cap X^{(j)}$, $j
= 1, \ldots, m(n)$. (The number of occurrence of each $x_{i} \in
{\cal X}$ in $w$ is the same in each $u_{j}$.) Thus $\Psi(w) =
\sum^{m(n)}_{j=1} \Psi(u_{j})$. By the definition of $\Psi$, we
have $\Psi(u_{j}) \in L({\cal W}_{\Psi})$ for $j = 1, \ldots,
m(n)$ and so, $\Psi(L(W)) \subseteq L({\cal W}_{\Psi})$. Since
$L^{\prime} = (\bigoplus_{i=1}^{3}L^{\prime}(V_{i})) \oplus L(W)$,
by the definition of $\Psi$, and $\Psi$ is an automorphism of
$L^{\prime}$, we have
$$
L^{\prime} = (\bigoplus_{i=1}^{3}L^{\prime}(V_{i})) \oplus
\Psi(L(W)).
$$
Since $\Psi(L(W)) \subseteq L({\cal W}_{\Psi})$, the modular law
and Lemma \ref{l9}, we get
$$
\Psi(L(W)) = L^{\prime} \cap L({\cal W}_{\Psi}) = L({\cal
W}_{\Psi}).
$$

\vskip .120 in

\emph{Proof of Proposition \ref{p3}}. Since $L(W)$ is a free Lie
algebra on $W$, $L(W)$ is an ideal of $L^{\prime}$ and $\Psi$ is
an automorphism of $L^{\prime}$, we
 have $\Psi(L(W))$ is a free Lie algebra, and $\Psi(L(W)) \subseteq L^{\prime}$.
 By Lemma \ref{l10} and Lemma \ref{l9}, we obtain the desired result. \qed

\subsection{An ordering on ${\cal W}_{\Psi}$}

The set ${\cal W}_{\Psi}$, defined in section 3.2, will play a
 significant role in the proof of our main result (Theorem
 \ref{th1}). We need to introduce a specific order on its
 elements. Let $A$ be a finite totally ordered alphabet.
 We order the free monoid $A^{*}$ with alphabetical
order, that is, $u < v$ if and only if either $v = ux$ for some $x
\in A^{+}$ (the free semigroup on $A$), or $u = xau^{\prime}$, $v
= xbv^{\prime}$ for some words $x, u^{\prime}, v^{\prime}$ and
some $a, b \in A$ with $a < b$. (Note that the empty word in
$A^{*}$ is regarded the smallest element in $A^{*}$.)We order the
elements of the natural $\mathbb{Z}$-basis of $W_{2,\Psi}$ as
follows:
\begin{flushleft}
$y_{4} \ll y_{5} \ll y_{8} \ll y_{11} \ll y_{13} \ll y_{14} \ll$
\end{flushleft}
\begin{flushright}
$\psi_{2}(y_{6}) \ll \psi_{2}(y_{7}) \ll \psi_{2}(y_{9}) \ll
\psi_{2}(y_{10}) \ll \psi_{2}(y_{12}) \ll \psi_{2}(y_{15})$.
\end{flushright}
Let ${\cal V}^{*}_{t} (t = 1, 2, 3)$ be the free monoid on ${\cal
V}_{t}$, and fix a positive integer $n$, with $n \geq 3$. Recall
that, for $\mu \in \{3,4,5,6\}$, $\nu \in \{1,2\}$,
$$
v^{(\mu,\nu,\Psi)}_{n,(\alpha,\beta,\gamma)} =
[\psi_{2}([x_{\mu},x_{\nu}]),x_{1,1}, \ldots,x_{\alpha,1},
x_{1,2}, \ldots,x_{\beta,2},x_{1,3}, \ldots,x_{\gamma,3}]
$$
with $x_{1,1}, \ldots, x_{\alpha,1} \in {\cal V}_{1},$ $x_{1,2},
\ldots,x_{\beta,2} \in {\cal V}_{2}$, $x_{1,3}, \ldots,
x_{\gamma,3} \in {\cal V}_{3}$, and for $\lambda \in \{5,6\}$,
$\tau \in \{3,4\}$,
$$
u^{(\lambda,\tau,\Psi)}_{n,(\delta,\epsilon)} =
[\psi_{2}([x_{\lambda},x_{\tau}]),z_{1,2},
\ldots,z_{\delta,2},z_{1,3}, \ldots,z_{\epsilon,3}]
$$
with $z_{1,2}, \ldots, z_{\delta,2} \in {\cal V}_{2},$ $z_{1,3},
\ldots, z_{\epsilon,3} \in {\cal V}_{3}$ and $\delta + \epsilon =
n-2$. For such $v^{(\mu,\nu,\Psi)}_{n,(\alpha,\beta,\gamma)},
u^{(\lambda,\tau,\Psi)}_{n,(\delta,\epsilon)}$, we write
$x^{(\mu,\nu)}_{n,\alpha} = x_{1,1} \cdots x_{\alpha,1} \in {\cal
V}^{*}_{1}$, $x^{(\mu,\nu)}_{n,\beta} = x_{1,2} \cdots
x_{\beta,2}, z^{(\lambda,\tau)}_{n,\delta} = z_{1,2} \cdots
z_{\delta,2} \in {\cal V}^{*}_{2}$ and $x^{(\mu,\nu)}_{n,\gamma} =
x_{1,3} \cdots x_{\gamma,3}, z^{(\lambda,\tau)}_{n,\epsilon} =
z_{1,3} \cdots z_{\epsilon,3} \in {\cal V}^{*}_{3}$. We write
$v^{(\mu_{1},\nu_{1},\Psi)}_{n,(\alpha_{1},\beta_{1},\gamma_{1})}
\ll
v^{(\mu_{2},\nu_{2},\Psi)}_{n,(\alpha_{2},\beta_{2},\gamma_{2})}$
if either $\psi_{2}([x_{\mu_{1}},x_{\nu_{1}}]) \ll
\psi_{2}([x_{\mu_{2}},x_{\nu_{2}}])$ or, if
$\psi_{2}([x_{\mu_{1}},x_{\nu_{1}}]) =
\psi_{2}([x_{\mu_{2}},x_{\nu_{2}}])$, and
$x^{(\mu_{1},\nu_{1})}_{n,\alpha_{1}} <
x^{(\mu_{2},\nu_{2})}_{n,\alpha_{2}}$ or, if
$\psi_{2}([x_{\mu_{1}},x_{\nu_{1}}]) =
\psi_{2}([x_{\mu_{2}},x_{\nu_{2}}])$,
$x^{(\mu_{1},\nu_{1})}_{n,\alpha_{1}} =
x^{(\mu_{2},\nu_{2})}_{n,\alpha_{2}}$, and
$x^{(\mu_{1},\nu_{1})}_{n,\beta_{1}} <
x^{(\mu_{2},\nu_{2})}_{n,\beta_{2}}$ or, if
$\psi_{2}([x_{\mu_{1}},x_{\nu_{1}}]) =
\psi_{2}([x_{\mu_{2}},x_{\nu_{2}}])$,
$x^{(\mu_{1},\nu_{1})}_{n,\alpha_{1}} =
x^{(\mu_{2},\nu_{2})}_{n,\alpha_{2}}$,
$x^{(\mu_{1},\nu_{1})}_{n,\beta_{1}} =
x^{(\mu_{2},\nu_{2})}_{n,\beta_{2}}$, and
$x^{(\mu_{1},\nu_{1})}_{n,\gamma_{1}} <
x^{(\mu_{2},\nu_{2})}_{n,\gamma_{2}}$. Similarly, we define
$u^{(\lambda_{1},\tau_{1},\Psi)}_{n,(\delta_{1},\epsilon_{1})} <
u^{(\lambda_{2},\tau_{2},\Psi)}_{n,(\delta_{2},\epsilon_{2})}$.
Thus the elements of ${\cal W}_{n,\Psi}$ are totally ordered.
Extend this ordering to ${\cal W}_{\Psi}$ by setting $u \ll v$ for
all $u \in {\cal W}_{n,\Psi}$ and $v \in {\cal W}_{m,\Psi}$ with
$n < m$.

\section{The ideal $J$}

\subsection{Preliminaries}

Recall that $J$ is the ideal of $L$ generated by the set $\cal V$,
where ${\cal V} = \{y_{1}, y_{2}, y_{3}, y_{6}+y_{7}, y_{7}+y_{5},
y_{9}+y_{8}, y_{10}+y_{8}, y_{12}+y_{13}, y_{15}+y_{13}\}$. For a
non-negative integer $m$, we write $[L^{2}, ~_{m}L]$ for the
$\mathbb{Z}$-module spanned by Lie commutators of the form $[v,
u_{1}, \ldots, u_{m}]$ with $u_{1}, \ldots, u_{m} \in L$ and $v
\in L^{2}$. By convention, $L^{2} = [L^{2},~_{0}L]$. Direct
calculations show that ${\cal V}$ is a $\mathbb{Z}$-basis of
$L^{2} \cap J$. Namely, we have the following $\mathbb{Z}$-module
decomposition of $L^{2}$
$$
L^{2} = (L^{2} \cap J) \oplus (L^{2})^{*},
$$
where $(L^{2})^{*}$ is the $\mathbb{Z}$-submodule of $L^{2}$
spanned by the set ${\cal V}^{*} = \{y_{4}, y_{5}, y_{8}, y_{11},
y_{13}, y_{14}\}$. Note that $(L^{2})^{*} = W^{(1)}_{2,\Psi}$. The
proof of the following result is elementary.

\vskip .120 in

\begin{lemma}\label{le14}
Let $\widetilde{J} = \sum_{m \geq 0} [L^{2} \cap J, ~_{m}L]$. Then
$\widetilde{J} = J$.
\end{lemma}

\vskip .120 in

Another useful technical result is the following.

\begin{lemma}\label{le15}
For each $c \geq 0$, let
$$
J_{c} = \sum_{m \geq c} [L^{2} \cap J, ~_{m}L].
$$
Then $J_{c}$ is an ideal of $L$ for all $c$,  $J_{c} = [L^{2} \cap
J, ~_{c}L^{1}] \oplus J_{c+1}$ and $[L^{2} \cap J, ~_{c}L^{1}] =
L^{c+2} \cap J = J^{c+2}$. Furthermore,
$$
J_{c} = \bigoplus_{m \geq c}([L^{2} \cap J, ~_{m}L^{1}]).
$$
\end{lemma}

\pf It is clear that $J_{c}$ is an ideal of $L$ for all $c$ and
$J_{0} = J$. Since $L^{2} \cap J \subseteq L^{2}$ and $J_{1}
\subseteq \gamma_{3}(L)$, we have $J = (L^{2} \cap J) \oplus
J_{1}$. It is easily verified that $J$ is the Lie subalgebra of
$L$ generated by the Lie commutators of the form  $h = [y, a_{1},
\ldots, a_{\kappa}]$ where $y \in L^{2} \cap J$ and $a_{1},
\ldots, a_{\kappa} \in L$. Since $L = \oplus_{m \geq 1}L^{m}$,
each $L^{m}$ is a $\mathbb{Z}$-module spanned by the Lie
commutators $[x_{j_{1}}, \ldots, x_{j_{m}}]$ with $j_{1}, \ldots,
j_{m} \in \{1, \ldots, 6\}$, and the Lie commutators are
multi-linear operations, we may assume that each $a_{\mu}$ is a
simple Lie commutator of the form $[x_{j_{1}}, \ldots,
x_{j_{m(\mu)}}]$ with $j_{1}, \ldots, j_{m(\mu)} \in \{1, \ldots,
6\}$ and $\mu \geq 1$. Thus $J$ is generated as a Lie subalgebra
by the set
\begin{flushleft}
$\{[y, [x_{j_{1,1}}, \ldots, x_{j_{m(1),1}}], \ldots,
[x_{j_{1,\kappa}}, \ldots, x_{j_{m(\kappa),\kappa}}]]: m(1) +
\cdots + m(\kappa) \geq 0,$ \end{flushleft}
\begin{flushright}
$[x_{j_{1,\lambda}}, \ldots, x_{j_{m(\lambda),\lambda}}] \in
L^{m(\lambda)}, x_{j_{1,\lambda}}, \ldots,
x_{j_{m(\lambda),\lambda}} \in \{x_{1}, \ldots, x_{6}\}, \lambda =
1, \ldots, \kappa\}.$
\end{flushright}
For $\nu \in \{1, \ldots, \kappa\}$, let $u_{\nu} =
[x_{j_{1,\nu}}, \ldots, x_{j_{m(\nu),\nu}}]$. Using the Jacobi
identity in the form $[x, [y, z]] = [x, y, z] - [x, z, y]$, we may
write each $[y, u_{1}, \ldots, u_{\kappa}]$, with $m(1) + \cdots +
m(\kappa) \geq c$, as a $\mathbb{Z}$-linear combination of Lie
commutators of the form
$$[y, x_{j_{1,1}}, \ldots, x_{j_{m(1),1}}, \ldots,
x_{j_{1,\kappa}}, \ldots, x_{j_{m(\kappa),\kappa}}]
$$ with $m(1) +
\cdots + m(\kappa) \geq c$. Let ${\cal J}_{\geq c} = \{[y,
x_{j_{1,1}}, \ldots, x_{j_{m(1),1}}, \ldots, x_{j_{1,\kappa}},
\ldots, x_{j_{m(\kappa),\kappa}}]: m(1) + \cdots + m(\kappa) \geq
c, j_{1,\lambda}, \ldots, j_{m(\lambda),\lambda} \in \{1, \ldots,
6\}, \lambda = 1, \ldots, \kappa\}$. It is clearly enough that the
Lie subalgebra of $L$ generated by ${\cal J}_{\geq c}$ is equal to
$J_{c}$. Let
\begin{flushleft}
${\cal J}_{ \geq c}^{c} = \{[y, x_{j_{1,1}}, \ldots,
x_{j_{m(1),1}}, \ldots, x_{j_{1,\kappa}}, \ldots,
x_{j_{m(\kappa),\kappa}}] \in {\cal J}_{\geq c}: m(1) + \cdots +
m(\kappa) = c,$
\end{flushleft}
\begin{flushright}
$j_{1,\lambda}, \ldots, j_{m(\lambda),\lambda} \in \{1, \ldots,
6\}, \lambda = 1, \ldots, \kappa\}$
\end{flushright}
and let $J^{c}_{\geq c}$ be the $\mathbb{Z}$-module spanned by
${\cal J}_{\geq c}^{c}$. Thus $J_{\geq c}^{c} = [L^{2} \cap J,
~_{c}L^{1}]$. Since $[L^{2} \cap J, ~_{c}L^{1}] \subseteq L^{c+2}$
and $J_{c+1} \subseteq \gamma_{c+3}(L)$, we have
$$
J_{c} = [L^{2} \cap J, ~_{c}L^{1}] \oplus J_{c+1}
$$
and so,
$$
J = \bigoplus_{c \geq 0} [L^{2} \cap J, ~_{c}L^{1}].
$$
On the other hand, since $J$ is generated by a set of homogenous
elements, we have
$$
J = \bigoplus_{m \geq 2} (L^{m} \cap J).
$$
Therefore, for all $c \geq 0$, $[L^{2} \cap J, ~_{c}L^{1}] = J
\cap L^{c+2} = J^{c+2}$. \qed

\subsection{A decomposition of $J$}

By Lemma \ref{l10}, Lemma
\ref{le4} and Corollary \ref{c2}, we have
$$
\begin{array}{lll}
L(W_{\Psi}) & = & L(W^{(1)}_{\Psi}) \oplus L(W^{(2)}_{\Psi})
\oplus L(\widetilde{W}_{\Psi,J}),
\end{array}
$$
where $\widetilde{W}_{\Psi,J} = \bigoplus_{n \geq
2}\widetilde{W}_{n,\Psi,J}$ and, for $n \geq 2$,
$$
\widetilde{W}_{n,\Psi,J} = \bigoplus_{\alpha + \beta = n-1 \atop
\alpha, \beta \geq 0}[W^{(2)}_{\Psi}, ~_{\alpha}W^{(1)}_{\Psi},
~_{\beta}W^{(2)}_{\Psi}].
$$
Note that $L(\widetilde{W}_{\Psi,J})$ is an ideal in
$L(W_{\Psi})$. Thus we have the following ${\mathbb{Z}}$-module
decomposition of $L(W_{\Psi})$ by means of the ideal $J$.

\begin{lemma}\label{11a}
The free Lie algebra $L(W_{\Psi})$ decomposes (as
${\mathbb{Z}}$-module) into a direct sum of the free Lie algebras
$L(W^{(1)}_{\Psi})$ and $L(W^{(2)}_{\Psi} \wr W^{(1)}_{\Psi}) =
L(W^{(2)}_{\Psi}) \oplus L(\widetilde{W}_{\Psi,J})$. In
particular, $L(W^{(2)}_{\Psi} \wr W^{(1)}_{\Psi}) \subseteq J$.
\end{lemma}

By the proof of Lemma \ref{l9} (Step 2),
$$
L^{\prime} = L^{\prime}(V_{1}) \oplus L^{\prime}(V_{2}) \oplus
L^{\prime}(V_{3}) \oplus L(W_{\Psi}).
$$
By Lemma \ref{11a} and the modular law, we have
$$
J = (L(W^{(1)}_{\Psi}) \cap J) \oplus
L^{\prime}(V_{1}) \oplus L^{\prime}(V_{2}) \oplus L^{\prime}(V_{3}) \oplus
L(W^{(2)}_{\Psi}) \oplus L(\widetilde{W}_{\Psi,J}). \eqno(2)
$$
Put
$$
J_{C} = L^{\prime}(V_{1}) \oplus L^{\prime}(V_{2}) \oplus L^{\prime}(V_{3}) \oplus
L(W^{(2)}_{\Psi}) \oplus L(\widetilde{W}_{\Psi,J}). \eqno (3)
$$
By Lemma \ref{11a},
$$
L^{\prime} = L(W^{(1)}_{\Psi}) \oplus
J_{C}. \eqno(4)
$$
Note that $J_{C} \subseteq J$, but $J_{C}$ is \emph{not} a Lie
subalgebra of $J$. (For example, $[y_{1}, y_{2}] \notin J_{C}$.)
By the equation (2), we have the following $\mathbb{Z}$-module
decomposition of $J$.

\begin{lemma}\label{le16}
$J = (L(W^{(1)}_{\Psi}) \cap J) \oplus
J_{C}$.
\end{lemma}

\section{$L/J$ is torsion free}

In this section we show that $L/J$ is torsion-free
$\mathbb{Z}$-module. For this we need to study the homogeneous
components of $J$ by means of Lyndon words, a filtration of the
tensor powers and the symmetric powers. For the necessary material
about Lyndon words and filtrations we refer the reader to
(\cite{loth}, \cite{reut}) and \cite{bsc}, respectively.

\subsection{Lyndon words}

In this subsection, we need some preliminaries results on Lyndon
words. Let $A$ be a totally ordered alphabet (not necessarily
finite). We order the free semigroup $A^{+}$ with alphabetical
order. By definition a word $w \in A^{+}$ is a Lyndon word if for
each non-trivial factorization $w = uv$ with $u, v \in A^{+}$, one has $w <
v$. The set of Lyndon words will be denoted by ${\mathbb{L}}_{A}$
(or briefly ${\mathbb{L}}$). If $v$ is the proper right factor of
maximal length of $w = uv \in {\mathbb{L}}$, then $u \in {\mathbb{L}}$, and $u < w < v$. Thus we
have a recursive algorithm to construct Lyndon words. The standard
factorization of each word $w$ of length $\geq 2$ is the
factorization $w = u \cdot v$, where $u$ is the smallest (proper)
right factor of $w$ for the alphabetical order. If $w$ is a Lyndon
word with standard factorization $u \cdot v$, then $u, v$ are
Lyndon words with $u < v$, $w < v$, and either $u$ is a letter
(i.e., an element of $A$), or the standard factorization of $u$ is
$x \cdot y$ with $y \geq v$.

Let $T(A)$ be the tensor algebra on
the free $\mathbb{Z}$-module with basis $A$. Let $L(A)$
denote the free Lie algebra (over $\mathbb{Z}$) on $A$. By the
Poincare-Birkhoff-Witt theorem, $T(A)$ is the universal enveloping
algebra of $L(A)$ (see, for example, \cite{jac}). Namely, giving
$T(A)$ the structure of Lie algebra, we may regard $L(A)$ as a Lie
subalgebra of $T(A)$. We let $q$ be the mapping of $\mathbb{L}$
into $L(A)$ defined inductively by $q(a) = a$, $a \in A$, and,
for $w \in {\mathbb{L}}\setminus \{A\}$, $q(w) =
[q(u),q(v)]$ where $w = u \cdot v$ is the standard
factorization of $w$. For a positive integer $m$, we write $A^{m}$
for the subset of $A^{+}$ consisting of all words of length $m$.
The following result is well known (see, for example,
\cite[Proposition 5.1.4, Lemma 5.3.2]{loth}).

\begin{lemma}\label{12c}
Let $\mathbb{L}$ be the set of Lyndon words on an alphabet $A$. If
$u, v \in {\mathbb{L}}$ with $u < v$, then $uv \in {\mathbb{L}}$.
Let $w \in {\mathbb{L}}\setminus A$ and its standard factorization
is $u \cdot v$. Then for any $y \in {\mathbb{L}}$ such that $w <
y$, the standard factorization of $wy \in {\mathbb{L}}$ is $w
\cdot y$ if and only if $y \leq v$. Let $w \in {\mathbb{L}}^{m} =
{\mathbb{L}} \cap A^{m}$ with $m \geq 2$. Then $q(w) = w + v$,
where $v$ belongs to the ${\mathbb{Z}}$-submodule of $T(A)$
spanned by those words $\widetilde{v} \in A^{m}$ such that $w <
\widetilde{v}$.
\end{lemma}

By Lemma \ref{12c}, $q({\mathbb{L}})$ is a set of linearly
independent elements and so, $q$ is injective. Furthermore, the
${\mathbb{Z}}$-module $L(A)$ is free with $q({\mathbb{L}})$ as
a basis. The elements of $q({\mathbb{L}})$ are called Lyndon
polynomials. We point out that the elements of
$q({\mathbb{L}})$ are simple Lie commutators. The proof of the
following result is straightforward.

\begin{corollary}\label{co3a}
For $i = 1, \ldots, m$, let $w_{i} \in {\mathbb{L}}_{A} \cap
A^{n_{i}}$. Then $q(w_{1}) \cdots q(w_{m})$ is equal to
$w_{1} \cdots w_{m} + v$, where $v$ belongs to the
${\mathbb{Z}}$-submodule of $T(A)$ spanned by those words
$\widetilde{v} \in A^{n_{1} + \cdots +n_{m}}$ such that $w_{1}
\cdots w_{m} < \widetilde{v}$.
\end{corollary}

\vskip .120 in

For $i = 1, 2, 3$, we recall that ${\cal V}_{i} = \{x_{2i-1},
x_{2i}\}$ with $x_{2i-1} < x_{2i}$, and let ${\mathbb{L}}_{V_{i}}$
denote the set of Lyndon words over ${\cal V}_{i}$. Let $T({\cal
V}_{i}) = T(V_{i})$ be the tensor algebra on $V_{i}$, and consider
$L(V_{i})$ as a Lie subalgebra of $T(V_{i})$. Furthermore, we
write $q({\mathbb{L}}_{V_{i}})$ for the $\mathbb{Z}$-basis of
$L(V_{i})$ corresponding to ${\mathbb{L}}_{V_{i}}$. For a positive
integer $m$, we write $q^{m}({\mathbb{L}}_{V_{i}})$ for
$L^{m}(V_{i}) \cap q({\mathbb{L}}_{V_{i}})$ which is a
$\mathbb{Z}$-basis for $L^{m}(V_{i})$. We recall that ${\cal
W}^{(\kappa)}_{\Psi} = \bigcup_{c \geq 0}{\cal
W}^{(\kappa)}_{c+2,\Psi}$ is a free generating set of
$L(W^{(\kappa)}_{\Psi})$ (with $\kappa = 1, 2$). We arbitrarily
order the elements of ${\cal W}^{(\kappa)}_{c+2,\Psi}$ ($\kappa =
1,2$) for all $c \geq 0$, and extend it to ${\cal
W}^{(\kappa)}_{\Psi}$ subject to $u < v$ if $u \in {\cal
W}^{(\kappa)}_{c+2,\Psi}$ and $v \in {\cal
W}^{(\kappa)}_{e+2,\Psi}$ with $c < e$. Let $T({\cal
W}^{(\kappa)}_{\Psi}) = T(W^{(\kappa)}_{\Psi})$ be the tensor
algebra on $W^{(\kappa)}_{\Psi}$. Regard $L(W^{(\kappa)}_{\Psi})$
as a Lie subalgebra of $T(W^{(\kappa)}_{\Psi})$ and the elements
of ${\cal W}^{(\kappa)}_{\Psi}$ are considered of degree $1$.
Write ${\mathbb{L}}_{W^{(\kappa)}_{\Psi}}$ for the set of Lyndon
words over ${\cal W}^{(\kappa)}_{\Psi}$. Thus
$q({\mathbb{L}}_{W^{(\kappa)}_{\Psi}})$ is a $\mathbb{Z}$-basis
for $L(W^{(\kappa)}_{\Psi})$. For $c \geq 0$, we write
$q^{c+2}_{\rm grad}({\mathbb{L}}_{W^{(\kappa)}_{\Psi}}) =
L^{c+2}_{\rm grad}(W^{(\kappa)}_{\Psi}) \cap
q({\mathbb{L}}_{W^{(\kappa)}_{\Psi}})$.

\subsection{Tensor and Symmetric Powers}

Let $T = T(L^{1}) =
T(V_{1} \oplus V_{2} \oplus V_{3})$ be the tensor algebra on
$L^{1}$. Note that $T(V_{i})$ ($i = 1, 2, 3)$ is a subalgebra of
$T(L^{1})$. For a nonnegative integer $c$, let $T^{c}$ denote the
$c$-th homogeneous component of $T$, that is, $T^{c}$ is the
$\mathbb{Z}$-submodule of $T$ spanned by all monomials $x_{i_{1}}
\cdots x_{i_{c}}$ with $i_{1}, \ldots, i_{c} \in \{1, \ldots,
6\}$. Thus $T = \bigoplus_{c \geq 0} T^{c}$ with $T^{0} =
{\mathbb{Z}}$. As before, we consider $L$ as a Lie subalgebra of
$T$. An analysis of the $c$-th homogeneous component $T^{c}$ of
$T$, with $c \geq 1$, will help us to understand better the
$(c+2)$-th homogeneous component $J^{c+2}$ of $J$. For $c \geq 1$,
we write ${\rm Part}(c)$ for the set of all partitions of $c$. By
a composition of $c$ we mean a sequence $\mu = (\mu_{1}, \ldots,
\mu_{\ell})$ of positive integers $\mu_{1}, \ldots, \mu_{\ell}$
satisfying $\mu_{1} + \cdots + \mu_{\ell} = c$. If $\mu =
(\mu_{1}, \ldots, \mu_{\ell})$ is a composition of $c$ and if we
rearrange $\mu_{1}, \ldots, \mu_{r}$ in such a way to obtain a
partition of $c$, we call it the associated partition to $\mu$.
(For example, $\mu = (1,2,1,2)$ is a composition of $6$ and
$(2,2,1,1)$ is the associated partition of $\mu$.) We use the
lexicographic order $\leq^{*}$ on the elements of ${\rm Part}(c)$.
Note that the smallest partition is $(1, \ldots, 1)$ written as
$(1^{c})$ and the largest partition is $(c)$.

For a free $\mathbb{Z}$-module $U$ of (finite) rank $r \geq 2$ with
a free generating set $\{u_{1}, \ldots, u_{r}\}$, we write $S(U)$
for the free symmetric algebra on $U$, that is, $S(U) =
{\mathbb{Z}}[u_{1}, \ldots, u_{r}]$. For a non-negative integer
$c$, we write $S^{c}(U)$ for the $c$-th homogeneous component of
$S(U)$ with $S^{0}(U) = \mathbb{Z}$. Thus $S(U) = \bigoplus_{c
\geq 0}S^{c}(U)$. The proof of the following result is elementary.

\begin{lemma}\label{12b}
Let $U_{1}, \ldots, U_{\kappa}$, with $\kappa \geq 2$, be free
$\mathbb{Z}$-modules of finite rank. Then, for $d \geq 1$
$$
S^{d}(U_{1} \oplus \cdots \oplus U_{\kappa}) \cong
\bigoplus_{\nu_{1} + \cdots + \nu_{\kappa} = d \atop \nu_{1},
\ldots, \nu_{\kappa \geq 0}} S^{\nu_{1}}(U_{1})
\otimes_{\mathbb{Z}} \cdots \otimes_{\mathbb{Z}}
S^{\nu_{\kappa}}(U_{\kappa})
$$
as $\mathbb{Z}$-modules in a natural way.
\end{lemma}

For a positive integer $m$, let ${\cal Y}_{m}$ be the set of all
elements of $T^{m}$ of the form $b_{1}b_{2} \cdots b_{\ell}$
($\ell$ arbitrary) where each $b_{i} \in L^{\mu_{i}}\setminus
\{0\}$ for some positive integer $\mu_{i}$. Thus $\mu_{i} = {\rm
deg}b_{i}$ for all $i$ and $(\mu_{1}, \ldots, \mu_{\ell})$ is a
composition of $m$. For each partition $\lambda$ of $m$ let ${\cal
Y}_{\lambda}$ denote the set of all such elements $b_{1}b_{2}
\cdots b_{\ell}$ where $({\rm deg}b_{1}, \ldots, {\rm
deg}b_{\ell})$ has $\lambda$ as its associated partition. For each
$\lambda \in {\rm Part}(m)$, let $\Phi_{\lambda}$ be the
$\mathbb{Z}$-module spanned by ${\cal Y}_{\theta}$ with $\lambda
\leq^{*} \theta$. For each $\lambda$ such that $\lambda \neq (m)$
let $\lambda+1$ be the partition of $m$ which is next bigger than
$\lambda$. Thus we have the filtration
$$
T^{m} = \Phi_{(1^{m})} \geq \cdots \geq \Phi_{\lambda} \geq
\Phi_{\lambda+1} \geq \cdots \geq \Phi_{(m)} \geq \{0\}.
$$
(Since $L$ is free Lie algebra of rank $6$, we have
$\Phi_{\lambda} > \Phi_{\lambda+1}$.) Thus, for all $\lambda$,
${\cal Y}_{\lambda}$ spans $\Phi_{\lambda}$ modulo
$\Phi_{\lambda+1}$. (We write $\Phi_{(m)+1} = \{0\}$.) Let $\mu =
(\mu_{1}, \mu_{2}, \ldots, \mu_{\ell})$ be a composition of $m$
with associated partition $\lambda$. For $i = 1, \ldots, \ell$,
let $b_{i} \in L^{\mu_{i}}\setminus \{0\}$, and let $\pi \in {\rm
Sym}(\ell)$. Thus $b_{1}b_{2} \cdots b_{\ell}$ and
$b_{\pi(1)}b_{\pi(2)} \cdots b_{\pi(\ell)}$ belong to ${\cal
Y}_{\lambda}$. As observed in \cite[Lemma 3.1]{bsc}, we have
$$
b_{1}b_{2} \cdots b_{\ell} + \Phi_{\lambda+1} =
b_{\pi(1)}b_{\pi(2)} \cdots b_{\pi(\ell)}+\Phi_{\lambda+1}. \eqno
(5)
$$

Recall that the set ${\cal V} = \{y_{1}, y_{2}, y_{3},$
$\psi_{2}(y_{6}), \psi_{2}(y_{7}), \psi_{2}(y_{9}), \psi_{2}(y_{10}), \psi_{2}(y_{12}), \psi_{2}(y_{15})\}$ is
a $\mathbb{Z}$-basis of $J^{2}$. Since $J^{c+2} = [J^{2},
~_{c}L^{1}]$ (by Lemma \ref{le15}) and the multi-linearity of the
Lie bracket, we have $J^{c+2}$ is the $\mathbb{Z}$-module spanned
by all Lie commutators of the form $[v, x_{i_{1}}, \ldots,
x_{i_{c}}]$ with $i_{1}, \ldots, i_{c} \in \{1, \ldots, 6\}$. For
$u = \sum *~x_{j_{1}} \cdots x_{j_{c}} \in T^{c}$ with the
coefficients $*$ are in $\mathbb{Z}$, we write $[v;u] = \sum
*~[v,x_{j_{1}}, \ldots, x_{j_{c}}]$. It is easily verified that
for $u, w \in T^{c}$ and $\alpha, \beta \in {\mathbb{Z}}$,
$[v;\alpha u + \beta w] = \alpha [v;u] + \beta [v;w]$. For each
partition $\lambda$ of $m$ let $[J^{2};\Phi_{\lambda}]$ be the
$\mathbb{Z}$-module spanned by all Lie commutators of the form
$[v, b_{1}, \ldots, b_{\ell}]$ where $v \in {\cal V}$, $b_{1}
\cdots b_{\ell} \in {\cal Y}_{\theta}$ with $\lambda \leq^{*}
\theta$. Note that $[J^{2};\Phi_{(1^{m})}] = [J^{2};T^{m}]$, and
$$
J^{m+2} = [J^{2};T^{m}] \geq [J^{2};\Phi_{(2,1^{m-1})}] \geq
\cdots \geq [J^{2};\Phi_{(m)}] > \{0\}.
$$
The following result can be proved using the Jacobi identity  in
the form $[x,y,z] = [x,z,y] + [x,[y,z]]$.

\begin{lemma}\label{12a}
Let $\mu = (\mu_{1}, \ldots, \mu_{\ell})$ be a composition of $m$
with associated partition $\lambda$. For $i = 1, \ldots, \ell$,
let $b_{i} \in L^{\mu_{i}} \setminus \{0\}$, and let $\pi \in {\rm
Sym}(\ell)$. Then, for $v \in J^{2}$, $[v, b_{1}, \ldots,
b_{\ell}]$ and $[v, b_{\pi(1)}, \ldots, b_{\pi(\ell)}]$ belong to
$[J^{2};\Phi_{\lambda}]$, and
$$
[v, b_{1}, \ldots, b_{\ell}] + [J^{2};\Phi_{\lambda+1}] = [v,
b_{\pi(1)}, \ldots, b_{\pi(\ell)}] + [J^{2};\Phi_{\lambda+1}].
$$
\end{lemma}

For $\kappa = 1, 2, 3$, we similarly define $J^{m+2}_{\kappa} =
[L^{2}(V_{\kappa});T^{m}]$ and $J^{m+2}_{\Psi} =
[W^{(2)}_{2,\Psi};T^{m}]$. Since $J^{2} = L^{2}(V_{1}) \oplus
L^{2}(V_{2}) \oplus L^{3}(V_{3}) \oplus W^{(2)}_{2,\Psi}$, we have
$$
J^{m+2} = J^{m+2}_{1} + J^{m+2}_{2} + J^{m+2}_{3} + J^{m+2}_{\Psi},
$$
and, by Lemma \ref{le15},
$$
J = J^{2} \oplus (\bigoplus_{m \geq 1}(J^{m+2}_{1} + J^{m+2}_{2} + J^{m+2}_{3} +
J^{m+2}_{\Psi})).
$$

\begin{corollary}\label{co3}
Let $\mu = (\mu_{1}, \ldots, \mu_{\ell})$ be a composition of $m$
with associated partition $\lambda$. For $i = 1, \ldots, \ell$,
let $b_{i} \in L^{\mu_{i}} \setminus \{0\}$, and let $\pi \in {\rm
Sym}(\ell)$. Then, for $j = 1, 2, 3$, $[y_{j}, b_{1}, \ldots,
b_{\ell}]$ and $[y_{j}, b_{\pi(1)}, \ldots, b_{\pi(\ell)}]$ belong
to $[L^{2}(V_{j});\Phi_{\lambda}]$, and
$$
[y_{j}, b_{1}, \ldots, b_{\ell}] + [L^{2}(V_{j});\Phi_{\lambda+1}]
= [y_{j}, b_{\pi(1)}, \ldots, b_{\pi(\ell)}] +
[L^{2}(V_{j});\Phi_{\lambda+1}].
$$
Furthermore, for $u \in {\cal W}^{(2)}_{2,\Psi}$, $[u, b_{1},
\ldots, b_{\ell}]$ and $[u, b_{\pi(1)}, \ldots, b_{\pi(\ell)}]$
belong to $[W^{(2)}_{2,\Psi};\Phi_{\lambda}]$, and
$$
[u, b_{1}, \ldots, b_{\ell}] + [W^{(2)}_{2,\Psi};\Phi_{\lambda+1}]
= [u, b_{\pi(1)}, \ldots, b_{\pi(\ell)}] +
[W^{(2)}_{2,\Psi};\Phi_{\lambda+1}].
$$
\end{corollary}

For a positive integer $n$, with $n \geq 2$, we write ${\cal
\widetilde{W}}_{n,\Psi,J}$ for the natural $\mathbb{Z}$-basis of
$\widetilde{W}_{n,\Psi,J}$, and ${\cal \widetilde{W}}_{\Psi,J} =
\bigcup_{n \geq 2}{\cal \widetilde{W}}_{n,\Psi,J}$. We arbitrarily
order the elements of ${\cal \widetilde{W}}_{n,\Psi,J}$ for all $n
\geq 2$, and extend it to ${\cal \widetilde{W}}_{\Psi,J}$ subject
to $u < v$ if $u \in {\cal \widetilde{W}}_{n,\Psi,J}$ and $v \in
{\cal \widetilde{W}}_{m,\Psi,J}$ with $n < m$. By Lemma \ref{le6}
(for ${\cal A} = {\cal \widetilde{W}}_{\Psi,J}$), $L({\cal
\widetilde{W}}_{\Psi,J})^{\prime}$ is freely generated by the set
${\cal \widetilde{W}}^{(1)}_{\Psi,J}$
$$
{\cal \widetilde{W}}^{(1)}_{\Psi,J} = \{[a_{i_{1}}, \ldots,
a_{i_{k}}]: k \geq 2, a_{i_{1}} > a_{i_{2}} \leq a_{i_{3}} \leq
\cdots \leq a_{i_{k}}, a_{i_{1}}, \ldots, a_{i_{k}} \in {\cal
\widetilde{W}}_{\Psi,J}\}.
$$
We arbitrarily order the elements of ${\cal
\widetilde{W}}^{(1)}_{\Psi,J}$ of the same degree, and elements of
degree $r$ are strictly less than the elements of degree $s$ with
$r < s$. Furthermore, for a positive integer $e$, let
$L^{(e)}({\cal \widetilde{W}}_{\Psi,J}) = (L^{(e-1)}({\cal
\widetilde{W}}_{\Psi,J}))^{\prime}$ with $L^{(0)}({\cal
\widetilde{W}}_{\Psi,J}) = L({\cal \widetilde{W}}_{\Psi,J})$ and
$L^{(1)}({\cal \widetilde{W}}_{\Psi,J}) = L({\cal
\widetilde{W}}_{\Psi})^{\prime}$. If ${\cal
\widetilde{W}}^{(e)}_{\Psi,J}$ is an ordered free generating set
for $L^{(e)}({\cal \widetilde{W}}_{\Psi,J})$, then ${\cal
\widetilde{W}}^{(e+1)}_{\Psi,J}$ is a free generating set for
$L^{(e+1)}({\cal \widetilde{W}}_{\Psi,J})$ with
$$
{\cal \widetilde{W}}^{(e+1)}_{\Psi,J} = \{[a^{(e)}_{i_{1}},
\ldots, a^{(e)}_{i_{k}}]: k \geq 2, a^{(e)}_{i_{1}} >
a^{(e)}_{i_{2}} \leq a^{(e)}_{i_{3}} \leq \cdots \leq
a^{(e)}_{i_{k}}, a^{(e)}_{i_{1}}, \ldots, a^{(e)}_{i_{k}} \in
{\cal \widetilde{W}}^{(e)}_{\Psi,J}\}.
$$
We write
$$
{\cal \widetilde{W}}^{\rm ext}_{\Psi,J} = \bigcup_{d \geq 0}{\cal
\widetilde{W}}^{(e)}_{\Psi,J}.
$$
Note that the set
$$
{\cal W}^{\star} = (\bigcup_{i=1}^{3}q({\mathbb{L}}_{{\cal
V}_{i}})) \cup (\bigcup_{\kappa =
1}^{2}q({\mathbb{L}}_{W^{(\kappa)}_{\Psi}})) \cup {\cal
\widetilde{W}}^{\rm ext}_{\Psi,J}
$$
is a ${\mathbb{Z}}$-basis of $L$. For a non-negative integer $e$,
let ${\cal \widetilde{W}}^{{\rm ext},e}_{c+2,\Psi,J} = L^{c+2}
\cap {\cal \widetilde{W}}^{(e)}_{\Psi,J}$. That is, ${\cal
\widetilde{W}}^{{\rm ext},e}_{c+2,\Psi,J}$ consists of all Lie
commutators in ${\cal \widetilde{W}}^{(e)}_{\Psi,J}$ of total
degree $c+2$. There exists a (unique) positive integer $n(c+2)$
such that the set ${\cal \widetilde{W}}^{{\rm ext}}_{c+2,\Psi,J} =
\bigcup_{e=0}^{n(c+2)}{\cal \widetilde{W}}^{{\rm
ext},e}_{c+2,\Psi,J}$ is a $\mathbb{Z}$-basis of $L^{c+2}_{\rm
grad}(\widetilde{W}_{\Psi,J})$. It is clear that $L^{c+2}_{\rm
grad}(\widetilde{W}_{\Psi,J}) =
\bigoplus_{e=0}^{n(c+2)}\widetilde{W}^{{\rm ext},e}_{c+2,\Psi}$
where $\widetilde{W}^{{\rm ext},e}_{c+2,\Psi}$ is the
$\mathbb{Z}$-span of ${\cal \widetilde{W}}^{{\rm
ext},e}_{c+2,\Psi}$. Thus, for $c \geq 0$, the set
$$
{\cal W}^{c+2} = (\bigcup_{i=1}^{3}q^{c+2}({\mathbb{L}}_{{\cal
V}_{i}})) \cup (\bigcup_{\kappa = 1}^{2}q_{\rm
grad}^{c+2}({\mathbb{L}}_{W^{(\kappa)}_{\Psi}})) \cup {\cal
W}^{\rm ext}_{c+2,\Psi,J}
$$
is a $\mathbb{Z}$-basis of $L^{c+2}$ and so,
$$
L^{c+2} = (\bigoplus_{i=1}^{3}L^{c+2}(V_{i})) \oplus
(\bigoplus_{\kappa = 1}^{2}L^{c+2}_{\rm
grad}(W^{(\kappa)}_{\Psi})) \oplus L^{c+2}_{\rm
grad}(\widetilde{W}_{\Psi,J}).
$$
The elements of ${\cal W}^{c+2}$ are ordered in such a way that
the elements in $q^{c+2}({\mathbb{L}}_{{\cal V}_{i}})$ ($i = 1,
2, 3)$, the elements in $q_{\rm
grad}^{c+2}({\mathbb{L}}_{W^{(\kappa)}_{\Psi}})$ ($\kappa = 1,2)$
and the elements in ${\cal W}^{\rm ext}_{c+2,\Psi,J}$ are preserve
their orderings, and $u_{1} \prec u_{2} \prec \cdots \prec u_{6}$
for all $u_{i} \in q^{c+2}({\mathbb{L}}_{{\cal V}_{i}}))$ ($i =
1, 2, 3)$, $u_{t} \in q_{\rm
grad}^{c+2}({\mathbb{L}}_{W^{(t-3)}_{\Psi}}))$ ($t = 4,5)$ and
$u_{6} \in {\cal W}^{\rm ext}_{c+2,\Psi,J}$. Choose a total
ordering $\preceq$ of ${\cal W}^{\star} = {\cal X} \cup
(\bigcup_{c \geq 0} {\cal W}^{c+2})$ in such a way all elements of
${\cal W}^{i}$ are smaller that all elements of ${\cal W}^{j}$
whenever $i < j$ with ${\cal W}^{1} = {\cal X}$. (Note that the
elements of each set ${\cal W}^{i}$ are already ordered.) By the
Poincare-Birkhoff-Witt Theorem (see \cite{bour}), $T$ has a
$\mathbb{Z}$-basis $\cal T$ consisting of all elements of the form
$$
a_{1}a_{2} \cdots a_{\ell}~~~~(\ell \geq 0, a_{1}, \ldots,
a_{\ell} \in {\cal W}^{\star}, a_{1} \preceq \cdots \preceq
a_{\ell}).
$$
The above elements of $\cal T$ are distinct as written. For $c
\geq 1$, we write ${\cal T}_{c}$ for the set of the above elements
of degree $c$. Clearly, ${\cal T}_{c}$ is a $\mathbb{Z}$-basis of
$T^{c}$. Also, for $\lambda \in {\rm Part}(c)$, we write ${\cal
T}_{\lambda}$ for the set of elements of the above form in ${\cal
T}_{c}$ such that the composition $({\rm deg}a_{1}, \ldots, {\rm
deg}a_{\ell})$ has $\lambda$ as associated partition. As pointed
out in \cite[at the bottom of p. 183]{bsc} the elements of ${\cal
T}_{\lambda}$ taken modulo $\Phi_{\lambda+1}$ form a
$\mathbb{Z}$-basis of $\Phi_{\lambda}/\Phi_{\lambda+1}$. Any
partition $\lambda \in {\rm Part}(c)$ is written as $\lambda =
(c^{n(c)}, (c-1)^{n(c-1)}, \ldots, 1^{n(1)})$ with $n(c), \ldots,
n(1)$ non-negative integers. For such a partition $\lambda$, we
define
$$
L^{\lambda} = S^{n(1)}(L^{1}) \otimes_{\mathbb{Z}} S^{n(2)}(L^{2})
\otimes_{\mathbb{Z}} \cdots \otimes_{\mathbb{Z}} S^{n(c)}(L^{c}).
$$
As observed in \cite[pp. 183-184]{bsc},
$\Phi_{\lambda}/\Phi_{\lambda + 1} \cong L^{\lambda}$ as
$\mathbb{Z}$-modules, and the basis $\{w+\Phi_{\lambda+1}: w \in
{\cal T}_{\lambda}\}$ of $\Phi_{\lambda}/\Phi_{\lambda+1}$
consists of all elements of the form
$$
a^{(1)}_{1} \cdots a^{(n(1))}_{1} a^{(1)}_{2} \cdots
a^{(n(2))}_{2} \cdots a^{(1)}_{c} \cdots a^{(n(c))}_{c} +
\Phi_{\lambda+1} \eqno (6)
$$
where $a^{(1)}_{1} \preceq \cdots \preceq a^{(n(1))}_{1} \preceq
a^{(1)}_{2} \preceq \cdots  \preceq a^{(n(2))}_{2} \preceq \cdots
\preceq a^{(1)}_{c} \preceq \cdots  \preceq a^{(n(c))}_{c}$ and
$a^{(1)}_{i}, \cdots, a^{(n(i))}_{i} \in {\cal W}^{i}$ for $i = 1,
\ldots, c$. (For $i = 1$, ${\cal W}^{1} = {\cal X}$.)

\subsection{An analysis of $L^{\lambda}$}

We define ${\cal G}_{\lambda} = {\cal G}_{1,n(1)} \times \cdots
\times {\cal G}_{c,n(c)}$ and ${\cal B}_{\lambda} = {\cal
B}_{1,n(1)} \times \cdots \times {\cal B}_{c,n(c)}$, where
$$
{\cal G}_{1,n(1)} = \{(\nu_{1,1},\nu_{2,1},\nu_{3,1}) \in
{\mathbb{N}}_{0}^{3}: \nu_{1,1}+\nu_{2,1}+\nu_{3,1} = n(1)\},
$$
for $\kappa = 2, 3$,
$$
{\cal G}_{\kappa,n(\kappa)} = \{(\nu_{1,\kappa},
\ldots,\nu_{5,\kappa}) \in {\mathbb{N}}_{0}^{5}: \nu_{1,\kappa}+
\cdots + \nu_{5,\kappa} = n(\kappa)\},
$$
and, for $4 \leq \kappa \leq c$,
$$
{\cal G}_{\kappa,n(\kappa)} = \{(\nu_{1,\kappa},
\ldots,\nu_{6,\kappa}) \in {\mathbb{N}}_{0}^{6}: \nu_{1,\kappa}+
\cdots + \nu_{6,\kappa} = n(\kappa)\},
$$
(${\mathbb{N}}_{0}$ denotes the set of all non-negative integers.)
$$
{\cal B}_{1,n(1)} =
\{(S^{\nu_{1,1}}(V_{1}),S^{\nu_{2,1}}(V_{2}),S^{\nu_{3,1}}(V_{3})):
(\nu_{1,1},\nu_{2,1},\nu_{3,1}) \in {\cal G}_{1,n(1)}\},
$$
for $\kappa = 2,3$,
$$
{\cal B}_{\kappa,n(\kappa)} =
\{(S^{\nu_{1,\kappa}}(L^{\kappa}(V_{1})), \ldots,
S^{\nu_{5,\kappa}}(L^{\kappa}(W^{(2)}_{\Psi})): (\nu_{1,\kappa},
\ldots, \nu_{5,\kappa}) \in {\cal G}_{\kappa,n(\kappa)}\},
$$
and, for $4 \leq \kappa \leq c$,
$$
{\cal B}_{\kappa,n(\kappa)} =
\{(S^{\nu_{1,\kappa}}(L^{\kappa}(V_{1})), \ldots,
S^{\nu_{6,\kappa}}(L^{\kappa}_{\rm grad}(\widetilde{W}_{\Psi,J})):
(\nu_{1,\kappa}, \ldots, \nu_{6,\kappa}) \in {\cal
G}_{\kappa,n(\kappa)}\}.
$$
By applying Lemma \ref{12b} (using the above $\mathbb{Z}$-module
decompositions of $L^{1}, \ldots, L^{c}$, respectively), by the
associativity of $\otimes_{\mathbb{Z}}$ and by the distributivity
of $\oplus_{\mathbb{Z}}$ with $\otimes_{\mathbb{Z}}$, we have
$$
L^{\lambda} \cong \bigoplus_{(X_{1}, \ldots, X_{c}) \in {\cal
B}_{\lambda}}(\overline{X}_{1} \otimes_{\mathbb{Z}} \cdots
\otimes_{\mathbb{Z}} \overline{X}_{c}).
$$
(For example, if $X_{1} =
(S^{\nu_{1,1}}(V_{1}),S^{\nu_{2,1}}(V_{2}),S^{\nu_{3,1}}(V_{3}))
\in {\cal B}_{1,n(1)}, \ldots, X_{c} =
(S^{\nu_{1,c}}(L^{c}(V_{1})),$ \linebreak $\ldots,
S^{\nu_{6,c}}(L^{c}_{\rm grad}(\widetilde{W}_{\Psi,J})) \in {\cal
B}_{c,n(c)}$, by $\overline{X}_{1} \otimes_{\mathbb{Z}} \cdots
\otimes_{\mathbb{Z}} \overline{X}_{c}$, we mean the expression
$$
S^{\nu_{1,1}}(V_{1}) \otimes_{\mathbb{Z}} S^{\nu_{2,1}}(V_{2})
\otimes_{\mathbb{Z}} S^{\nu_{3,1}}(V_{3}) \otimes_{\mathbb{Z}}
\cdots \otimes_{\mathbb{Z}} S^{\nu_{1,c}}(L^{c}(V_{1}))
\otimes_{\mathbb{Z}} \cdots \otimes_{\mathbb{Z}}
S^{\nu_{6,c}}(L^{c}_{\rm grad}(\widetilde{W}_{\Psi,J})).
$$
The elements of ${\cal G}_{\kappa,n(\kappa)}$ are denoted by
$Y_{1,n(1),3}$ (for $\kappa = 1$), $Y_{\kappa,n(\kappa),5}$ (for
$\kappa = 2, 3$) and $Y_{\kappa,n(\kappa),6}$ (for $4 \leq \kappa
\leq c$). For $Y = (Y_{1,n(1),3}, \ldots, Y_{c,n(c),6}) \in {\cal
G}_{\lambda}$, we correspond a unique element of the basis of
$\Phi_{\lambda}/\Phi_{\lambda+1}$ of the form (6) where the first
$\nu_{1,1}$ elements are in ${\cal V}_{1}$ (if $\nu_{1,1} \geq
1$), the second $\nu_{2,1}$ elements are in ${\cal V}_{2}$ (if
$\nu_{2,1} \geq 1$), the next $\nu_{3,1}$ elements are in ${\cal
V}_{3}$ (if $\nu_{3,1} \geq 1$) with
$(\nu_{1,1},\nu_{2,1},\nu_{3,1}) \in Y_{1,n(1),3}$ and so on, and
vice versa. Furthermore, for such $Y \in {\cal G}_{\lambda}$, we
associate the nonnegative integers $m_{\lambda,Y}(1), \ldots,
m_{\lambda,Y}(6)$ defined as follows:
$$
\begin{array}{lll}
m_{\lambda,Y}(t) & = & \nu_{t,1}\cdot 1 + \nu_{t,2} \cdot 2 +
\cdots + \nu_{t,c} \cdot c ~~(t = 1, 2, 3), \\
m_{\lambda,Y}(t) & = & \nu_{t,2} \cdot 2 + \cdots + \nu_{t,c}
\cdot
c ~~(t = 4, 5) \\
{\rm and} & & \\
m_{\lambda,Y}(6) & = & \nu_{6,4} \cdot 4 + \cdots + \nu_{6,c}
\cdot c.
\end{array}
$$
Note that $m_{\lambda,Y}(1) + \cdots + m_{\lambda,Y}(6) = c$ for
all $Y \in {\cal G}_{\lambda}$. By the equation (5), we rearrange
the above elements in such a way the first $m_{\lambda,Y}(1)$
elements are in $q({\mathbb{L}}_{V_{1}})$, the second
$m_{\lambda,Y}(2)$ elements are in $q({\mathbb{L}}_{V_{2}})$
and so on. More precisely, for $Y = (Y_{1,n(1),3}, \ldots,
Y_{c,n(c),6}) \in {\cal G}_{\lambda}$ we have the partitions
$p_{m_{\lambda,Y}(t)}$ of $m_{\lambda,Y}(t)$, $t = 1, \ldots, 6$,
where $p_{m_{\lambda,Y}(t)} = (c^{\nu_{t,c}}, (c-1)^{\nu_{t,c-1}},
\ldots, 2^{\nu_{t,2}}, 1^{\nu_{t,1}})$ (for $t = 1, 2, 3)$,
$p_{m_{\lambda,Y}(t)} = (c^{\nu_{t,c}}, (c-1)^{\nu_{t,c-1}},
\ldots, 2^{\nu_{t,2}})$ (for $t = 4, 5)$ and $p_{m_{\lambda,Y}(6)}
= (c^{\nu_{6,c}}, (c-1)^{\nu_{6,c-1}}, \ldots, 4^{\nu_{6,c}})$. It
is clearly enough that for $\lambda = (c^{n(c)}, (c-1)^{n(c-1)},
\ldots, 1^{n(1)}) \in {\rm Part}(c)$ with $n(c), \ldots, n(1)$
nonnegative integers
$$
L^{\lambda} \cong \bigoplus_{Y \in {\cal
G}_{\lambda}}(L^{p_{m_{\lambda,Y}(1)}}(V_{1}) \otimes_{\mathbb{Z}}
\cdots \otimes_{\mathbb{Z}}
L^{p_{m_{\lambda,Y}(6)}}(\widetilde{W}_{\Psi,J})), \eqno (7)
$$
where $L^{p_{m_{\lambda,Y}(t)}}(V_{t}) = S^{\nu_{t,1}}(V_{t})
\otimes_{\mathbb{Z}} \cdots \otimes_{\mathbb{Z}}
S^{\nu_{t,c}}(L^{c}(V_{t}))$ $(t = 1, 2, 3)$,
$L^{p_{m_{\lambda,Y}(t)}}(W^{(t-3)}_{\Psi}) =
S^{\nu_{t,2}}(L^{2}_{\rm grad}(W^{(t-3)}_{\Psi}))
\otimes_{\mathbb{Z}} \cdots \otimes_{\mathbb{Z}}
S^{\nu_{t,c}}(L^{c}_{\rm grad}(W^{(t-3)}_{\Psi}))$ $(t = 4, 5)$
and $L^{p_{m_{\lambda,Y}(6)}}(\widetilde{W}_{\Psi,J}) =$
\linebreak $S^{\nu_{6,4}}(L^{4}_{\rm
grad}(\widetilde{W}_{\Psi,J})) \otimes_{\mathbb{Z}} \cdots
\otimes_{\mathbb{Z}} S^{\nu_{6,c}}(L^{c}_{\rm
grad}(\widetilde{W}_{\Psi,J}))$.

\subsection{A~ $\mathbb{Z}$-basis of
$\Phi_{\lambda}/\Phi_{\lambda+1}$}

For $t = 1, 2, 3$, we write
$$
w(p_{m_{\lambda,Y}(t)}) = a_{1,t}^{(1)} \cdots
a_{1,t}^{(\nu_{t,1})}a_{2,t}^{(1)} \cdots a_{2,t}^{(\nu_{t,2})}
\cdots a_{c,t}^{(1)} \cdots a_{c,t}^{(\nu_{t,c})} \in T(L^{1}),
$$
where $a_{j,t}^{(1)}, \ldots, a_{j,t}^{(\nu_{t,j})},
 \in
q^{j}({\mathbb{L}}_{V_{t}})$, $j = 1, \ldots, c$, for $t = 4,
5$
$$
w(p_{m_{\lambda,Y}(t)}) = a_{2,t}^{(1)} \cdots
a_{2,t}^{(\nu_{t,2})} \cdots a_{c,t}^{(1)} \cdots
a_{c,t}^{(\nu_{t,c})} \in T(L^{1}),
$$
where $a_{j,t}^{(1)}, \ldots, a_{j,t}^{(\nu_{t,j})} \in
q^{j}_{\rm grad}({\mathbb{L}}_{W^{(t-3)}_{\Psi}})$, $j = 2,
\ldots, c$, and
$$
w(p_{m_{\lambda,Y}(6)}) = a_{4,6}^{(1)} \cdots
a_{4,6}^{(\nu_{6,4})} \cdots a_{c,6}^{(1)} \cdots
a_{c,6}^{(\nu_{6,c})} \in T(L^{1}),
$$
where $a_{j,6}^{(1)}, \ldots, a_{j,6}^{(\nu_{6,j})} \in {\cal
W}^{\rm ext}_{j,\Psi,J}$, $j = 4, \ldots, c$. Furthermore, we
write
$$
w(p_{m_{\lambda,Y}(1)}, \ldots, p_{m_{\lambda,Y}(6)}) =
w(p_{m_{\lambda,Y}(1)}) \cdots w(p_{m_{\lambda,Y}(6)}) \in
T(L^{1}).
$$

\begin{lemma}
The set $\{w(p_{m_{\lambda,Y}(1)}, \ldots, p_{m_{\lambda,Y}(6)}) +
\Phi_{\lambda+1}: Y \in {\cal G}_{\lambda}\}$ spans
$\Phi_{\lambda}/\Phi_{\lambda+1}$. In particular, it is a
$\mathbb{Z}$-basis of $\Phi_{\lambda}/\Phi_{\lambda+1}$.
\end{lemma}

\pf Indeed, writing $w = a^{(1)}_{1} \cdots a^{(n(1))}_{1} \cdots
a^{(1)}_{c} \cdots a^{(n(c))}_{c} = b_{1} \cdots b_{\ell}$, then,
by a suitable permutation $\pi \in {\rm Sym}(\ell)$, $b_{\pi(1)}
\cdots b_{\pi(\ell)} + \Phi_{\lambda+1} = w(p_{m_{\lambda,Y}(1)},
\ldots, p_{m_{\lambda,Y}(6)}) + \Phi_{\lambda+1}$ for a unique $Y
\in {\cal G}_{\lambda}$. Since the set $\{w+\Phi_{\lambda+1}: w
\in {\cal T}_{\lambda}\}$ is a $\mathbb{Z}$-basis, we have the
required result. \qed

\begin{remark}\label{r1}
We point out that, for any $\pi \in {\rm Sym}(6)$, we have
$w(p_{m_{\lambda,Y}(1)}, \ldots, p_{m_{\lambda,Y}(6)}) +
\Phi_{\lambda+1} = w(p_{m_{\lambda,Y}(\pi(1))}, \ldots,
p_{m_{\lambda,Y}(\pi(6))}) + \Phi_{\lambda+1}$.
\end{remark}

\subsection{A presentation of $J^{c+2}$ via a filtration}

For $v \in {\cal V}$ and $w(p_{m_{\lambda,Y}(1)}, \ldots,
p_{m_{\lambda,Y}(6)})$ ($Y \in {\cal G}_{\lambda}$) as above, we
denote
$$
[v;w(p_{m_{\lambda,Y}(1)}, \ldots, p_{m_{\lambda,Y}(6)})] = [v,
w(p_{m_{\lambda,Y}(1)}), \ldots, w(p_{m_{\lambda,Y}(6)})] \in
L^{c+2}.
$$
We write $[J^{2};\Phi_{\lambda}/\Phi_{\lambda+1}]$ for the
$\mathbb{Z}$-submodule of $[J^{2};T^{c}]$ spanned by all Lie
commutators of the form $[v;w(p_{m_{\lambda,Y}(1)}, \ldots,
p_{m_{\lambda,Y}(6)})]$ with $v \in {\cal V}$ and $Y \in {\cal
G}_{\lambda}$. Clearly, for $c \geq 0$,
$$
J^{c+2} = \sum_{\lambda \in {\rm Part}(c)}
[J^{2};\Phi_{\lambda}/\Phi_{\lambda+1}].
$$
For $\kappa = 1, 2, 3$, we similarly define
$[L^{2}(V_{\kappa});\Phi_{\lambda}/\Phi_{\lambda+1}]$ and
$[W^{(2)}_{2,\Psi};\Phi_{\lambda}/\Phi_{\lambda+1}]$. Since $J^{2}
= L^{2}(V_{1}) \oplus L^{2}(V_{2}) \oplus L^{2}(V_{3}) \oplus W^{(2)}_{2,\Psi}$, we
get
$$
[J^{2};\Phi_{\lambda}/\Phi_{\lambda+1}] =
[L^{2}(V_{1});\Phi_{\lambda}/\Phi_{\lambda+1}] +
[L^{2}(V_{2});\Phi_{\lambda}/\Phi_{\lambda+1}] + [L^{2}(V_{3});\Phi_{\lambda}/\Phi_{\lambda+1}] +
[W^{(2)}_{2,\Psi};\Phi_{\lambda}/\Phi_{\lambda+1}].
$$
Writing
$$
J^{c+2}_{\kappa} = \sum_{\lambda \in {\rm Part}(c)}
[L^{2}(V_{\kappa});\Phi_{\lambda}/\Phi_{\lambda+1}]~~(\kappa
=1,2,3)~~{\rm and}~~ J^{c+2}_{\Psi} = \sum_{\lambda \in {\rm
Part}(c)} [W^{(2)}_{2,\Psi};\Phi_{\lambda}/\Phi_{\lambda+1}],
$$
we have
$$
J^{c+2} = J^{c+2}_{1} + J^{c+2}_{2} + J^{c+2}_{3} + J^{c+2}_{\Psi}.
$$
For $\lambda \in {\rm Part}(c)$, and $\kappa = 1, 2, 3$, we let
$$
J^{c+2}_{\kappa,\lambda} = \sum_{\theta \in {\rm Part}(c) \atop
\lambda \leq^{*}
\theta}[L^{2}(V_{\kappa});\Phi_{\lambda}/\Phi_{\lambda+1}]~~{\rm
and}~~J^{c+2}_{\Psi,\lambda} = \sum_{\theta \in {\rm Part}(c)
\atop \lambda \leq^{*}
\theta}[W^{(2)}_{2,\Psi};\Phi_{\lambda}/\Phi_{\lambda+1}].
$$
Note that $J^{c+2}_{\kappa,\lambda} =
[L^{2}(V_{\kappa});\Phi_{\lambda}/\Phi_{\lambda+1}] +
J^{c+2}_{\kappa,\lambda+1}$ and $J^{c+2}_{\Psi,\lambda} =
[W^{(2)}_{2,\Psi};\Phi_{\lambda}/\Phi_{\lambda+1}] +
J^{c+2}_{\Psi,\lambda+1}$.  Furthermore, $J^{c+2}_{\kappa} =
J^{c+2}_{\kappa,(1^{c})}$ ($\kappa = 1, 2, 3)$ and $J^{c+2}_{\Psi} =
J^{c+2}_{\Psi,(1^{c})}$ and so,
$$
J^{c+2} = J^{c+2}_{1,(1^{c})} + J^{c+2}_{2,(1^{c})} + J^{c+2}_{3,(1^{c})} +
J^{c+2}_{\Psi,(1^{c})}.
$$

\subsection{A decomposition of $L^{2}, L^{3}$}

In this section we give a $\mathbb{Z}$-module decomposition of
$L^{2}$ and $L^{3}$.

\begin{lemma}\label{le10a}
We have the following ${\mathbb{Z}}$-module decompositions \\
{\rm (I)} $L^{2} = W^{(1)}_{2,\Psi} \oplus
J^{2}$. In particular, $L^{2}/J^{2}$ is torsion-free. \\
{\rm (II)} $J^{3} = (\bigoplus^{3}_{i=2}[L^{2}(V_{1}),V_{i}])
\oplus (\bigoplus^{3}_{i=1 \atop i\neq 2}[L^{2}(V_{2}),V_{i}]) \oplus (\bigoplus_{i=1}^{2}[L^{2}(V_{3}),V_{i}])
\oplus [[V_{3},V_{2}]^{(2)},V_{1}] \oplus J^{3}_{C}.$ \\
{\rm (III)} $L^{3} = (W^{(1)}_{3,\Psi})^{*}
\oplus J^{3}$, where $(W^{(1)}_{3,\Psi})^{*}$ denotes the
${\mathbb{Z}}$-submodule of $W^{(1)}_{3,\Psi}$ spanned by the set
$({\cal W}^{(1)}_{3,\Psi})^{*} = \{[y_{4},x_{i}], [y_{4},x_{5}], [y_{5},x_{3}], [y_{5},x_{5}], [y_{5},x_{6}],
[y_{8},x_{3}], [y_{8},x_{5}], [y_{11},x_{2}], [y_{11},x_{4}],
[y_{11},x_{6}],$ \linebreak $[y_{13},x_{5}], [y_{14},x_{3}],$
$[y_{14},x_{4}], [y_{14},x_{5}], i = 1,2,3\}$. In particular, $L^{3}/J^{3}$
is torsion-free.
\end{lemma}

\pf (I) This is straightforward.

(II) Note that $J^{3}_{\kappa} = [L^{2}(V_{\kappa}),L^{1}]$ ($\kappa = 1, 2, 3$) and $J^{3}_{\Psi} =
[W^{(2)}_{2,\Psi},L^{1}]$. Since $L^{1} = V_{1} \oplus V_{2}
\oplus V_{3}$, $[W^{(2)}_{2,\Psi},L^{1}] =
[[V_{3},V_{2}]^{(2)},V_{1}] + W^{(2)}_{3,\Psi}$ and $J^{3} =
J^{3}_{1} + J^{3}_{2} + J^{3}_{3} + J^{3}_{\Psi}$, we have
$$
J^{3} = (\sum^{3}_{i=2}[L^{2}(V_{1}),V_{i}]) + (\sum^{3}_{i=1
\atop i\neq 2}[L^{2}(V_{2}),V_{i}]) + (\sum_{i=1}^{2}[L^{2}(V_{3}),V_{i}]) + [[V_{3},V_{2}]^{(2)},V_{1}] +
J^{3}_{C}.
$$
Using the Jacobi identity and the definition of $\psi_{2}$, we have
$$
\begin{array}{lll}
~[y_{1},x_{3}] & = & -[y_{5},x_{1}] + [y_{4},x_{2}] - [\psi_{2}(y_{6}),x_{1}] + [\psi_{2}(y_{7}),x_{1}], \\
~[y_{1},x_{4}] & = & [y_{5},x_{1}] + [y_{5},x_{2}] - [\psi_{2}(y_{7}),x_{1}], \\
~[y_{1},x_{5}] & = & -[\psi_{2}(y_{10}),x_{1}]+[y_{8},x_{1}]-[y_{8},x_{2}], \\
~[y_{1},x_{6}] & = & -[y_{11},x_{1}]-[y_{8},x_{2}]-[\psi_{2}(y_{9}),x_{2}],
~~~~~(C1)
\end{array}
$$
$$
\begin{array}{lll}
~[y_{2},x_{1}] & = & [y_{5},x_{3}]-[y_{4},x_{4}], \\
~[y_{2},x_{2}] & = & -[y_{5},x_{3}]-[y_{5},x_{4}]+[\psi_{2}(y_{7}),x_{3}]-[\psi_{2}(y_{6}),x_{4}]+[\psi_{2}(y_{7}),x_{4}], \\
~[y_{2},x_{5}] & = & -[y_{14},x_{3}]+[\psi_{2}(y_{12}),x_{4}]-[y_{13},x_{4}], \\
~[y_{2},x_{6}] & = & [y_{13},x_{4}]-[\psi_{2}(y_{15}),x_{3}]+[y_{13},x_{3}], ~~~~(C2) \\
\end{array}
$$
$$
\begin{array}{lll}
~[y_{3},x_{1}] & = & [\psi_{2}(y_{9}),x_{5}] - [y_{8},x_{5}] - [y_{8},x_{6}], \\
~[y_{3},x_{2}] & = & -[\psi_{2}(y_{10}),x_{6}] + [y_{11},x_{5}] + [y_{8},x_{6}], \\
~[y_{3},x_{3}] & = & [y_{13},x_{5}] - [\psi_{2}(y_{12}),x_{6}] + [y_{13},x_{6}], \\
~[y_{3},x_{4}] & = & [\psi_{2}(y_{15}),x_{5}] - [y_{13},x_{5}] - [y_{14},x_{6}], ~~~~(C3) \\
\end{array}
$$
$$
\begin{array}{lll}
~[\psi_{2}(y_{12}),x_{1}] & = & [\psi_{2}(y_{9}),x_{3}] - [y_{4},x_{5}] - [y_{8},x_{3}] - [y_{4},x_{6}], \\
~[\psi_{2}(y_{12}),x_{2}] & = & [\psi_{2}(y_{10}),x_{3}]-[y_{8},x_{3}] - [y_{5},x_{5}] - [\psi_{2}(y_{6}),x_{5}]
+ [\psi_{2}(y_{7}),x_{5}]~ +\\
& &  [y_{11},x_{3}] - [\psi_{2}(y_{6}),x_{6}] + [\psi_{2}(y_{7}),x_{6}] - [y_{5},x_{6}], ~~~~~(C4)
\end{array}
$$
$$
\begin{array}{lll}
~[\psi_{2}(y_{15}),x_{1}] & = & [\psi_{2}(y_{9}),x_{4}] - [y_{8},x_{4}] - [y_{5},x_{6}] + [\psi_{2}(y_{9}),x_{3}]-[y_{8},x_{3}]-[y_{4},x_{6}], \\
{\rm and} & & \\
~[\psi_{2}(y_{15}),x_{2}] & =
&[y_{11},x_{4}] - [\psi_{2}(y_{6}),x_{6}] + [y_{11},x_{3}] ~~~~~(C5).
\end{array}
$$
Write $\widetilde{{\cal W}}_{3}$ for the set of the above
elements. Working in $W_{3,\Psi}$, by direct calculations, we show
that the elements of $\widetilde{{\cal W}}_{3}$ are
$\mathbb{Z}$-linear independent. Hence, the sum
$$
(\sum^{3}_{i=2}[L^{2}(V_{1}),V_{i}]) + (\sum^{3}_{i=1
\atop i\neq 2}[L^{2}(V_{2}),V_{i}]) + (\sum_{i=1}^{2}[L^{2}(V_{3}),V_{i}]) + [[V_{3},V_{2}]^{(2)},V_{1}]
$$
is direct, and $\widetilde{{\cal W}}_{3}$ is a
${\mathbb{Z}}$-basis of the aforementioned direct sum.

(III) Since $J^{3} = \langle \widetilde{{\cal W}}_{3} \rangle + J^{3}_{C}$
and $\langle \widetilde{{\cal W}}_{3} \rangle \cap J^{3}_{C} = \{0\}$, we have $
J^{3} = \langle \widetilde{{\cal W}}_{3} \rangle \oplus J^{3}_{C}$.
Next, we claim that
$$
J^{3} = (J \cap W^{(1)}_{3,\Psi}) \oplus J^{3}_{C}.
$$
It is enough to show that
$$
\langle \widetilde{{\cal W}}_{3} \rangle \subseteq (J \cap
W^{(1)}_{3,\Psi}) \oplus J^{3}_{C}
$$
which follows from the equations $(C1) - (C5)$. From the above
equations $(C1) - (C5)$, and since $J \cap W^{(1)}_{3,\Psi} \cong
\langle \widetilde{{\cal W}}_{3} \rangle$ as $\mathbb{Z}$-modules,
we may easily describe the $\mathbb{Z}$-module $J \cap
W^{(1)}_{3,\Psi}$. Let $(W^{(1)}_{3,\Psi})^{*}$ be the
${\mathbb{Z}}$-submodule of $W^{(1)}_{3,\Psi}$ spanned by the set
$({\cal W}^{(1)}_{3,\Psi})^{*} = \{[y_{4},x_{i}], [y_{4},x_{5}], [y_{5},x_{3}], [y_{5},x_{5}], [y_{5},x_{6}],$ \linebreak
$[y_{8},x_{3}], [y_{8},x_{5}], [y_{11},x_{2}], [y_{11},x_{4}],
[y_{11},x_{6}], [y_{13},x_{5}], [y_{14},x_{3}],$
$[y_{14},x_{4}], [y_{14},x_{5}], i = 1,2,3\}$. It is easily verified that
$$
W^{(1)}_{3,\Psi} = (J \cap W^{(1)}_{3,\Psi}) \oplus
(W^{(1)}_{3,\Psi})^{*}.
$$
Thus
$$
L^{3} = (W^{(1)}_{3,\Psi})^{*} \oplus J^{3}
$$
that is, the required result. \qed

\subsection{Three technical results}

In this section we shall show three technical results which are
using in the proof of our main result Theorem \ref{th1}. Let $u \in
q^{c+2}_{\rm grad}({\mathbb{L}}_{W^{(\kappa)}_{\Psi}})$. Then
there exists a unique $w \in {\mathbb{L}}_{W^{(\kappa)}_{\Psi}}$
such that $u = q(w)$ and $w$ has length $c+2$ in $x_{1}, \ldots,
x_{6}$. Since $w \in {\mathbb{L}}_{W^{(\kappa)}_{\Psi}}$, $w$ has
a \lq \lq degree \rq \rq in terms of elements of ${\cal
W}^{(\kappa)}_{\Psi}$ as well. To express these information, we
write $w^{(\kappa)}_{[s,t]}$ for a Lyndon word $w \in ({\cal
W}^{(\kappa)}_{\Psi})^{+}$ of degree $s$ (in terms of $x_{1},
\ldots, x_{6}$) and degree $t$ (in terms of the elements of ${\cal
W}^{(\kappa)}_{\Psi}$). (For example, $y_{4}, y_{11} \in {\cal
W}^{(1)}_{\Psi}$ are Lyndon words of type $w^{(1)}_{[2,1]}$, and
$[\psi_{2}(y_{10}), \psi_{2}(y_{12})]$ is a Lyndon word over
${\cal W}^{(2)}_{\Psi}$ of type $w^{(2)}_{[4,2]}$.) For positive
integers $s$ and $t$, with $s \geq 2$, let
$W^{(\kappa),t}_{s,\Psi}$ be the $\mathbb{Z}$-submodule of
$L^{s}_{\rm grad}(W^{(\kappa)}_{\Psi})$ spanned by all Lyndon
polynomials $\zeta^{(\kappa)}_{[s,t]} = q^{s}_{\rm
grad}(w^{(\kappa)}_{[s,t]})$. By the definition of $L^{s}_{\rm
grad}(W^{(\kappa)}_{\Psi})$, there exists a (unique) positive
integer $d_{s}$ such that
$$
L^{s}_{\rm grad}(W^{(\kappa)}_{\Psi}) = \bigoplus_{t = 1}^{d_{s}}
W^{(\kappa),t}_{s,\Psi}.
$$
Furthermore, we write for $e \in \{1, \ldots, d_{s}\}$
$$
L^{s}_{\rm grad,e}(W^{(\kappa)}_{\Psi}) = \bigoplus_{t =
e}^{d_{s}} W^{(\kappa),t}_{s,\Psi}.
$$
Note that $L^{s}_{\rm grad,1}(W^{(\kappa)}_{\Psi}) = L^{s}_{\rm
grad}(W^{(\kappa)}_{\Psi})$.

\begin{lemma}\label{17b}
Let $r, s$ be positive integers, with $r, s \geq 2$, and $\kappa,
\mu \in \{1, 2\}$. Let $b$ be a Lyndon polynomial in
$L^{r}(V_{\mu})$. Then, for $t \in \{1, \ldots, d_{s}\}$,
$[\zeta^{(\kappa)}_{[s,t]}, b] = \zeta^{(\kappa)}_{[s+r,t]} +
w^{\prime}_{s+r,t}$ with $w^{(\kappa)}_{[s,t]} <
w^{(\kappa)}_{[s+r,t]}$ and $w^{\prime}_{s+r,t} \in L^{s+r}_{\rm
grad}(W_{\Psi})$. Furthermore, $w^{(\kappa)}_{[s+r,t]}$ is the
smallest word in the expression $\zeta^{(\kappa)}_{[s+r,t]} +
w^{\prime}_{s+r,t}$.
\end{lemma}

\pf We shall show our claim for $\mu = 2$ and $\kappa = 1$.
Similar arguments may be applied to the other cases. Let $s \geq
2$. We induct on $t$. For $t = 1$, $W^{(1),1}_{s,\Psi} =
W^{(1)}_{s,\Psi}$. Thus $\zeta^{(1)}_{[s,1]} = q^{s}_{\rm
grad}(w^{(1)}_{[s,1]}) = w^{(1)}_{[s,1]}$ ($\in {\cal
W}^{(1)}_{s,\Psi}$) has the form
$v^{(\mu,1,\Psi)}_{s,(\alpha,\beta,\gamma)}$, $v^{(6,\nu,\Psi)}_{s,(\alpha,\beta,\gamma)}$, or
$v^{(5,4,\Psi)}_{s,(\alpha,\beta,\gamma)}$ with $\alpha + \beta +
\gamma = s-2$, $\mu \in \{3,4,5\}$ and $\nu \in \{2,3\}$. Let $b =
q^{r}(u)$, $u \in {\mathbb{L}}^{r}_{V_{2}}$, a basis element of
$L^{r}(V_{2})$. By Lemma \ref{12c}, $b = u + v$, where $v$ belongs
to the $\mathbb{Z}$-submodule of $T(V_{2})$ spanned by
$\widetilde{v} \in {\cal V}^{r}_{2}$ and $u < \widetilde{v}$. Thus
$$
[\zeta^{(1)}_{[s,1]}, b] = [\zeta^{(1)}_{[s,1]}, x_{3}, \ldots,
x_{4}] + \sum *~[\zeta^{(1)}_{[s,1]}, x_{i_{1}}, \ldots,
x_{i_{r}}], \eqno(\dag)
$$
where $u = x_{3} \cdots x_{4}$, $\widetilde{v} = x_{i_{1}} \cdots
x_{i_{r}}$, $i_{1}, \ldots, i_{r} \in \{3,4\}$ and the
coefficients $*$ are in $\mathbb{Z}$. Using the Jacobi identity in
the form $[x,y,z] = [x,z,y] + [x, [y,z]]$, by the equations $(D)$
and the definition of $\psi_{2}$, we have
$$
[\zeta^{(1)}_{[s,1]}, x_{3}, \ldots, x_{4}] =
\zeta^{(1)}_{[s+r,1]} + w_{s+r,1}
$$
where $\zeta^{(1)}_{[s+r,1]}$ has the form
$v^{(\mu,1,\Psi)}_{s+r,(\alpha,\beta + r,\gamma)}$, $v^{(6,\nu,\Psi)}_{s+r,(\alpha,\beta + r,\gamma)}$, or
$v^{(5,4,\Psi)}_{s+r,(\alpha,\beta + r,\gamma)}$, $w_{s+r,1} \in
L^{s+r}_{\rm grad,e}(W^{(1)}_{\Psi}) \oplus L^{s+r}_{\rm
grad}(\widetilde{W}_{\Psi,J})$ for some $e \geq 2$. Thus
$\zeta^{(1)}_{[s+r,1]}$ does not occur in the expression of
$w_{s+r,1}$. Note that, by the ordering on ${\cal W}_{\Psi}$
(section 3.3), $w^{(1)}_{[s,1]} < w^{(1)}_{[s+r,1]}$. We apply
similar arguments to each $[\zeta^{(1)}_{[s,1]}, x_{i_{1}},
\ldots, x_{i_{r}}]$ appearing in the right hand side of the
equation $(\dag)$. More precisely,
$$
[\zeta^{(1)}_{[s,1]}, x_{i_{1}}, \ldots, x_{i_{r}}] =
\zeta^{(1),i_{1}, \ldots,i_{r}}_{[s+r,1]} +
w_{s+r,1,i_{1},\ldots,i_{r}}
$$
where $\zeta^{(1),i_{1}, \ldots,i_{r}}_{[s+r,1]}$ has the form
$v^{(\mu,1,\Psi)}_{s+r,(\alpha,\beta + r,\gamma)}$, $v^{(6,\nu,\Psi)}_{s+r,(\alpha,\beta + r,\gamma)}$, or
$v^{(5,4,\Psi)}_{s+r,(\alpha,\beta + r,\gamma)}$,
$w_{s+r,1,i_{1},\ldots,i_{r}} \in L^{s+r}_{\rm
grad,e}(W^{(1)}_{\Psi}) \oplus L^{s+r}_{\rm
grad}(\widetilde{W}_{\Psi,J})$ for some $e \geq 2$. Since $u <
\widetilde{v}$, we obtain by the ordering on ${\cal W}_{\Psi}$,
$w^{(1)}_{[s+r,1]} < w^{(1)}_{[s+r,1,i_{1}, \ldots,i_{r}]}$ and
so, $w^{(1)}_{[s+r,1]}$ is the smallest Lyndon word occurring in
the right hand side of the equation $(\dag)$. Therefore, for all
$s, r \geq 2$,
$$
[\zeta^{(1)}_{[s,1]},b] = \zeta^{(1)}_{[s+r,1]} +
w^{\prime}_{s+r,1} ~~{\rm with} ~~w^{(1)}_{[s,1]} <
w^{(1)}_{[s+r,1]}, \eqno(\ddag)
$$
where $w^{\prime}_{s+r,1} \in L^{s+r}_{\rm grad}(W^{(1)}_{\Psi})
\oplus L^{s+r}_{\rm grad}(\widetilde{W}_{\Psi,J})$. Furthermore,
$w^{(1)}_{[s+r,1]}$ is the smallest word occurring in the right
hand side of the equation $(\ddag)$.

We assume that for $t \in \{1, \ldots, d_{s}\}$, our claim is
valid for all $t^{\prime} < t$, all $s_{1}, s_{2} < s$ and all $r
\geq 2$. We shall show our claim for $t$. Let $\zeta^{(1)}_{[s,t]}
= q^{s}_{\rm grad}(w^{(1)}_{[s,t]})$, and let
$w^{(1)}_{[s_{1},t_{1}]} \cdot w^{(1)}_{[s_{2},t_{2}]}$ ($s =
s_{1} + s_{2}$ and $t = t_{1} + t_{2}$) be the standard
factorization of $w^{(1)}_{[s,t]}$. Then $w^{(1)}_{[s_{1},t_{1}]},
w^{(1)}_{[s_{2},t_{2}]}$ are Lyndon words with
$w^{(1)}_{[s_{1},t_{1}]} < w^{(1)}_{[s_{2},t_{2}]}$ and
$w^{(1)}_{[s,t]} < w^{(1)}_{[s_{2},t_{2}]}$. Then
$$
\begin{array}{rll}
[\zeta^{(1)}_{[s,t]}, b] & = & [q^{s}_{\rm grad}(w^{(1)}_{[s,t]}),b] \\
& = & [[\zeta^{(1)}_{[s_{1},t_{1}]},\zeta^{(1)}_{[s_{2},t_{2}]}],b] \\
({\rm Jacobi~~identity})& = &
[[\zeta^{(1)}_{[s_{1},t_{1}]},b],\zeta^{(1)}_{[s_{2},t_{2}]}] -
[[\zeta^{(1)}_{[s_{2},t_{2}]},b],\zeta^{(1)}_{[s_{1},t_{1}]}].
\end{array}
$$
By our inductive hypothesis, for $j = 1,2$,
$$
[\zeta^{(1)}_{[s_{j},t_{j}]}, b] = \zeta^{(1)}_{[s_{j}+r,t_{j}]} +
w^{\prime}_{s_{j}+r,t_{j}} ~{\rm with}~ w^{(1)}_{[s_{j},t_{j}]} <
w^{(1)}_{[s_{j}+r,t_{j}]},
$$
$w^{\prime}_{s_{j}+r,t_{j}} \in L^{s_{j} + r}_{\rm
grad}(W_{\Psi})$ and $w^{(1)}_{[s_{j}+r,t_{j}]}$ is the smallest
word in the expression of $\zeta^{(1)}_{[s_{j}+r,t_{j}]} +
w^{\prime}_{s_{j}+r,t_{j}}$. Therefore
$$
[[\zeta^{(1)}_{[s_{1},t_{1}]},b],\zeta^{(1)}_{[s_{2},t_{2}]}] =
[\zeta^{(1)}_{[s_{1}+r,t_{1}]},\zeta^{(1)}_{[s_{2},t_{2}]}] +
[w^{\prime}_{s_{1}+r,t_{1}},\zeta^{(1)}_{[s_{2},t_{2}]}]
$$
and
$$
[[\zeta^{(1)}_{[s_{2},t_{2}]},b],\zeta^{(1)}_{[s_{1},t_{1}]}] =
[\zeta^{(1)}_{[s_{2}+r,t_{2}]},\zeta^{(1)}_{[s_{1},t_{1}]}] +
[w^{\prime}_{s_{2}+r,t_{2}},\zeta^{(1)}_{[s_{1},t_{1}]}].
$$
Since $[u,v] = uv - vu$ for all $u,v \in T(L^{1})$, and by using
Lemma \ref{12c}, we have
$$
[\zeta^{(1)}_{[s_{1},t_{1}]},\zeta^{(1)}_{[s_{2}+r,t_{2}]}] =
w^{(1)}_{[s_{1},t_{1}]}w^{(1)}_{[s_{2}+r,t_{2}]} + \sum *~
\widetilde{w} \eqno(\flat)
$$
where $w^{(1)}_{[s_{1},t_{1}]}w^{(1)}_{[s_{2}+r,t_{2}]} \in
{\mathbb{L}}^{s+r}_{W^{(1)}_{\Psi}}$, the coefficients $*$ are in
$\mathbb{Z}$ and $\widetilde{w} \in ({\cal W}_{\Psi})^{s+r}$ with
$w^{(1)}_{[s_{1},t_{1}]}w^{(1)}_{[s_{2}+r,t_{2}]} <
\widetilde{w}$. Since
$[\zeta^{(1)}_{[s_{1},t_{1}]},\zeta^{(1)}_{[s_{2}+r,t_{2}]}] \in
W^{(1),t}_{s+r,\Psi}$, it is (uniquely) written as a
$\mathbb{Z}$-linear combination of Lyndon polynomials of the form
$\zeta^{(1)}_{[s+r,t]}$. That is,
$$
[\zeta^{(1)}_{[s_{1},t_{1}]},\zeta^{(1)}_{[s_{2}+r,t_{2}]}] =
\sum^{m_{s,r}}_{\rho = 1} \alpha_{\rho} \zeta^{(1),\rho}_{[s+r,t]}
\eqno(\flat^{\prime})
$$
with $\alpha_{\rho} \in \mathbb{Z}$. Without loss the generality,
we assume that the corresponding Lyndon words are ordered as
$w^{(1),1}_{[s+r,t]} < \cdots < w^{(1),m_{s,t}}_{[s+r,t]}$. Since
the Lyndon word
$w^{(1)}_{[s_{1},t_{1}]}w^{(1)}_{[s_{2}+r,t_{2}]}$, say
$w^{(1)}_{[s+r,t]}$, occurs in the right hand side of the equation
$(\flat)$, it should be occurred in the right hand side of the
equation $(\flat^{\prime})$. By the ordering of
$w^{(1),\rho}_{[s+r,t]}$ $(\rho = 1, \ldots, m_{s,t})$, the
equation $(\flat)$ and Lemma \ref{12c}, we must have $\alpha_{1} =
1$ and $\zeta^{(1),1}_{[s+t,r]} = q^{s+r}_{\rm
grad}(w^{(1)}_{[s+r,t]})$. Having in mind the conditions which are
satisfied by $w^{(1)}_{[s_{1}+r,t_{1}]}$ and
$w^{(1)}_{[s_{2}+r,t_{2}]}$, and the equation $(\flat^{\prime})$,
we have
$$
[\zeta^{(1)}_{[s,t]}, b] = \zeta^{(1)}_{[s+r,t]} +
w^{\prime}_{s+r,t} ~~{\rm with}~w^{(1)}_{[s,t]} <
w^{(1)}_{[s+r,t]}
$$
and  $w^{\prime}_{s+r,t}$ satisfies the conditions of our claim.
\qed

\vskip .120 in

The following technical result can be easily shown by an inductive
argument.

\begin{lemma}\label{16a}
Let $K$ be a commutative ring with unit and characteristic $0$.
Let $T(V)$ be the tensor algebra on a free $K$-module $V$.
Consider $T(V)$ as a Lie algebra in a natural way and $L(V)$ be
the Lie subalgebra of $T(V)$ generated by $V$.

{\rm (I)} For homogeneous elements $a, b \in L(V)$ and $m \in
{\mathbb{N}}$
$$
[a, ~_{m}b] = \sum^{m}_{\kappa = 0} (-1)^{\kappa} {m \choose
\kappa} b^{\kappa}ab^{m-\kappa}.
$$

{\rm (II)} For homogeneous elements $a, b_{1}, \ldots, b_{\tau}
\in L(V)$ and $m_{1}, \ldots, m_{\tau} \in {\mathbb{N}}$
$$
[a,_{m_{1}}b_{1}, \ldots,_{m_{\tau}}b_{\tau}] =
\sum^{m_{\tau}}_{\kappa_{\tau} = 0} \cdots
\sum^{m_{1}}_{\kappa_{1} = 0}  (-1)^{\kappa_{\tau}+\cdots +
\kappa_{1}} {m_{\tau} \choose \kappa_{\tau}} \cdots {m_{1} \choose
\kappa_{1}} b^{\kappa_{\tau}}_{\tau} \cdots
b^{\kappa_{1}}_{1}ab^{m_{1}-\kappa_{1}}_{1} \cdots
b^{m_{\tau}-\kappa_{\tau}}_{\tau}.
$$
\end{lemma}

Next, we prove a technical result giving us some information about
the Lyndon polynomials over ${\cal W}^{(1)}_{\Psi}$.

\begin{lemma}\label{18a}
Let ${\cal W}^{(1)}_{\Psi} = \{z_{1}, z_{2}, \ldots \}$ with
$z_{1} \ll z_{2} \ll \cdots$ (as in section 3.3). For $i = 1,
\ldots, \tau$, let $w_{i} \in {\mathbb{L}}_{W^{(1)}_{\Psi}} \cap
({\cal W}^{(1)}_{\Psi})^{r_{i}}$, $r_{i} \in {\mathbb{N}}$, such
that $w_{\tau} < w_{\tau-1} < \cdots < w_{1}$, and let $w_{\zeta}
= \zeta_{1}^{\mu_{1}} \cdots \zeta_{\tau}^{\mu_{\tau}}$ where
$\zeta_{i} = q(w_{i})$ and $\mu_{i} \in \mathbb{N}$, $i = 1,
\ldots, \tau$. Then, for $z \in {\cal W}^{(1)}_{\Psi}$ with $z
\neq \zeta_{1}$, either $[z, ~_{\mu_{1}}\zeta_{1}, \ldots,
~_{\mu_{\tau}}\zeta_{\tau}] = q(zw^{\mu_{1}}_{1} \cdots
w^{\mu_{\tau}}_{\tau})$ with $zw^{\mu_{1}}_{1} \cdots
w^{\mu_{\tau}}_{\tau} \in {\mathbb{L}}_{W^{(1)}_{\Psi}}$ or, $[z,
~_{\mu_{1}}\zeta_{1}, \ldots, ~_{\mu_{\tau}}\zeta_{\tau}] =
w^{\mu_{\tau}}_{\tau} \cdots w^{\mu_{1}}_{1}z + v$, with
$w^{\mu_{\tau}}_{\tau} \cdots w^{\mu_{1}}_{1}z \in
{\mathbb{L}}_{W^{(1)}_{\Psi}}$ and where $v$ is a
$\mathbb{Z}$-linear combination of words $\widetilde{v} \in {\cal
W}^{+}_{\Psi}$ such that $w^{\mu_{\tau}}_{\tau} \cdots
w^{\mu_{1}}_{1}z < \widetilde{v}$ or, for unique $i \in \{2,
\ldots, \tau-1\}$, $[z, ~_{\mu_{1}}\zeta_{1}, \ldots,
~_{\mu_{\tau}}\zeta_{\tau}] = w^{\mu_{\tau}}_{\tau} \cdots
w^{\mu_{i}}_{i}zw^{\mu_{1}}_{1} \cdots w^{\mu_{i-1}}_{i-1} +
v^{\prime}$, with $w^{\mu_{\tau}}_{\tau} \cdots
w^{\mu_{i}}_{i}zw^{\mu_{1}}_{1} \cdots w^{\mu_{i-1}}_{i-1} \in
{\mathbb{L}}_{W^{(1)}_{\Psi}}$ and where $v^{\prime}$ is a
$\mathbb{Z}$-linear combination of words $\widetilde{v^{\prime}}
\in {\cal W}^{+}_{\Psi}$ such that $w^{\mu_{\tau}}_{\tau} \cdots
w^{\mu_{i}}_{i}zw^{\mu_{1}}_{1} \cdots w^{\mu_{i-1}}_{i-1} <
\widetilde{v^{\prime}}$.
\end{lemma}

\pf Since $({\cal W}^{(1)}_{\Psi})^{+}$ is a totally ordered set
(see section 3.3), we have either $z < w_{\tau}$ or $w_{1} < z$ or
there exists a unique $i \in \{2, \ldots, \tau-1\}$ such that
$w_{\tau} < \cdots < w_{i} < z < w_{i-1} < \cdots < w_{1}$.
Suppose that $z < w_{\tau}$. By Lemma \ref{12c}, $zw^{\mu_{1}}_{1}
\in {\mathbb{L}}_{W^{(1)}_{\Psi}}$, and its standard factorization
is $(zw^{\mu_{1}-1}_{1}) \cdot w_{1}$ (since $z \in {\cal
W}^{(1)}_{\Psi}$). Since $z < w_{2} < w_{1}$, we have
$(zw^{\mu_{1}}_{1})w^{\mu_{2}}_{2} \in
{\mathbb{L}}_{W^{(1)}_{\Psi}}$, and its standard factorization is
$(zw^{\mu_{1}}_{1}w^{\mu_{2}-1}_{2}) \cdot w_{2}$. Continuing in
this way, we obtain $zw^{\mu_{1}}_{1} \cdots w^{\mu_{\tau}}_{\tau}
\in {\mathbb{L}}_{W^{(1)}_{\Psi}}$, and its standard factorization
is $(zw^{\mu_{1}}_{1} \cdots w^{\mu_{\tau}-1}_{\tau}) \cdot
w_{\tau}$. Thus
$$
[z,~_{\mu_{1,4}}\zeta_{1,4},
\ldots,~_{\mu_{\tau,4}}\zeta_{\tau,4}] = q(zw^{\mu_{1,4}}_{1,4}
\cdots w^{\mu_{\tau(4),4}}_{\tau(4),4}).
$$
Next, we assume that $w_{1,4} < z$. By Lemma \ref{16a},
\begin{flushleft}
$[z,~_{\mu_{1}}\zeta_{1}, \ldots,~_{\mu_{\tau}}\zeta_{\tau}] =
(-1)^{\mu_{1}+ \cdots +\mu_{\tau}}\zeta_{\tau}^{\mu_{\tau}} \cdots
\zeta_{1}^{\mu_{1}}z + z\zeta_{1}^{\mu_{1}} \cdots
\zeta_{\tau}^{\mu_{\tau}} +$
\end{flushleft}
\begin{flushright}
$\sum *~\zeta^{\kappa_{\tau}}_{\tau} \cdots
\zeta^{\kappa_{1}}_{1}z\zeta^{\mu_{1}-\kappa_{1}}_{1} \cdots
\zeta^{\mu_{\tau}-\kappa_{\tau}}_{\tau} ~~~~(8)$
\end{flushright}
with $\kappa_{1} + \cdots + \kappa_{\tau} \geq 1$ and
$\mu_{j}-\kappa_{j} \geq 1$ for at least one $j \in \{1, \ldots,
\tau\}$, and the coefficients $*$ are in $\mathbb{Z} \setminus
\{-1,0,1\}$. By Lemma \ref{12c}, for $i \in \{1, \ldots, \tau\}$,
we have
$$
\zeta_{i} = q(w_{i}) = w_{i} + \sum_{w_{i} < \widetilde{w}_{i}}
*~\widetilde{w}_{i} \eqno (9)
$$
where $\widetilde{w}_{i} \in ({\cal W}^{(1)}_{\Psi})^{r_{i}}$, $i
= 1, \ldots, \tau$, and the coefficients $*$ are in $\mathbb{Z}$,
whose exact values are not important. By the equation (9),
$$
\zeta_{\tau}^{\mu_{\tau}} \cdots \zeta_{1}^{\mu_{1}}z =
w_{\tau}^{\mu_{\tau}} \cdots w_{1}^{\mu_{1}}z + \sum *~u_{1}u_{2}
\cdots u_{m_{\lambda,Y}(4)}z \eqno (10)
$$
where the summation runs over words $u_{1},u_{2}, \ldots,
u_{m_{\lambda,Y}(4)}$ with the first $u_{1}, \ldots,
u_{\mu_{\tau}}$ \linebreak $\geq w_{\tau}$ and $u_{1}, \ldots,
u_{\mu_{\tau}} \in ({\cal W}^{(1)}_{\Psi})^{r_{\tau}}$, the next
$u_{\mu_{\tau}+1}, \ldots, u_{\mu_{\tau}+\mu_{\tau-1}} \geq
w_{\tau-1}$ and $u_{\mu_{\tau}+1}, \ldots,$ \linebreak
$u_{\mu_{\tau}+\mu_{\tau-1}} \in ({\cal
W}^{(1)}_{\Psi})^{r_{\tau-1}}$ and so on, and $u_{j} > w_{j}$ for
at least one $j$; in this case, $u_{1}u_{2} \cdots
u_{m_{\lambda,Y}(4)}z > w_{\tau}^{\mu_{\tau}} \cdots
w_{1}^{\mu_{1}}z$. Thus, $w_{\tau}^{\mu_{\tau}} \cdots
w_{1}^{\mu_{1}}z$ does not occur in the summation in the right
hand side of the equation (10). Similar analysis, we have for
$\zeta^{\kappa_{\tau}}_{\tau} \cdots
\zeta^{\kappa_{1}}_{1}z\zeta^{\mu_{1}-\kappa_{1}}_{1} \cdots
\zeta^{\mu_{\tau}-\kappa_{\tau}}_{\tau}$ with $\kappa_{1} + \cdots
+ \kappa_{\tau} \geq 1$ and $\mu_{j}-\kappa_{j} \geq 1$ for at
least one $j \in \{1, \ldots, \tau\}$, and for
$z\zeta_{1}^{\mu_{1}} \cdots \zeta_{\tau}^{\mu_{\tau}}$. Since $z
> w_{1} > \cdots
> w_{\tau}$, we have by Lemma \ref{12c},
$w_{\tau}^{\mu_{\tau}} \cdots w_{1}^{\mu_{1}}z \in
{\mathbb{L}}_{W^{(1)}_{\Psi}}$. By the above discussion, the
Lyndon word $w_{\tau}^{\mu_{\tau}} \cdots w_{1}^{\mu_{1}}z$ is the
smallest (in the ordering of $({\cal W}^{(1)}_{\Psi})^{+}$) of the
words occurring in the equation (8). Finally, we assume that there exists a unique $i \in \{2, \ldots,
\tau-1\}$ such that
$$
w_{\tau} < \cdots < w_{i} < z < w_{i-1} < \cdots < w_{1}. \eqno
(11)
$$
Since $z < w_{i-1} < \cdots < w_{1}$, we obtain, as in the case $z
< w_{\tau}$, $zw^{\mu_{1}}_{1} \cdots w^{\mu_{i-1}}_{i-1} \in
{\mathbb{L}}_{W^{(1)}_{\Psi}}$, and its standard factorization is
$(zw^{\mu_{1}}_{1} \cdots w^{\mu_{i-1}-1}_{i-1}) \cdot w_{i-1}$.
Thus
$$
[z,~_{\mu_{1}}\zeta_{1}, \ldots,~_{\mu_{i-1}}\zeta_{i-1}] =
q(zw^{\mu_{1}}_{1} \cdots w^{\mu_{i-1}}_{i-1}).
$$
For the next few lines, we write $a = q(zw^{\mu_{1}}_{1} \cdots
w^{\mu_{i-1}}_{i-1})$ and $w = zw^{\mu_{1}}_{1} \cdots
w^{\mu_{i-1}}_{i-1}$. By the equation (11), $w^{\mu_{\tau}}_{\tau}
\cdots w^{\mu_{i}}_{i}w \in {\mathbb{L}}_{W^{(1)}_{\Psi}}$. By
Lemma \ref{16a},
\begin{flushleft}
$[a,~_{\mu_{i}}\zeta_{i}, \ldots,~_{\mu_{\tau}}\zeta_{\tau}] =
(-1)^{\mu_{i}+ \cdots +\mu_{\tau}}\zeta_{\tau}^{\mu_{\tau}} \cdots
\zeta_{1}^{\mu_{1}}a + a\zeta_{i}^{\mu_{i}} \cdots
\zeta_{\tau}^{\mu_{\tau}} +$
\end{flushleft}
\begin{flushright}
$\sum *~\zeta^{\kappa_{\tau}}_{\tau} \cdots
\zeta^{\kappa_{i}}_{i}a\zeta^{\mu_{i}-\kappa_{i}}_{i} \cdots
\zeta^{\mu_{\tau}-\kappa_{\tau}}_{\tau} ~~~~(12)$
\end{flushright}
with $\kappa_{i} + \cdots + \kappa_{\tau} \geq 1$ and
$\mu_{j}-\kappa_{j} \geq 1$ for at least one $j \in \{i, \ldots,
\tau\}$, and the coefficients $*$ are in $\mathbb{Z} \setminus
\{-1,0,1\}$. By applying similar arguments as in the case $w_{1} <
z$, we obtain the Lyndon word $w^{\mu_{\tau}}_{\tau} \cdots
w^{\mu_{i}}_{i}w$ is the smallest of the words occurring in the
equation (12).  \qed

\subsection{The main result}

In this subsection we show one of the main results of this paper.
It gives us a $\mathbb{Z}$-module decomposition of $J^{c+2}$ in
terms of $J^{c+2}_{C}$ helping us to deduce $L^{c+2}/J^{c+2}$ is
torsion-free.

Before stating the main result and giving its proof, we need to
recall some notation and definitions. Let $c \geq 2$, and $\lambda
= (c^{n(c)}, (c-1)^{n(c-1)}, \ldots, 2^{n(2)}, 1^{n(1)}) \in {\rm
Part}(c)$. Each element $w(p_{m_{\lambda,Y}(t)})$ (with $t = 1, 2,
3)$ has the form
$$
w(p_{m_{\lambda,Y}(t)}) = a_{1,t}^{(1)} \cdots
a_{1,t}^{(\nu_{t,1})}a_{2,t}^{(1)} \cdots a_{2,t}^{(\nu_{t,2})}
\cdots a_{d,t}^{(1)} \cdots a_{d,t}^{(\nu_{t,d})},
$$
where $a_{1,t}^{(1)} \preceq \cdots \preceq a_{1,t}^{(\nu_{t,1})}
\preceq \cdots \preceq a_{d,t}^{(1)} \preceq \cdots \preceq
a_{d,t}^{(\nu_{t,d})}$ with $a_{r,t}^{(1)}, \ldots,
a_{r,t}^{(\nu_{t,r})} \in q^{r}({\mathbb{L}}_{V_{t}})$ ($r = 1,
\ldots, d)$ with $d \leq c$. Let us assume that $\nu_{t,r} \geq 1$
(for $r = 1, \ldots, d$). Then, there are unique Lyndon words
$w^{(1)}_{r,t}, \ldots, w^{(\nu_{t,r})}_{r,t} \in
{\mathbb{L}}^{r}_{V_{t}}$ $(r = 1, \ldots, d)$ such that
$q(w^{(s)}_{r,t}) = a^{(s)}_{r,t}$, $r = 1, \ldots, t$, $s = 1,
\ldots, \nu_{t,r}$. By Lemma \ref{12c} (for $A = {\cal V}_{t})$
each $q(w^{(s)}_{r,t}) = w^{(s)}_{r,t} + v^{(s)}_{r,t}$, where
$v^{(s)}_{r,t}$ belongs to the $\mathbb{Z}$-submodule of
$T(V_{t})$ spanned by those words $\widetilde{v}^{(s)}_{r,t} \in
{\cal V}_{t}^{r}$ such that $w^{(s)}_{r,t} <
\widetilde{v}^{(s)}_{r,t}$. Note that for $r = 1$, $w^{(s)}_{1,t}
= a^{(s)}_{1,t} \in {\cal V}_{t}$, $s = 1, \ldots, \nu_{t,1}$ (and
so, $v^{(s)}_{1,t}$ is the empty word, $s = 1, \ldots,
\nu_{t,1}$). Define the corresponding word of
$w(p_{m_{\lambda,Y}(t)})$ as
$$
w_{\ell}(p_{m_{\lambda,Y}(t)}) = w_{1,t}^{(1)} \cdots
w_{1,t}^{(\nu_{t,1})}w_{2,t}^{(1)} \cdots w_{2,t}^{(\nu_{t,2})}
\cdots w_{d,t}^{(1)} \cdots w_{d,t}^{(\nu_{t,d})} \in {\cal
V}_{t}^{m_{\lambda,Y}(t)}.
$$
Thus we may write
$$
w(p_{m_{\lambda,Y}(t)}) = b^{\mu_{1,t}}_{1,t} b^{\mu_{2,t}}_{2,t}
\cdots b^{\mu_{\tau(t),t}}_{\tau(t),t} \eqno(13)
$$
with $\mu_{1,t}, \ldots, \mu_{\tau(t),t}$ positive integers and
$b_{i,t} \in q^{r_{i,t}}({\mathbb{L}}_{V_{t}})$, $i = 1,
\ldots, \tau(t)$ and $r_{i,t} \leq r_{j,t}$ with $i < j$ and so,
the corresponding word $w_{\ell}(p_{m_{\lambda,Y}(t)})$ is written
$w_{\ell}(p_{m_{\lambda,Y}(t)}) = w^{\mu_{1,t}}_{1,t} \cdots
w^{\mu_{\tau(t),t}}_{\tau(t),t}$ where $b_{i,t} = q(w_{i,t})$,
$i = 1, \ldots, \tau(t)$. Having in mind the basic property of
Lyndon words (see Lemma \ref{12c} for $A = {\cal V}_{t}$) we
rearrange the elements in the expression (13) as follows: there
exists a unique $\pi \in {\rm Sym}(\tau(t))$ such that
$w_{\pi(\tau(t)),t} < \cdots < w_{\pi(1),t}$ in the
lexicographical order and so, the desired rearrangement of
$w(p_{m_{\lambda,Y}(t)})$ is $w_{\pi}(p_{m_{\lambda,Y}(t)}) =
b_{\pi(1),t}^{\mu_{\pi(1),t}}b_{\pi(2),t}^{\mu_{\pi(2),t}}b_{\pi(3),t}^{\mu_{\pi(3),t}}
\cdots b_{\pi(\tau(t)),t}^{\mu_{\pi(\tau(t)),t}}$. For simplicity,
we write
$$
w(p_{m_{\lambda,Y}(t)}) =
b_{1,t}^{\mu_{1,t}}b_{2,t}^{\mu_{2,t}}b_{3,t}^{\mu_{3,t}} \cdots
b_{\tau(t),t}^{\mu_{\tau(t),t}}, ~~w_{\ell}(p_{m_{\lambda,Y}(t)})
= w_{1,t}^{\mu_{1,t}}w_{2,t}^{\mu_{2,t}}w_{3,t}^{\mu_{3,t}} \cdots
w_{\tau(t),t}^{\mu_{\tau(t),t}} \eqno(14)
$$
with $w_{\tau(t),t} < w_{\tau(t)-1,t} < \cdots < w_{1,t}$.
Replacing each $b_{r,t}$, occurring in the equation (14), by
$w_{r,t} + v_{r,t}$, we obtain from Corollary \ref{co3a}
$$
w(p_{m_{\lambda,Y}(t)}) = w_{\ell}(p_{m_{\lambda,Y}(t)}) +
\sum_{\rm finite} *~ w_{f,\ell}(p_{m_{\lambda,Y}(t)}), \eqno (15)
$$
where the coefficients $*$ are in $\mathbb{Z}$,
$w_{\ell}(p_{m_{\lambda,Y}(t)}) <
w_{f,\ell}(p_{m_{\lambda,Y}(t)})$ in the alphabetical order for
all $w_{f,\ell}(p_{m_{\lambda,Y}(t)})$ appearing in the sum of the
equation (15), and
$$
w_{f,\ell}(p_{m_{\lambda,Y}(t)}) =
f_{1,t}^{\mu_{1,t}}f_{2,t}^{\mu_{2,t}}f_{3,t}^{\mu_{3,t}} \cdots
f_{\tau(t),t}^{\mu_{\tau(t),t}}
$$
with $f_{s,t} \in
\{w_{s,t},\widetilde{v}_{s,t}\}$, $s = 1, \ldots, \tau(t)$, and at
least one $f_{s,t}$ is equal to $\widetilde{v}_{s,t}$ in the
aforementioned expression. For $\kappa = 1, 2, 3$,
$[L^{2}(V_{\kappa});\Phi_{\lambda}/\Phi_{\lambda+1}]$ is the
$\mathbb{Z}$-submodule of $[J^{2};T^{c}]$ spanned by all Lie
commutators of the form $[y_{\kappa};w(p_{m_{\lambda,Y}(1)},
\ldots, p_{m_{\lambda,Y}(6)})]$ with $Y \in {\cal G}_{\lambda}$.
We introduce some notation: Fix $Y \in {\cal G}_{\lambda}$. For
$t_{1} \in \{1,2,3\}$, with $m_{\lambda,Y}(t_{1}) \geq 1$, we
write
$$
w_{\lambda,Y}(y_{\kappa};t_{1}) =
[y_{\kappa};w(p_{m_{\lambda,Y}(t_{1})})] = [y_{\kappa},
~_{\mu_{1,t_{1}}}b_{1,t_{1}}, ~_{\mu_{2,t_{1}}}b_{2,t_{1}},
\ldots, ~_{\mu_{\tau(t_{1}),t_{1}}}b_{\tau(t_{1}),t_{1}}].
$$
For $t_{2} \in \{1,2,3\} \setminus \{t_{1}\}$ with
$m_{\lambda,Y}(t_{1}), m_{\lambda,Y}(t_{2}) \geq 1$, we write
$$
w_{\lambda,Y}(y_{\kappa};t_{1},t_{2}) =
[w_{\lambda,Y}(y_{\kappa};t_{1});w(p_{m_{\lambda,Y}(t_{2})})] =
[w_{\lambda,Y}(y_{\kappa};t_{1}), ~_{\mu_{1,t_{2}}}b_{1,t_{2}},
\ldots, ~_{\mu_{\tau(t_{2}),t_{2}}}b_{\tau(t_{2}),t_{2}}]
$$
and, for $t_{3} \neq t_{1}, t_{2}$ with $m_{\lambda,Y}(t_{3}) \geq
1$, $w_{\lambda,Y}(y_{\kappa};t_{1},t_{2},t_{3}) =
[w_{\lambda,Y}(y_{\kappa};t_{1},t_{2});w(p_{m_{\lambda,Y}(t_{3})}]$.
Similarly, we write for $w_{\lambda,Y,\ell}(y_{\kappa};t_{1})$,
$w_{\lambda,Y,\ell}(y_{\kappa};t_{1},t_{2})$,
$w_{\lambda,Y,\ell}(y_{\kappa};t_{1},t_{2},t_{3})$,
$w_{\lambda,Y,\ell,f}(y_{\kappa};t_{1})$, \linebreak
$w_{\lambda,Y,\ell,f}(y_{\kappa};t_{1},t_{2})$ and
$w_{\lambda,Y,\ell,f}(y_{\kappa};t_{1},t_{2},t_{3})$. Furthermore,
for $u \in L^{\prime}$ and $v = z_{1} \cdots z_{m} \in
\bigcup_{j=1}^{3}{\cal V}^{+}_{j}$, we write $f(u;v) = [u,z_{1},
\ldots, z_{m}]$. It is clearly enough that $f(u;v_{1}v_{2}) =
f(f(u;v_{1});v_{2})$ for all $u \in L^{\prime}$ and $v_{1}, v_{2}
\in \bigcup_{j=1}^{3}{\cal V}^{+}_{j}$.

\vskip .120 in

Our main result in this section is the following.

\begin{theorem}\label{th1}
Let $c \geq 2$. There are ${\mathbb{Z}}$-submodules $U^{c+2}_{1}$,
$U^{c+2}_{2}$, $U^{c+2}_{3}$ and $U^{c+2}_{\Psi}$ of $J^{c+2}_{1}$, $J^{c+2}_{2}$, $J^{c+2}_{3}$
and $J^{c+2}_{\Psi}$, respectively, such that
$$
J^{c+2} = U^{c+2}_{1} \oplus U^{c+2}_{2} \oplus U^{c+2}_{3} \oplus U^{c+2}_{\Psi}
\oplus J^{c+2}_{C} = (J \cap L^{c+2}_{\rm grad}(W^{(1)}_{\Psi}))
\oplus J^{c+2}_{C}.
$$
Furthermore, $J \cap L^{c+2}_{\rm grad}(W^{(1)}_{\Psi})$ is a
direct summand of $L^{c+2}_{\rm grad}(W^{(1)}_{\Psi})$. In
particular,
$$
L^{c+2} = (L^{c+2}_{\rm
grad}(W^{(1)}_{\Psi}))^{*} \oplus J^{c+2},
$$
where $L^{c+2}_{\rm grad}(W^{(1)}_{\Psi}) = (L^{c+2}_{\rm
grad}(W^{(1)}_{\Psi}))^{*} \oplus (J \cap L^{c+2}_{\rm
grad}(W^{(1)}_{\Psi}))$.
\end{theorem}

\vskip .120 in

\pf We shall prove our claim into four cases.

\vskip .120 in

\emph{Case 1. An analysis of $J^{c+2}_{1,\lambda}$ modulo
$(J^{c+2}_{1,\lambda+1}+J^{c+2}_{C})$.} Recall that, for $\lambda
\in {\rm Part}(c)$,
$$
J^{c+2}_{1,\lambda} + J^{c+2}_{C} =
[L^{2}(V_{1});\Phi_{\lambda}/\Phi_{\lambda+1}] +
J^{c+2}_{1,\lambda+1} + J^{c+2}_{C}
$$
and
$$
[L^{2}(V_{1});\Phi_{\lambda}/\Phi_{\lambda+1}] = {\mathbb{Z}}-{\rm
span} \{[y_{1};w(p_{m_{\lambda,Y}(1)}, \ldots,
p_{m_{\lambda,Y}(6)})], Y \in {\cal G}_{\lambda}\}.
$$
In the next subcases, we analyze the Lie commutator
$[y_{1};w(p_{m_{\lambda,Y}(1)}, \ldots, p_{m_{\lambda,Y}(6)})]$
with $Y \in {\cal G}_{\lambda}$. So, we fix $Y \in {\cal
G}_{\lambda}$.

\vskip .120 in

\emph{Subcase 1a. $m_{\lambda,Y}(1) \geq 1$.} Since
$w_{\lambda,Y}(y_{1};1) \neq 0$, we have $x_{1}x_{2} \neq
w_{1,1}$. By the equation (15) (for $t = 1$), we get
$$
w_{\lambda,Y}(y_{1};1) = w_{\lambda,Y,\ell}(y_{1};1) + \sum_{\rm
finite} *~ w_{\lambda,Y,\ell,f}(y_{1};1).
$$
If each $m_{\lambda,Y}(j) = 0$ for $j = 2, \ldots, 6$, we have
$[y_{1};w(p_{m_{\lambda,Y}(1)}, \ldots, p_{m_{\lambda,Y}(6)})]$
$= w_{\lambda,Y}(y_{1};1) \in J^{c+2}_{C}$. Thus we assume $m_{\lambda,Y}(j) \geq 1$ for some $j \in \{2,
\ldots, 6\}$. Suppose that $m_{\lambda,Y}(2) \geq 1$. If $w_{1,2}
= x_{4}$, then $w_{\ell}(p_{m_{\lambda,Y}(2)}) =
x^{\mu_{1,2}}_{4}$ and $w(p_{m_{\lambda,Y}(2)}) =
x^{\mu_{1,2}}_{4}$. Using the Jacobi identity in the form $[x,y,z]
= [x,z,y] + [x,[y,z]]$ and since $\psi_{2}(y_{7}) = y_{7}+y_{5}$, we have
$$
[w_{\lambda,Y}(y_{1};1), x_{4}] =
f([y_{1},x_{4}];w_{\ell}(p_{m_{\lambda,Y}(1)}))+w,
$$
where $w \in L^{3+m_{\lambda,Y}(1)}_{\rm grad, e}(W_{\Psi})$
with $e \geq 2$. By the appropriate equation of (C1),
\begin{flushleft}
$[w_{\lambda,Y}(y_{1};1),~_{\mu_{1,2}} x_{4}] =
f([y_{5},x_{1}];w_{\ell}(p_{m_{\lambda,Y}(1)})x^{\mu_{1,2}-1}_{4})+f([y_{5},x_{2}];w_{\ell}(p_{m_{\lambda,Y}(1)})x^{\mu_{1,2}-1}_{4})$
\end{flushleft}
\begin{flushright}
$-f([\psi_{2}(y_{7}),x_{1}];w_{\ell}(p_{m_{\lambda,Y}(1)})x^{\mu_{1,2}-1}_{4})+[w,~_{(\mu_{1,2}-1)}x_{4}].
~~~~~~~~~~(16)$
\end{flushright}
We notice that $[w,~_{(\mu_{1,2}-1)}x_{4}] = w^{(1)}_{(1,4)}+w^{(2)}_{(1,4)}+\widetilde{w}_{(1,4)}$
where $w^{(\kappa)}_{(1,4)} \in L^{2+m_{\lambda,Y}(1)+\mu_{1,2}}_{{\rm grad},e_{\kappa}}(W^{(\kappa)}_{\Psi})$ $(e_{\kappa} \geq 2,
\kappa = 1, 2)$ and $\widetilde{w}_{(1,4)} \in L^{2+m_{\lambda,Y}(1)+\mu_{1,2}}_{{\rm grad},e_{\kappa}}(\widetilde{W}_{\Psi,J})$. Furthermore,
$$
f([y_{5},x_{1}];w_{\ell}(p_{m_{\lambda,Y}(1)})x^{\mu_{1,2}-1}_{4}), f([y_{5},x_{2}];w_{\ell}(p_{m_{\lambda,Y}(1)})x^{\mu_{1,2}-1}_{4})
\in {\cal W}^{(1)}_{2+m_{\lambda,Y}(1)+\mu_{1,2},\Psi},
$$
$$
f([\psi_{2}(y_{7}),x_{1}];w_{\ell}(p_{m_{\lambda,Y}(1)})x^{\mu_{1,2}-1}_{4}) \in {\cal W}^{(2)}_{2+m_{\lambda,Y}(1)+\mu_{1,2},\Psi},
$$
and
$$
f([y_{5},x_{1}];w_{\ell}(p_{m_{\lambda,Y}(1)})x^{\mu_{1,2}-1}_{4}) + f([y_{5},x_{2}];w_{\ell}(p_{m_{\lambda,Y}(1)})x^{\mu_{1,2}-1}_{4}) +
w^{(1)}_{(1,4)} \in J \cap
L^{2+m_{\lambda,Y}(1)+\mu_{1,2}}_{\rm grad}(W^{(1)}_{\Psi}).
$$
Thus, we assume that $w_{1,2} \neq x_{4}$. Each element in
${\mathbb{L}}_{V_{2}}$, but not $x_{4}$, starts with $x_{3}$ and
ends with $x_{4}$. Using the Jacobi identity as before, the
suitable equations of (C1) and $y_{6} (= [x_{3},x_{2}]) = \psi_{2}(y_{6})-\psi_{2}(y_{7})+y_{5}$,
$$
\begin{array}{lll}
w_{\lambda,Y}(y_{1};1,2) & = &
-f([y_{5},x_{1}];w_{\ell}(p_{m_{\lambda,Y}(1)})w^{\prime}_{1,2}w^{{(\mu_{1,2}-1)}}_{1,2}
\ldots w^{\mu_{\tau(2),2}}_{\tau(2),2}) + \\
& & f([y_{4},x_{2}];w_{\ell}(p_{m_{\lambda,Y}(1)})w^{\prime}_{1,2}w^{{(\mu_{1,2}-1)}}_{1,2}
\ldots w^{\mu_{\tau(2),2}}_{\tau(2),2})~- \\
& & f([\psi_{2}(y_{6}),x_{1}];w_{\ell}(p_{m_{\lambda,Y}(1)})w^{\prime}_{1,2}w^{{(\mu_{1,2}-1)}}_{1,2}
\ldots w^{\mu_{\tau(2),2}}_{\tau(2),2})~+ \\
& & f([\psi_{2}(y_{7}),x_{1}];w_{\ell}(p_{m_{\lambda,Y}(1)})w^{\prime}_{1,2}w^{{(\mu_{1,2}-1)}}_{1,2}
\ldots w^{\mu_{\tau(2),2}}_{\tau(2),2}) + w^{(1,2)},~~~~~~~~(17)
\end{array}
$$
where $w^{\prime}_{1,2}$ is the unique element in ${\cal
V}^{*}_{2}$ (the free monoid on ${\cal V}_{2}$) such that $w_{1,2}
= x_{3}w^{\prime}_{1,2}$, and $w^{(1,2)} \in
L^{\rho_{\lambda,Y}(1,2)+2}_{\rm grad, e}(W_{\Psi})$ with $e
\geq 2$ and $\rho_{\lambda,Y}(1,2) =
m_{\lambda,Y}(1)+m_{\lambda,Y}(2)$. In particular, $w^{(1,2)} = w^{(1,2,1)}_{(1,3)} + w^{(1,2,2)}_{(1,3)}+\widetilde{w}^{(1,2)}_{(1,3)}$,
where $w^{(1,2,\kappa)}_{(1,3)} \in L^{2+\rho_{\lambda,Y}(1,2)}_{{\rm grad},e_{\kappa}}(W^{(\kappa)}_{\Psi})$
$(e_{\kappa} \geq 2, \kappa = 1,2)$ and $\widetilde{w}^{(1,2)}_{(1,3)} \in L^{2+\rho_{\lambda,Y}(1,2)}_{{\rm grad},e}(\widetilde{W}_{\Psi,J})$, $e \geq 2$.
Note that
$$
\begin{array}{ll}
f([y_{5},x_{1}];w_{\ell}(p_{m_{\lambda,Y}(1)})w^{\prime}_{1,2}w^{{(\mu_{1,2}-1)}}_{1,2}
\ldots w^{\mu_{\tau(2),2}}_{\tau(2),2}), & \\
f([y_{4},x_{2}];w_{\ell}(p_{m_{\lambda,Y}(1)})w^{\prime}_{1,2}w^{{(\mu_{1,2}-1)}}_{1,2}
\ldots w^{\mu_{\tau(2),2}}_{\tau(2),2}) \in {\cal W}^{(1)}_{2+\rho_{\lambda,Y}(1,2),\Psi},
&
\end{array}
$$
$$
\begin{array}{ll}
f([\psi_{2}(y_{6}),x_{1}];w_{\ell}(p_{m_{\lambda,Y}(1)})w^{\prime}_{1,2}w^{{(\mu_{1,2}-1)}}_{1,2}
\ldots w^{\mu_{\tau(2),2}}_{\tau(2),2}), & \\
f([\psi_{2}(y_{7}),x_{1}];w_{\ell}(p_{m_{\lambda,Y}(1)})w^{\prime}_{1,2}w^{{(\mu_{1,2}-1)}}_{1,2}
\ldots w^{\mu_{\tau(2),2}}_{\tau(2),2}) \in {\cal W}^{(2)}_{2+\rho_{\lambda,Y}(1,2),\Psi},
&
\end{array}
$$
and
$$
\begin{array}{ll}
-f([y_{5},x_{1}];w_{\ell}(p_{m_{\lambda,Y}(1)})w^{\prime}_{1,2}w^{{(\mu_{1,2}-1)}}_{1,2}
\ldots w^{\mu_{\tau(2),2}}_{\tau(2),2}) & ~+ \\
f([y_{4},x_{2}];w_{\ell}(p_{m_{\lambda,Y}(1)})w^{\prime}_{1,2}w^{{(\mu_{1,2}-1)}}_{1,2}
\ldots w^{\mu_{\tau(2),2}}_{\tau(2),2}) & ~+ \\
& w^{(1,2,1)}_{(1,3)} \in J \cap
L^{2+\rho_{\lambda,Y}(1,2)}_{\rm grad}(W^{(1)}_{\Psi}).
\end{array}
$$
If $m_{\lambda,Y}(3) \geq 1$ and $m_{\lambda,Y}(2) = 0$, then, by
applying similar arguments as above and $y_{9} = \psi_{2}(y_{9}) - y_{8}$, we have by the suitable equations of
(C1),
\begin{flushleft}
$[w_{\lambda,Y}(y_{1};1),~_{\mu_{1,3}} x_{6}] =
-f([y_{11},x_{1}];w_{\ell}(p_{m_{\lambda,Y}(1)})x^{\mu_{1,3}-1}_{6})-f([y_{8},x_{2}];w_{\ell}(p_{m_{\lambda,Y}(1)})x^{\mu_{1,3}-1}_{6})+$
\end{flushleft}
\begin{flushright}
$-f([\psi_{2}(y_{9}),x_{2}];w_{\ell}(p_{m_{\lambda,Y}(1)})x^{\mu_{1,3}-1}_{6})+
[w_{1},~_{(\mu_{1,3}-1)}x_{6}],~~~~~(18)$
\end{flushright}
where $w_{1} \in L^{3+m_{\lambda,Y}(1)}_{\rm grad, e}(W_{\Psi})$
with $e \geq 2$. In particular, $[w_{1},~_{(\mu_{1,3}-1)}x_{6}] =
w^{(1)}_{(1,6)} + w^{(2)}_{(1,6)} + \widetilde{w}_{(1,6)}$ where
$w^{(\kappa)}_{(1,6)} \in L^{2+m_{\lambda,Y}(1)+\mu_{1,3}}_{\rm
grad, e_{\kappa}}(W^{(\kappa)}_{\Psi})$ ($e_{\kappa} \geq 2,
\kappa = 1,2$) and $\widetilde{w}_{(1,6)} \in
L^{2+m_{\lambda,Y}(1)+\mu_{1,3}}_{\rm grad,
e}(\widetilde{W}_{\Psi,J})$, $e \geq 2$. Note that
$$
f([y_{11},x_{1}];w_{\ell}(p_{m_{\lambda,Y}(1)})x^{\mu_{1,3}-1}_{6}),
~f([y_{8},x_{2}];
w_{\ell}(p_{m_{\lambda,Y}(1)})x^{\mu_{1,3}-1}_{6}) \in {\cal
W}^{(1)}_{m_{\lambda,Y}(1)+\mu_{1,3}+2,\Psi}
$$
and
$$
\begin{array}{ll}
f([\psi_{2}(y_{9}),x_{2}];w_{\ell}(p_{m_{\lambda,Y}(1)})x^{\mu_{1,3}-1}_{6})
\in {\cal W}^{(2)}_{m_{\lambda,Y}(1)+\mu_{1,3}+2,\Psi}. &
\end{array}
$$
Furthermore,
\begin{flushleft}
$-f([y_{11},x_{1}];w_{\ell}(p_{m_{\lambda,Y}(1)})x^{\mu_{1,3}-1}_{6})
- f([y_{8},x_{2}];
w_{\ell}(p_{m_{\lambda,Y}(1)})x^{\mu_{1,3}-1}_{6}) +$
\end{flushleft}
\begin{flushright}
$w^{(1)}_{(1,6)} \in J \cap L^{2+m_{\lambda,Y}(1)+\mu_{1,3}}_{\rm
grad}(W^{(1)}_{\Psi})$.
\end{flushright}
Thus, we assume that $w_{1,3} \neq x_{6}$. Each element in
${\mathbb{L}}_{V_{5}}$, but not $x_{6}$, starts with $x_{5}$ and
ends with $x_{6}$. By using the Jacobi identity as before, $y_{10} = \psi_{2}(y_{10}) - y_{8}$ and the
suitable equations of (C1), we have
$$
\begin{array}{lll}
w_{\lambda,Y}(y_{1};1,3) & = &
f([y_{8},x_{1}];w_{\ell}(p_{m_{\lambda,Y}(1)})w^{\prime}_{1,3}w^{{(\mu_{1,3}-1)}}_{1,3}
\ldots w^{\mu_{\tau(3),3}}_{\tau(3),3})~- \\
& & f([y_{8},x_{2}];w_{\ell}(p_{m_{\lambda,Y}(1)})w^{\prime}_{1,3}w^{{(\mu_{1,3}-1)}}_{1,3}
\ldots w^{\mu_{\tau(3),3}}_{\tau(3),3})~ + \\
& & -f([\psi_{2}(y_{10}),x_{1}];w_{\ell}(p_{m_{\lambda,Y}(1)})w^{\prime}_{1,3}w^{{(\mu_{1,3}-1)}}_{1,3}
\ldots w^{\mu_{\tau(3),3}}_{\tau(3),3})+w^{(1,3)},~~~~~~(19)
\end{array}
$$
where $w^{\prime}_{1,3}$ is the unique element in ${\cal
V}^{*}_{3}$ (the free monoid on ${\cal V}_{3}$) such that $w_{1,3}
= x_{5}w^{\prime}_{1,3}$, and $w^{(1,3)} \in
L^{\rho_{\lambda,Y}(1,3)+2}_{\rm grad}(W_{\Psi})$ with
$\rho_{\lambda,Y}(1,3) = m_{\lambda,Y}(1)+m_{\lambda,Y}(3)$. In
fact, $w^{(1,3)} = w^{(1,3)}_{1} + w^{(1,3)}_{2} +
w^{(1,3)}_{\Psi,J}$ where $w^{(1,3)}_{\kappa} \in
L^{\rho_{\lambda,Y}(1,3)+2}_{\rm grad,
e_{\kappa}}(W^{(\kappa)}_{\Psi})$ ($e_{\kappa} \geq 2, \kappa =
1,2$) and $w^{(1,3)} \in L^{\rho_{\lambda,Y}(1,3)+2}_{\rm grad,
e}(\widetilde{W}_{\Psi,J})$, $e \geq 2$. Note that
$$
\begin{array}{ll}
f([y_{8},x_{1}];w_{\ell}(p_{m_{\lambda,Y}(1)})w^{\prime}_{1,3}w^{{(\mu_{1,3}-1)}}_{1,3}
\ldots w^{\mu_{\tau(3),3}}_{\tau(3),3}), & \\
f([y_{8},x_{2}];w_{\ell}(p_{m_{\lambda,Y}(1)})w^{\prime}_{1,3}w^{{(\mu_{1,3}-1)}}_{1,3}
\ldots w^{\mu_{\tau(3),3}}_{\tau(3),3}) \in {\cal
W}^{(1)}_{\rho_{\lambda,Y}(1,3)+2,\Psi},
\end{array}
$$
$$
f([\psi_{2}(y_{10}),x_{1}];w_{\ell}(p_{m_{\lambda,Y}(1)})w^{\prime}_{1,3}w^{{(\mu_{1,3}-1)}}_{1,3}
\ldots w^{\mu_{\tau(3),3}}_{\tau(3),3}) \in {\cal
W}^{(2)}_{\rho_{\lambda,Y}(1,3)+2,\Psi}
$$
and
\begin{flushleft}
$f([y_{8},x_{1}];w_{\ell}(p_{m_{\lambda,Y}(1)})w^{\prime}_{1,3}
\ldots w^{\mu_{\tau(3),3}}_{\tau(3),3}) -
f([y_{8},x_{2}];w_{\ell}(p_{m_{\lambda,Y}(1)})w^{\prime}_{1,3}
\ldots w^{\mu_{\tau(3),3}}_{\tau(3),3}) +$
\end{flushleft}
\begin{flushright}
$w^{(1,3)}_{1} \in J \cap L^{\rho_{\lambda,Y}(1,3)+2}_{\rm
grad}(W^{(1)}_{\Psi})$.
\end{flushright}
>From the equations (16) and (15) (for $t = 3)$,
$$
\begin{array}{l}
f([w_{\lambda,Y}(y_{1};1),~_{\mu_{1,2}}
x_{4}];w_{\ell}(p_{\lambda,Y}(3))) =
f([y_{5},x_{1}];w_{\ell}(p_{m_{\lambda,Y}(1)})x^{\mu_{1,2}-1}_{4}w_{\ell}(p_{\lambda,Y}(3)))~+ \\
 f([y_{5},x_{2}];w_{\ell}(p_{m_{\lambda,Y}(1)})x^{\mu_{1,2}-1}_{4}w_{\ell}(p_{\lambda,Y}(3)))~- \\
 f([\psi_{2}(y_{7}),x_{1}];w_{\ell}(p_{m_{\lambda,Y}(1)})x^{\mu_{1,2}-1}_{4}w_{\ell}(p_{\lambda,Y}(3)))+w_{2},
 ~~~~~~~~~~~~~~~~~~~~~~~~~~~~~~~~~~~~~~~~(20)
\end{array}
$$
where $w_{2} \in L^{2+\rho_{\lambda,Y}(1,2,3)}_{\rm grad,
e}(W_{\Psi})$, $e \geq 2$, $\rho_{\lambda,Y}(1,2,3) =
m_{\lambda,Y}(1)+m_{\lambda,Y}(2)+m_{\lambda,Y}(3)$. In
particular, $w_{2} = w^{(1)}_{2} + w^{(2)}_{2} +
\widetilde{w}_{2}$, where $w^{(\kappa)}_{2} \in
L^{2+\rho_{\lambda,Y}(1,2,3)}_{{\rm
grad},e_{\kappa}}(W^{(\kappa)}_{\Psi})$ $( e_{\kappa} \geq 2,
\kappa = 1, 2$) and $\widetilde{w}_{2} \in
L^{2+\rho_{\lambda,Y}(1,2,3)}_{{\rm
grad},e}(\widetilde{W}_{\Psi,J})$. The elements
$$
\begin{array}{ll}
f([y_{5},x_{1}];w_{\ell}(p_{m_{\lambda,Y}(1)})x^{\mu_{1,2}-1}_{4}w_{\ell}(p_{\lambda,Y}(3))),
& \\
f([y_{5},x_{2}];w_{\ell}(p_{m_{\lambda,Y}(1)})x^{\mu_{1,2}-1}_{4}w_{\ell}(p_{\lambda,Y}(3)))
\in {\cal W}^{(1)}_{\rho_{\lambda,Y}(1,2,3)+2,\Psi}, &
\end{array}
$$
$$
\begin{array}{ll}
f([\psi_{2}(y_{7}),x_{1}];w_{\ell}(p_{m_{\lambda,Y}(1)})x^{\mu_{1,2}-1}_{4}w_{\ell}(p_{\lambda,Y}(3)))
\in {\cal W}^{(2)}_{\rho_{\lambda,Y}(1,2,3)+2,\Psi} &
\end{array}
$$
and
$$
\begin{array}{l}
f([y_{5},x_{1}];w_{\ell}(p_{m_{\lambda,Y}(1)})x^{\mu_{1,2}-1}_{4}w_{\ell}(p_{\lambda,Y}(3)))
+
\\ f([y_{5},x_{2}];w_{\ell}(p_{m_{\lambda,Y}(1)})x^{\mu_{1,2}-1}_{4}w_{\ell}(p_{\lambda,Y}(3)))
+ w^{(1)}_{2} \in J \cap L^{2+\rho_{\lambda,Y}(1,2,3)}_{\rm
grad}(W^{(1)}_{\Psi}).
\end{array}
$$
Finally, from the equations (17) and (15) (for $t = 3$), we have
$$
\begin{array}{l}
w_{\lambda,Y}(y_{1};1,2,3) =
-f([y_{5},x_{1}];w_{\ell}(p_{m_{\lambda,Y}(1)})w^{\prime}_{1,2}w^{{(\mu_{1,2}-1)}}_{1,2}
\ldots w^{\mu_{\tau(2),2}}_{\tau(2),2}w^{\mu_{1,3}}_{1,3} \ldots
w^{\mu_{\tau(3),3}}_{\tau(3),3}) + \\
f([y_{4},x_{2}];w_{\ell}(p_{m_{\lambda,Y}(1)})w^{\prime}_{1,2}w^{{(\mu_{1,2}-1)}}_{1,2}
\ldots w^{\mu_{\tau(2),2}}_{\tau(2),2}w^{\mu_{1,3}}_{1,3} \ldots
w^{\mu_{\tau(3),3}}_{\tau(3),3}) - \\
f([\psi_{2}(y_{6}),x_{1}];w_{\ell}(p_{m_{\lambda,Y}(1)})w^{\prime}_{1,2}w^{{(\mu_{1,2}-1)}}_{1,2}
\ldots w^{\mu_{\tau(2),2}}_{\tau(2),2}w^{\mu_{1,3}}_{1,3} \ldots
w^{\mu_{\tau(3),3}}_{\tau(3),3}) + \\
f([\psi_{2}(y_{7}),x_{1}];w_{\ell}(p_{m_{\lambda,Y}(1)})w^{\prime}_{1,2}w^{{(\mu_{1,2}-1)}}_{1,2}
\ldots w^{\mu_{\tau(2),2}}_{\tau(2),2}w^{\mu_{1,3}}_{1,3} \ldots
w^{\mu_{\tau(3),3}}_{\tau(3),3})+w^{(1,2,3)}, ~~~~~~(21)
\end{array}
$$
where $w^{(1,2,3)} \in L^{2+\rho_{\lambda,Y}(1,2,3)}_{\rm grad,
e}(W_{\Psi})$, $e \geq 2$. In particular, $w^{(1,2,3)} =
w^{(1,2,3)}_{1}+w^{(1,2,3)}_{2}+\widetilde{w}^{(1,2,3)}$, where
$w^{(1,2,3)}_{\kappa} \in L^{2+\rho_{\lambda,Y}(1,2,3)}_{{\rm
grad},e_{\kappa}}(W^{(\kappa)}_{\Psi})$ $(e_{\kappa} \geq 2,
\kappa = 1, 2)$ and $\widetilde{w}^{(1,2,3)} \in
L^{2+\rho_{\lambda,Y}(1,2,3)}_{{\rm
grad},e}(\widetilde{W}_{\Psi,J})$. The elements
\begin{flushleft}
$f([y_{5},x_{1}];w_{\ell}(p_{m_{\lambda,Y}(1)})w^{\prime}_{1,2}w^{{(\mu_{1,2}-1)}}_{1,2}
\ldots w^{\mu_{\tau(2),2}}_{\tau(2),2}w^{\mu_{1,3}}_{1,3} \ldots
w^{\mu_{\tau(3),3}}_{\tau(3),3}),$
\end{flushleft}
\begin{flushright}
$f([y_{4},x_{2}];w_{\ell}(p_{m_{\lambda,Y}(1)})w^{\prime}_{1,2}w^{{(\mu_{1,2}-1)}}_{1,2}
\ldots w^{\mu_{\tau(2),2}}_{\tau(2),2}w^{\mu_{1,3}}_{1,3} \ldots
w^{\mu_{\tau(3),3}}_{\tau(3),3}) \in {\cal
W}^{(1)}_{2+\rho_{\lambda,Y}(1,2,3),\Psi},$
\end{flushright}
\begin{flushleft}
$f([\psi_{2}(y_{6}),x_{1}];w_{\ell}(p_{m_{\lambda,Y}(1)})w^{\prime}_{1,2}w^{{(\mu_{1,2}-1)}}_{1,2}
\ldots w^{\mu_{\tau(2),2}}_{\tau(2),2}w^{\mu_{1,3}}_{1,3} \ldots
w^{\mu_{\tau(3),3}}_{\tau(3),3}),$
\end{flushleft}
\begin{flushright}
$f([\psi_{2}(y_{7}),x_{1}];w_{\ell}(p_{m_{\lambda,Y}(1)})w^{\prime}_{1,2}w^{{(\mu_{1,2}-1)}}_{1,2}
\ldots w^{\mu_{\tau(2),2}}_{\tau(2),2}w^{\mu_{1,3}}_{1,3} \ldots
w^{\mu_{\tau(3),3}}_{\tau(3),3}) \in {\cal
W}^{(2)}_{2+\rho_{\lambda,Y}(1,2,3),\Psi},$
\end{flushright}
and
$$
\begin{array}{l}
-f([y_{5},x_{1}];w_{\ell}(p_{m_{\lambda,Y}(1)})w^{\prime}_{1,2}w^{{(\mu_{1,2}-1)}}_{1,2}
\ldots w^{\mu_{\tau(2),2}}_{\tau(2),2}w^{\mu_{1,3}}_{1,3} \ldots
w^{\mu_{\tau(3),3}}_{\tau(3),3}) + \\
f([y_{4},x_{2}];w_{\ell}(p_{m_{\lambda,Y}(1)})w^{\prime}_{1,2}w^{{(\mu_{1,2}-1)}}_{1,2}
\ldots w^{\mu_{\tau(2),2}}_{\tau(2),2}w^{\mu_{1,3}}_{1,3} \ldots
w^{\mu_{\tau(3),3}}_{\tau(3),3}) + \\ w^{(1,2,3)}_{1} \in J \cap
L^{2+\rho_{\lambda,Y}(1,2,3)}_{\rm grad}(W^{(1)}_{\Psi}).
\end{array}
$$
If $m_{\lambda,Y}(5)$ or $m_{\lambda,Y}(6) \geq 1$, we obtain from
the equations $(16)-(21)$ and, by the construction of
$L(\widetilde{W}_{\Psi,J})$, the corresponding Lie commutator $
[y_{1};w(p_{m_{\lambda,Y}(1)}, \ldots,p_{m_{\lambda,Y}(6)})] \in
J^{c+2}_{C}.$ Thus we assume that $m_{\lambda,Y}(5) = m_{\lambda,Y}(6) = 0$. Let
$m_{\lambda,Y}(4) \geq 1$. Without loss of generality, we consider
the case where both $m_{\lambda,Y}(2)$ and $m_{\lambda,Y}(3) \geq
1$. (We deal with the same method the cases either
$m_{\lambda,Y}(2) \geq 1$ or $m_{\lambda,Y}(3) \geq 1$.) Thus, by
the equations (20), (21), we may write
$$
w_{\lambda,Y}(y_{1};1,2,3) = \pm v_{1} + v_{2} + w_{(1,2)},
$$
with $v_{1}, v_{2} \in {\cal W}^{(1)}_{m,\Psi}$ (that is, $v_{1},
v_{2}$ are Lyndon words of length $1$) with $m = 2 +
\rho_{\lambda,Y}(1,2,3)$ and $w_{(1,2)} \in W^{(2)}_{\Psi} \oplus
L^{m}_{\rm grad, e}(W_{\Psi})$, $e \geq 2$. In particular,
$w_{(1,2)} =
w^{(1)}_{(1,2)}+w^{(1,2)}_{(1,2)}+w^{(2)}_{(1,2)}+\widetilde{w}_{(1,2)}$
where $w^{(\kappa)}_{(1,2)} \in L^{m}_{{\rm
grad},e_{\kappa}}(W^{(\kappa)}_{\Psi})$ $(e_{\kappa} \geq 2,
\kappa = 1, 2)$, $w^{(1,2)}_{(1,2)} \in {\cal W}^{(2)}_{m,\Psi}$
and $\widetilde{w}_{(1,2)} \in L^{m}_{{\rm
grad},e}(\widetilde{W}_{\Psi,J})$, $e \geq 2$. Recall that
$w(p_{m_{\lambda,Y}(4)}) = \zeta^{\mu_{1,4}}_{1,4}
\zeta^{\mu_{2,4}}_{2,4} \cdots
\zeta^{\mu_{\tau(4),4}}_{\tau(4),4}$ with $q(w_{i,4}) =
\zeta_{i,4}$ and $w_{\tau(4),4} < \cdots < w_{1,4}$  (where $<$
denotes the alphabetical order in $({\cal W}^{(1)}_{\Psi})^{+}$),
and $w_{i,4} \in ({\cal W}^{(1)}_{\Psi})^{r_{i}}$ with $r_{i} \in
{\mathbb{N}}$ and $i = 1, \ldots, \tau(4)$. So,
$$
\begin{array}{lll}
[w_{\lambda,Y}(y_{1};1,2,3);w(p_{m_{\lambda,Y}(4)})] & = & \pm
[v_{1};w(p_{m_{\lambda,Y}(4)})] + [v_{2};w(p_{m_{\lambda,Y}(4)})]
~+ \\
& &
[w^{(1)}_{(1,2)};w(p_{m_{\lambda,Y}(4)})]+[w^{(2)}_{(1,2)};w(p_{m_{\lambda,Y}(4)})]+\\
& &
[w^{(1,2)}_{(1,2)};w(p_{m_{\lambda,Y}(4)})]+[\widetilde{w}_{(1,2)};w(p_{m_{\lambda,Y}(4)})].
\end{array}
$$
Note that $[w^{(1)}_{(1,2)};w(p_{m_{\lambda,Y}(4)})] \in
L^{m+m_{\lambda,Y}(4)}_{\rm grad, e}(W^{(1)}_{\Psi})$, $e \geq 2 +
\rho$, where $\rho$ is the minimum of the degrees of $w_{1,4},
\ldots, w_{\tau(4),4}$ in terms of elements of ${\cal
W}^{(1)}_{\Psi}$, and
$$
 \pm
[v_{1};w(p_{m_{\lambda,Y}(4)})] + [v_{2};w(p_{m_{\lambda,Y}(4)})]
+ [w^{(1)}_{(1,2)};w(p_{m_{\lambda,Y}(4)})] \in J \cap
L^{m+m_{\lambda,Y}(4)}_{\rm grad}(W^{(1)}_{\Psi}).
$$
Now either $v_{1} = \zeta_{1,4}$ or $v_{2} = \zeta_{1,4}$ or
$v_{1}, v_{2} \neq \zeta_{1,4}$. Let us assume that $v_{1}, v_{2}
\neq \zeta_{1,4}$. Then, by Lemma \ref{18a}, different Lyndon
words occur in the expressions of
$[v_{1};w(p_{m_{\lambda,Y}(4)})]$ and
$[v_{2};w(p_{m_{\lambda,Y}(4)})]$. We replace the (unique) Lyndon
polynomial corresponding to the smallest Lyndon word (by means of
the ordering of ${\cal W}^{(1)}_{\Psi}$) by
$[w_{\lambda,Y}(y_{1};1,2,3);w(p_{m_{\lambda,Y}(4)})]$. Similar
arguments may be applied if $v_{1} = \zeta_{1,4}$ or $v_{2} =
\zeta_{1,4}$. Furthermore, we deal with similar arguments the
cases either $m_{\lambda,Y}(2) \geq 1$ or $m_{\lambda,Y}(3) \geq
1$. We write $W_{\lambda,Y}(y_{1};1,m_{\lambda,Y}(2) +
m_{\lambda,Y}(3),4)$ for the $\mathbb{Z}$-submodule of
$J^{c+2}_{1,\lambda}$ spanned by all Lie commutators of the form
$[w_{\lambda,Y}(y_{1};1,2,3);w(p_{m_{\lambda,Y}(4)})]$ mentioned
above with $m_{\lambda,Y}(1) \geq 1$, $m_{\lambda,Y}(2) +
m_{\lambda,Y}(3) \geq 1$ and $m_{\lambda,Y}(4)$.

Next we assume that $m_{\lambda,Y}(2) = m_{\lambda,Y}(3) = 0$ and
$m_{\lambda,Y}(4) \geq 1$. Recall that
$$
[w_{\lambda,Y}(y_{1};1),w(p_{m_{\lambda,Y}(4)})] =
[w_{\lambda,Y}(y_{1};1), ~_{\mu_{1,4}}\zeta_{1,4}, \ldots,
~_{\mu_{\tau(4),4}}\zeta_{\tau(4),4}].
$$
Note that $w_{\lambda,Y}(y_{1};1) \in
L^{2+m_{\lambda,Y}(1)}(V_{1})$. Thus
$$
[w_{\lambda,Y}(y_{1};1),w(p_{m_{\lambda,Y}(4)})] = -
[\zeta_{1,4},w_{\lambda,Y}(y_{1};1), ~_{(\mu_{1,4}-1)}\zeta_{1,4},
\ldots, ~_{\mu_{\tau(4),4}}\zeta_{\tau(4),4}].
$$
If $m_{\lambda,Y}(5)$ or $m_{\lambda,Y}(6) \geq 1$, then, by Lemma
\ref{17b} (for $\kappa = \mu = 1)$, we have
$$
[w_{\lambda,Y}(y_{1};1),w(p_{m_{\lambda,Y}(4)}),w(p_{m_{\lambda,Y}(5)}),w(p_{m_{\lambda,Y}(6)})]
\in J^{c+2}_{C}.
$$
Thus we assume that $m_{\lambda,Y}(5) = m_{\lambda,Y}(6) = 0$. Let
$W_{\lambda,Y}(y_{1};1,4)$ be the $\mathbb{Z}$-submodule of
$J^{c+2}_{1,\lambda}$ spanned by all Lie commutators of the form
$[w_{\lambda,Y}(y_{1};1),w(p_{m_{\lambda,Y}(4)})]$. Note that the
Lie commutators of the form
$[w_{\lambda,Y}(y_{1};1),w(p_{m_{\lambda,Y}(4)})]$ are not
effected by the equations $(C1)$, and
$$
[w_{\lambda,Y}(y_{1};1),w(p_{m_{\lambda,Y}(4)})] \in J \cap
L^{2+m_{\lambda,Y}(1) +m_{\lambda,Y}(4)}_{\rm
grad}(W^{(1)}_{\Psi}).
$$
By the analysis of
$[y_{1},w(p_{m_{\lambda,Y}(1)}),w(p_{m_{\lambda,Y}(2)}),w(p_{m_{\lambda,Y}(3)})]$
in the equations (16)-(21) and the definition of
$W_{\lambda,Y}(y_{1};1,4)$ (having in mind Lemma \ref{17b} (for
$\kappa = \mu = 1$) and the ordering on ${\cal W}_{\Psi}$ (section
3.3)), we have the sum of
$W_{\lambda,Y}(y_{1};1,m_{\lambda,Y}(2)+m_{\lambda,Y}(3),4)$ and
$W_{\lambda,Y}(y_{1};1,4)$ is direct.

Let $m_{\lambda,Y}(4) = 0$ and $m_{\lambda,Y}(5) \geq 1$. If
$m_{\lambda,Y}(2)$ or $m_{\lambda,Y}(3) \geq 1$, then, by the
equations (16)-(21), we have
$$
[w_{\lambda,Y}(y_{1};t_{1},t_{2},t_{3}),w(p_{m_{\lambda,Y}(5)}),w(p_{m_{\lambda,Y}(6)})]
\in J^{c+2}_{C}.
$$
So, we assume that $m_{\lambda,Y}(2) = m_{\lambda,Y}(3) = 0$. By
Lemma \ref{17b} (for $\kappa = 2, \mu = 1$), we get
$$
[w_{\lambda,Y}(y_{1};1),w(p_{m_{\lambda,Y}(5)}),w(p_{m_{\lambda,Y}(6)})]
\in J^{c+2}_{C}.
$$
Finally, we assume that $(m_{\lambda,Y}(4) =)~~ m_{\lambda,Y}(5) =
0$. Suppose that $m_{\lambda,Y}(6) \geq 1$. If $m_{\lambda,Y}(2)$
or $m_{\lambda,Y}(3) \geq 1$, then, by the equations (16)-(21),
and since $L(\widetilde{W}_{\Psi,J})$ is an ideal $L(W_{\Psi})$,
we have
$$
[w_{\lambda,Y}(y_{1};t_{1},t_{2},t_{3}),w(p_{m_{\lambda,Y}(6)})]
\in J^{c+2}_{C}.
$$
Thus, $m_{\lambda,Y}(2) = m_{\lambda,Y}(3) = 0$. By the
construction of $L(\widetilde{W}_{\Psi,J})$, we get
$$
[w_{\lambda,Y}(y_{1};1),w(p_{m_{\lambda,Y}(6)})] \in J^{c+2}_{C}.
$$

\vskip .120 in

\emph{Subcase 1b. $m_{\lambda,Y}(1) = 0$.} This case is, in fact,
a special case of Subcase $1a$. Note that $w_{\lambda,Y}(y_{1};1)
= y_{1}$. Following the Subcase $1a$ step-by-step we have: Let
$m_{\lambda,Y}(2)$ or $m_{\lambda,Y}(3) \geq 1$. If
$m_{\lambda,Y}(5)$ or $m_{\lambda,Y}(6) \geq 1$, then the
equations (16)-(21), and, by the construction of
$L(\widetilde{W}_{\Psi,J})$ (It is an ideal in $L(W_{\Psi})$.),
the corresponding Lie commutator $[y_{1};w(p_{m_{\lambda,Y}(2)}),
\ldots, w(p_{m_{\lambda,Y}(6)})] \in J^{c+2}_{C}$. Thus we assume
that $m_{\lambda,Y}(5) = m_{\lambda,Y}(6) = 0$. Let
$m_{\lambda,Y}(4) \geq 1$. As in the Subcase $1a$, a similar
analysis of
$[y_{1},w(p_{m_{\lambda,Y}(2)}),w(p_{m_{\lambda,Y}(3)})]$ (in the
equations (16)-(21)), we obtain a $\mathbb{Z}$-submodule
$W_{\lambda,Y}(y_{1};m_{\lambda,Y}(2)+m_{\lambda,Y}(3),4)$ of
$J^{c+2}_{1,\lambda}$. Furthermore, a $\mathbb{Z}$-submodule
$W_{\lambda,Y}(y_{1};4)$ of $J^{c+2}_{1,\lambda}$ is constructed.

Define
$$
\begin{array}{lll}
V^{c+2}_{1,\lambda} & = & \sum_{Y \in {\cal G}_{\lambda}}(
W_{\lambda,Y}(y_{1};m_{\lambda,Y}(1),m_{\lambda,Y}(2)+m_{\lambda,Y}(3),4)
+ W_{\lambda,Y}(y_{1};1,4) ~+ \\
& &  W_{\lambda,Y}(y_{1};m_{\lambda,Y}(2)+m_{\lambda,Y}(3),4) +
W_{\lambda,Y}(y_{1};4)).
\end{array}
$$
We call the Lie commutators that generate the above
$\mathbb{Z}$-modules and are coming from the equations $(C1)$
\emph{the extended Lie commutators from the equations $(C1)$}.
Each extended Lie commutator produces an element of $J \cap
L^{c+2}_{\rm grad}(W^{(1)}_{\Psi})$. Such an element is called
\emph{the corresponding extended element from the equations
$(C1)$}. By the analysis in the Subcases $1a$ and $1b$,
$V^{c+2}_{1,\lambda}$ is the direct sum of the aforementioned
$\mathbb{Z}$-modules. Thus
$$
\begin{array}{lll}
V^{c+2}_{1,\lambda} & = & \bigoplus_{Y \in {\cal G}_{\lambda}}(
W_{\lambda,Y}(y_{1};m_{\lambda,Y}(1),m_{\lambda,Y}(2)+m_{\lambda,Y}(3),4)
\oplus W_{\lambda,Y}(y_{1};1,4) ~\oplus \\
& &  W_{\lambda,Y}(y_{1};m_{\lambda,Y}(2)+m_{\lambda,Y}(3),4)
\oplus W_{\lambda,Y}(y_{1};4)).
\end{array}
$$
Furthermore, it is clearly enough that $V^{c+2}_{1,\lambda} \cap
J^{c+2}_{C} = \{0\}$, and $J^{c+2}_{1,\lambda} +
J^{c+2}_{1,\lambda+1} + J^{c+2}_{C} = V^{c+2}_{1,\lambda} \oplus
(J^{c+2}_{1,\lambda+1} + J^{c+2}_{C})$. By applying an inverse
induction on the ${\rm Part}(c)$ (starting from $\lambda = (c)$),
we have
$$
J^{c+2}_{1} + J^{c+2}_{C} = U^{c+2}_{1} \oplus J^{c+2}_{C},
$$
where $U^{c+2}_{1} = \bigoplus_{\lambda \in {\rm Part}(c)}
V^{c+2}_{1,\lambda}$.

\vskip .120 in

\emph{Case 2. An analysis of $J^{c+2}_{2,\lambda}$ modulo
$(J^{c+2}_{2,\lambda+1}+J^{c+2}_{C})$.} Recall that, for $\lambda
\in {\rm Part}(c)$,
$$
J^{c+2}_{2,\lambda} + J^{c+2}_{C} =
[L^{2}(V_{2});\Phi_{\lambda}/\Phi_{\lambda+1}] +
J^{c+2}_{2,\lambda+1} + J^{c+2}_{C}.
$$
By Corollary \ref{co3}, we let
$[L^{2}(V_{2});\Phi_{\lambda}/\Phi_{\lambda+1}]$ be
$\mathbb{Z}$-spanned by all Lie commutators of the form
$[y_{2};w(p_{m_{\lambda,Y}(2)},p_{m_{\lambda,Y}(1)},
\ldots,p_{m_{\lambda,Y}(6)})]$ with $Y \in {\cal G}_{\lambda}$. We
separate two subcases, namely, $m_{\lambda,Y}(2) \geq 1$ and
$m_{\lambda,Y}(2) = 0$. Instead of the equations $(C1)$, we use
the equations $(C2)$. We point out that in order to construct
similar expressions as in the Case 1 (the equations (16)-(21)), we
apply the Jacobi identity several times. Following the analysis of
case $1$ step-by-step, we construct a $\mathbb{Z}$-submodule
$V^{c+2}_{2,\lambda}$ of $J^{c+2}_{2,\lambda}$ such that
$J^{c+2}_{2,\lambda} + J^{c+2}_{2,\lambda+1} + J^{c+2}_{C} =
V^{c+2}_{2,\lambda} \oplus (J^{c+2}_{2,\lambda+1} + J^{c+2}_{C})$.
In fact,
$$
\begin{array}{lll}
V^{c+2}_{2,\lambda} & = & \bigoplus_{Y \in {\cal G}_{\lambda}}(
W_{\lambda,Y}(y_{2};m_{\lambda,Y}(2),m_{\lambda,Y}(1)+m_{\lambda,Y}(3),4)
\oplus W_{\lambda,Y}(y_{2};2,4) ~\oplus \\
& &  W_{\lambda,Y}(y_{2};m_{\lambda,Y}(1)+m_{\lambda,Y}(3),4)
\oplus W_{\lambda,Y}(y_{2};4)).
\end{array}
$$
As in case 1, by applying an inverse induction on the ${\rm
Part}(c)$ (starting from $\lambda = (c)$), we have
$$
J^{c+2}_{2} + J^{c+2}_{C} = U^{c+2}_{2} \oplus J^{c+2}_{C},
$$
where $U^{c+2}_{2} = \bigoplus_{\lambda \in {\rm Part}(c)}
V^{c+2}_{2,\lambda}$.

\vskip .120 in

\emph{Case 3. An analysis of $J^{c+2}_{3,\lambda}$ modulo
$(J^{c+2}_{3,\lambda+1}+J^{c+2}_{C})$.} Recall that, for $\lambda
\in {\rm Part}(c)$,
$$
J^{c+2}_{3,\lambda} + J^{c+2}_{C} =
[L^{2}(V_{3});\Phi_{\lambda}/\Phi_{\lambda+1}] +
J^{c+2}_{3,\lambda+1} + J^{c+2}_{C}.
$$
By Corollary \ref{co3}, we let
$[L^{2}(V_{3});\Phi_{\lambda}/\Phi_{\lambda+1}]$ be
$\mathbb{Z}$-spanned by all Lie commutators of the form
$[y_{3};w(p_{m_{\lambda,Y}(3)},p_{m_{\lambda,Y}(1)},
\ldots,p_{m_{\lambda,Y}(6)})]$ with $Y \in {\cal G}_{\lambda}$. We
separate two subcases, namely, $m_{\lambda,Y}(3) \geq 1$ and
$m_{\lambda,Y}(3) = 0$. Instead of the equations $(C1)$, we use
the equations $(C3)$. We point out that in order to construct
similar expressions as in the Case 1 (the equations (16)-(21)), we
apply the Jacobi identity several times. Following the analysis of
case $1$ step-by-step, we construct a $\mathbb{Z}$-submodule
$V^{c+2}_{3,\lambda}$ of $J^{c+2}_{3,\lambda}$ such that
$J^{c+2}_{3,\lambda} + J^{c+2}_{3,\lambda+1} + J^{c+2}_{C} =
V^{c+2}_{3,\lambda} \oplus (J^{c+2}_{3,\lambda+1} + J^{c+2}_{C})$.
In fact,
$$
\begin{array}{lll}
V^{c+2}_{3,\lambda} & = & \bigoplus_{Y \in {\cal G}_{\lambda}}(
W_{\lambda,Y}(y_{3};m_{\lambda,Y}(3),m_{\lambda,Y}(1)+m_{\lambda,Y}(2),4)
\oplus W_{\lambda,Y}(y_{3};3,4) ~\oplus \\
& &  W_{\lambda,Y}(y_{3};m_{\lambda,Y}(1)+m_{\lambda,Y}(2),4)
\oplus W_{\lambda,Y}(y_{3};4)).
\end{array}
$$
As in case 1, by applying an inverse induction on the ${\rm
Part}(c)$ (starting from $\lambda = (c)$), we have
$$
J^{c+2}_{3} + J^{c+2}_{C} = U^{c+2}_{3} \oplus J^{c+2}_{C},
$$
where $U^{c+2}_{3} = \bigoplus_{\lambda \in {\rm Part}(c)}
V^{c+2}_{3,\lambda}$.

\vskip .120 in

\emph{Case 4. An analysis of $J^{c+2}_{\Psi,\lambda}$ modulo
$(J^{c+2}_{\Psi,\lambda+1}+J^{c+2}_{C})$.} Write ${\cal
V}^{(2)}_{1} = \{\psi_{2}(y_{6}), \psi_{2}(y_{7}),$ \linebreak
$\psi_{2}(y_{9}), \psi_{2}(y_{10})\}$ and ${\cal V}^{(2)}_{2} =
\{\psi_{2}(y_{12}), \psi_{2}(y_{15})\}$. That is, ${\cal V} =
{\cal V}^{(2)}_{1} \cup {\cal V}^{(2)}_{2}$ is a
$\mathbb{Z}$-basis for $W^{(2)}_{2,\Psi}$. Recall that
$[W^{(2)}_{2,\Psi};\Phi_{\lambda}/\Phi_{\lambda+1}]$ is the
$\mathbb{Z}$-submodule of $[J^{2};T^{c}]$ spanned by all Lie
commutators of the form $[v;w(p_{m_{\lambda,Y}(1)}, \ldots,
p_{m_{\lambda,Y}(6)})]$ with $v \in {\cal V}$ and $Y \in {\cal
G}_{\lambda}$. Thus
$$
[W^{(2)}_{2,\Psi};\Phi_{\lambda}/\Phi_{\lambda+1}] =
[[V_{2},V_{1}]^{(2)};\Phi_{\lambda}/\Phi_{\lambda+1}] +
[[V_{3},V_{1}]^{(2)};\Phi_{\lambda}/\Phi_{\lambda+1}] +
[[V_{3},V_{2}]^{(2)};\Phi_{\lambda}/\Phi_{\lambda+1}].
$$
By the definition of $J^{c+2}_{C}$, we have
$$
[[V_{2},V_{1}]^{(2)};\Phi_{\lambda}/\Phi_{\lambda+1}]+[[V_{3},V_{1}]^{(2)};\Phi_{\lambda}/\Phi_{\lambda+1}]
\subseteq J^{c+2}_{C}.
$$
Thus
$$
J^{c+2}_{\Psi,\lambda} + J^{c+2}_{\Psi,\lambda+1} + J^{c+2}_{C} =
[[V_{3},V_{2}]^{(2)};\Phi_{\lambda}/\Phi_{\lambda+1}] +
J^{c+2}_{\Psi,\lambda+1} + J^{c+2}_{C}.
$$
If $m_{\lambda,Y}(1) = 0$, then, by the definition of
$J^{c+2}_{C}$, each Lie commutator
$$
[v;w(p_{m_{\lambda,Y}(2)}, \ldots, p_{m_{\lambda,Y}(6)})] \in
J^{c+2}_{C}
$$
for all $v \in {\cal V}^{(2)}_{2}$. Thus we assume that
$m_{\lambda,Y}(1) \geq 1$. If $m_{\lambda,Y}(5)$ or
$m_{\lambda,Y}(6) \geq 1$, then, by the equations $(C3)-(C5)$, we
have
$$
[v,w(p_{m_{\lambda,Y}(1)}), \ldots ,w(p_{m_{\lambda,Y}(6)})] \in
J^{c+2}_{C}.
$$
Thus we assume that $m_{\lambda,Y}(5) = m_{\lambda,Y}(6) = 0$. For
$v \in {\cal V}_{2}^{(2)}$, let $\langle v \rangle$ denote the
cyclic $\mathbb{Z}$-module generated by $v$. Since
$[V_{3},V_{2}]^{(2)} = \bigoplus_{v \in {\cal V}_{2}^{(2)}}\langle
v \rangle$, we have
$$
[[V_{3},V_{2}]^{(2)};\Phi_{\lambda}/\Phi_{\lambda+1}] = \sum_{v
\in {\cal V}^{(2)}_{2}}[\langle v
\rangle;\Phi_{\lambda}/\Phi_{\lambda+1}]. \eqno(22)
$$
First we shall work with $v = \psi_{2}(y_{12})$. Similar arguments
may be applied to $v = \psi_{2}(y_{15})$. Suppose that $w_{1,1} =
x_{2}$. Then $w_{\ell}(p_{m_{\lambda,Y}(1)}) = x^{\mu_{1,1}}_{2}$
and $w(p_{m_{\lambda,Y}(1)}) = x^{\mu_{1,1}}_{2}$. By the
equations $(C4)$, the Jacobi identity and the equations (D), we
have
$$
\begin{array}{lll}
[\psi_{2}(y_{12});x^{\mu_{1,1}}_{2}] & = &
[\psi_{2}(y_{10}),~_{(\mu_{1,1}-1)}x_{2},x_{3}] -
[y_{8},~_{\mu_{1,1}}x_{2},x_{3}]-[y_{5},~_{(\mu_{1,1}-1)}x_{2},x_{5}] ~-  \\
& & [\psi_{2}(y_{6}),~_{(\mu_{1,1}-1)}x_{2},x_{5}] +
[\psi_{2}(y_{7}),~_{(\mu_{1,1}-1)}x_{2},x_{5}] + [y_{11}~_{(\mu_{1,1}-1)}x_{2},x_{3}] ~-  \\
& & [\psi_{2}(y_{6}),~_{(\mu_{1,1}-1)}x_{2},x_{6}] + [\psi_{2}(y_{7}),~_{(\mu_{1,1}-1)}x_{2},x_{6}]~- \\
& & [y_{5},~_{(\mu_{1,1}-1)}x_{2},x_{6}]+ w_{(12,2)} +
w_{(12,2,J)}, ~~~~~~~~~~~~~~~~~~~~~~~~~(23)
\end{array}
$$
where $w_{(12,2)} \in L^{2+m_{\lambda,Y}(1)}_{\rm grad,
e}(W^{(1)}_{\Psi})$, with $e \geq 2$, and $w_{(12,2,J)} \in
J_{C}$. Thus we assume that $w_{1,1} \neq x_{2}$. By using similar
arguments as before, we have
$$
\begin{array}{lllll}
[\psi_{2}(y_{12});w(p_{m_{\lambda,Y}(1)})] & = &
f(\psi_{2}(y_{9});w^{\prime}_{1,1}w^{(\mu_{1,1}-1)}_{1,1} \cdots
w^{\mu_{\tau(1),1}}_{\tau(1),1}x_{3}) - & & \\
& & f(y_{4};w^{\prime}_{1,1}w^{(\mu_{1,1}-1)}_{1,1} \cdots
w^{\mu_{\tau(1),1}}_{\tau(1),1}x_{5}) + & & \\
& & f(y_{8};w^{\prime}_{1,1}w^{(\mu_{1,1}-1)}_{1,1} \cdots
w^{\mu_{\tau(1),1}}_{\tau(1),1}x_{3}) + & & \\
& & f(y_{4};w^{\prime}_{1,1}w^{(\mu_{1,1}-1)}_{1,1} \cdots
w^{\mu_{\tau(1),1}}_{\tau(1),1}x_{6}) + & & \\
& & w_{(12,m_{\lambda,Y}(1))} + w_{(12,m_{\lambda,Y}(1),J)}, & &
~~~~~~~~~~~(24)
\end{array}
$$
where $w^{\prime}_{1,1}$ is the unique element in ${\cal
V}^{*}_{1}$ (the free monoid on ${\cal V}_{1}$) such that $w_{1,1}
= x_{1}w^{\prime}_{1,1}$, $w_{(12,m_{\lambda,Y}(1))} \in
L^{2+m_{\lambda,Y}(1)}_{\rm grad, e}(W^{(1)}_{\Psi})$ (with $e
\geq 2$) and $w_{(12,m_{\lambda,Y}(1),J)} \in J_{C}$. By the
equation (15) (for $t = 2,3$), the equation (23), the Jacobi
identity and the equations (D), we get
$$
\begin{array}{ll}
[[\psi_{2}(y_{12});x^{\mu_{1,1}}_{2}];w(p_{m_{\lambda,Y}(2)}),
w(p_{m_{\lambda,Y}(3)})] = &
f([\psi_{2}(y_{10}),x_{2}];x^{\mu_{1,1}-2}_{2}x_{3}w_{\ell}(p_{m_{\lambda,Y}(2)}),w_{\ell}(p_{m_{\lambda,Y}(3)}))
- \\
&
f([y_{8},x_{2}];x^{\mu_{1,1}-1}_{2}x_{3}w_{\ell}(p_{m_{\lambda,Y}(2)}),w_{\ell}(p_{m_{\lambda,Y}(3)}))
- \\
&
f([y_{5},x_{2}];x^{\mu_{1,1}-2}_{2}w_{\ell}(p_{m_{\lambda,Y}(2)}),x_{5}w_{\ell}(p_{m_{\lambda,Y}(3)}))
- \\
&
f([\psi_{2}(y_{6}),x_{2}];x^{\mu_{1,1}-2}_{2}w_{\ell}(p_{m_{\lambda,Y}(2)}),x_{5}w_{\ell}(p_{m_{\lambda,Y}(3)}))
+ \\
&
f([\psi_{2}(y_{7}),x_{2}];x^{\mu_{1,1}-2}_{2}w_{\ell}(p_{m_{\lambda,Y}(2)}),x_{5}w_{\ell}(p_{m_{\lambda,Y}(3)}))
+ \\
&
f([y_{11},x_{2}];x^{\mu_{1,1}-1}_{2}x_{3}w_{\ell}(p_{m_{\lambda,Y}(2)}),w_{\ell}(p_{m_{\lambda,Y}(3)}))
- \\
&
f([\psi_{2}(y_{6}),x_{2}];x^{\mu_{1,1}-2}_{2}w_{\ell}(p_{m_{\lambda,Y}(2)}),x_{6}w_{\ell}(p_{m_{\lambda,Y}(3)}))
+ \\
&
f([\psi_{2}(y_{7}),x_{2}];x^{\mu_{1,1}-2}_{2}w_{\ell}(p_{m_{\lambda,Y}(2)}),x_{6}w_{\ell}(p_{m_{\lambda,Y}(3)}))
- \\
&
f([y_{5},x_{2}];x^{\mu_{1,1}-2}_{2}w_{\ell}(p_{m_{\lambda,Y}(2)}),x_{6}w_{\ell}(p_{m_{\lambda,Y}(3)}))
+ \\
& w_{(12,1)} + w_{(12,1,J)},
~~~~~~~~~~~~~~~~~~~~~~~~~~~~~~~~~~~(25)
\end{array}
$$
where $w_{(12,1)} \in
L^{2+\mu_{1,1}+m_{\lambda,Y}(2)+m_{\lambda,Y}(3)}_{\rm grad,
e}(W^{(1)}_{\Psi})$, $e \geq 2$, and $w_{(12,1,J)} \in J_{C}$. By
applying similar arguments as above and using the equation (24),
we obtain a similar expression for
$[[\psi_{2}(y_{12});w(p_{m_{\lambda,Y}(1)})];w(p_{m_{\lambda,Y}(2)}),w(p_{m_{\lambda,Y}(3)})]$
(when $w_{1,1} \neq x_{2}$).

For the next few lines, we write
$$
w^{(12)}_{(1,2,3)} =
[[\psi_{2}(y_{12});w(p_{m_{\lambda,Y}(1)})];w(p_{m_{\lambda,Y}(2)}),
w(p_{m_{\lambda,Y}(3)})].
$$
For $m_{\lambda,Y}(4) = 0$, let
$W^{(12)}_{\lambda,Y,\Psi}(m_{\lambda,Y}(1),m_{\lambda,Y}(2),m_{\lambda,Y}(3))$
be the $\mathbb{Z}$-submodule of $J^{c+2}_{\Psi,\lambda}$ spanned
by all Lie commutators $w^{(12)}_{(1,2,3)}$. Finally, we assume
that $m_{\lambda,Y}(4) \geq 1$. The analysis of
$w^{(12)}_{(1,2,3)}$ allows us to assume that it is written as
$-v_{1}\pm v_{2}+v_{3}+\alpha v_{4}
+\widetilde{w}^{(12)}_{(1,2,3)}$ with $\alpha \in \{-1,0\}$,
$\widetilde{w}^{(12)}_{(1,2,3)} \in L^{c+2}_{\rm grad}(W_{\Psi})$,
$v_{1}, v_{2}, v_{3}, v_{4} \in {\cal W}^{(1)}_{\Psi}$ and $v_{1},
v_{2}, v_{3}, v_{4}$ do not occur in the expression of
$\widetilde{w}^{(12)}_{(1,2,3)}$. We proceed as in Subcase 1a.
Then either $v_{1} = \zeta_{1,4}$ or $v_{2} = \zeta_{1,4}$ or
$v_{3} = \zeta_{1,4}$ or $v_{4} = \zeta_{1,4}$ or $v_{1}, v_{2},
v_{3}, v_{4} \neq \zeta_{1,4}$. Let us assume that $v_{1}, v_{2},
v_{3}, v_{4} \neq \zeta_{1,4}$. Then, by Lemma \ref{18a},
different Lyndon words occur in the expressions of
$[v_{1};w(p_{m_{\lambda,Y}(4)})]$,
$[v_{2};w(p_{m_{\lambda,Y}(4)})]$,
$[v_{3};w(p_{m_{\lambda,Y}(4)})]$ and
$[v_{4};w(p_{m_{\lambda,Y}(4)})]$. We replace the (unique) Lyndon
polynomial corresponding to the smallest Lyndon word (by means of
the ordering of ${\cal W}^{(1)}_{\Psi}$) by
$[w^{(12)}_{(1,2,3)};w(p_{m_{\lambda,Y}(4)})]$. We point out the
above (different) Lyndon words do not appear in the expression of
$[\widetilde{w}^{(12)}_{(1,2,3)};w(p_{m_{\lambda,Y}(4)})]$.
Similar arguments may be applied if $v_{1} = \zeta_{1,4}$ or
$v_{2} = \zeta_{1,4}$, $v_{3} = \zeta_{1,4}$ or $v_{4} =
\zeta_{1,4}$. Furthermore, we deal with similar arguments the
cases either $m_{\lambda,Y}(2) \geq 1$ or $m_{\lambda,Y}(3) \geq
1$. We write $W^{(12)}_{\lambda,Y,\Psi}(1,m_{\lambda,Y}(2) +
m_{\lambda,Y}(3),4)$ for the $\mathbb{Z}$-submodule of
$J^{c+2}_{\Psi,\lambda}$ spanned by all Lie commutators of the
form $[w^{(12)}_{(1,2,3)};w(p_{m_{\lambda,Y}(4)})]$ mentioned
above with $m_{\lambda,Y}(1) \geq 1$, $m_{\lambda,Y}(2) +
m_{\lambda,Y}(3) \geq 1$ and $m_{\lambda,Y}(4)$.

Next we assume that $m_{\lambda,Y}(2) = m_{\lambda,Y}(3) = 0$.
Note that
$$
[\psi_{2}(y_{12}),w(p_{m_{\lambda,Y}(1)}),w(p_{m_{\lambda,Y}(4)})]
\in L_{\rm grad}(\widetilde{W}_{\Psi,J})
$$
and so, if $m_{\lambda,Y}(5) + m_{\lambda,Y}(6) \geq 0$, we have
$$
[\psi_{2}(y_{12}),w(p_{m_{\lambda,Y}(1)}),w(p_{m_{\lambda,Y}(4)}),w(p_{m_{\lambda,Y}(5)}),w(p_{m_{\lambda,Y}(6)})]
\in J^{c+2}_{C}.
$$
Define
$$
V^{c+2}_{12,\Psi,\lambda} = \sum_{Y \in {\cal G}_{\lambda}}(
W^{(12)}_{\lambda,Y,\Psi}(m_{\lambda,Y}(1),m_{\lambda,Y}(2),m_{\lambda,Y}(3))
+
W^{(12)}_{\lambda,Y,\Psi}(1,m_{\lambda,Y}(2)+m_{\lambda,Y}(3),4)).
$$
By the above analysis, $V^{c+2}_{12,\Psi,\lambda}$ is the direct
sum of the aforementioned $\mathbb{Z}$-modules. Thus
$$
V^{c+2}_{12,\Psi,\lambda} = \bigoplus_{Y \in {\cal G}_{\lambda}}(
W^{(12)}_{\lambda,Y,\Psi}(m_{\lambda,Y}(1),m_{\lambda,Y}(2),m_{\lambda,Y}(3))
\oplus
W^{(12)}_{\lambda,Y,\Psi}(1,m_{\lambda,Y}(2)+m_{\lambda,Y}(3),4)).
$$
It is clearly enough that $V^{c+2}_{12,\Psi,\lambda} \cap
J^{c+2}_{C} = \{0\}$. Following the same method (and giving the
same arguments) as above for the case $v = \psi_{2}(y_{15})$, and
using the analogous equations, we have
$$
V^{c+2}_{15,\Psi,\lambda} = \bigoplus_{Y \in {\cal G}_{\lambda}}(
W^{(15)}_{\lambda,Y,\Psi}(m_{\lambda,Y}(1),m_{\lambda,Y}(2),m_{\lambda,Y}(3))
\oplus
W^{(15)}_{\lambda,Y,\Psi}(1,m_{\lambda,Y}(2)+m_{\lambda,Y}(3),4)).
$$
Having in mind that the elements in the equations $(C3)-(C5)$ are
$\mathbb{Z}$-linearly independent, we may show that the above (in
all possible cases) extended Lie commutators from the equations
$(C3)-(C6)$ are $\mathbb{Z}$-linearly independent. This shows that
the sum $\sum_{Y \in {\cal G}_{\lambda}}(V^{c+2}_{12,\Psi,\lambda}
+ V^{c+2}_{15,\Psi,\lambda})$ is direct. Define
$$
V^{c+2}_{\Psi,\lambda} = \bigoplus_{Y \in {\cal
G}_{\lambda}}(V^{c+2}_{12,\Psi,\lambda} \oplus
V^{c+2}_{15,\Psi,\lambda}).
$$
Since $V^{c+2}_{\Psi,\lambda} \cap J^{c+2}_{C} = \{0\}$, we get
from the equation (22) and the above analysis,
$$
J^{c+2}_{\Psi,\lambda} + J^{c+2}_{\Psi,\lambda+1} + J^{c+2}_{C} =
V^{c+2}_{\Psi,\lambda} \oplus (J^{c+2}_{\Psi,\lambda+1} +
J^{c+2}_{C}).
$$
By applying an inverse induction on the ${\rm Part}(c)$ (starting
from $\lambda = (c)$), we have
$$
J^{c+2}_{\Psi} + J^{c+2}_{C} = U^{c+2}_{\Psi} \oplus J^{c+2}_{C},
$$
where $U^{c+2}_{\Psi} = \bigoplus_{\lambda \in {\rm Part}(c)}
V^{c+2}_{\Psi,\lambda}$.

Write $U^{c+2} = U^{c+2}_{1} + U^{c+2}_{2} + U^{c+2}_{3} +
U^{c+2}_{\Psi}$. Let ${\cal E}_{c+2}$ be the set consisting of the
extended Lie commutators (in all possible cases) which create
$U^{c+2}_{1}, U^{c+2}_{2}, U^{c+2}_{3}, U^{c+2}_{\Psi}$,
respectively, and the Lie commutators which are not coming from
the equations $(C1)-(C5)$. By the cases 1, 2 and 3, the ordering
on ${\cal W}^{(1)}_{\Psi}$ and having in mind that the elements in
the equations $(C1)-(C5)$ are $\mathbb{Z}$-linearly independent,
we may show that ${\cal E}_{c+2}$ is $\mathbb{Z}$-linearly
independent. So, $U^{c+2}$ is a direct sum. But
$$
\begin{array}{lll}
J^{c+2} & = & J^{c+2}_{1} + J^{c+2}_{2} + J^{c+2}_{3} + J^{c+2}_{\Psi} \\
& = & U^{c+2} + J^{c+2}_{C}.
\end{array}
$$
Since $U^{c+2} \cap J^{c+2}_{C} = \{0\}$, we have
$$
J^{c+2} = U^{c+2} \oplus J^{c+2}_{C}. \eqno(26)
$$
We claim that
$$
J^{c+2} = (J \cap L^{c+2}_{\rm grad}(W^{(1)}_{\Psi})) \oplus
J^{c+2}_{C}. \eqno(27)
$$
>From the analysis of the Lie commutators in the cases 1, 2, 3 and
4, it is clearly enough that
$$
U^{c+2} \subseteq (J \cap L^{c+2}_{\rm grad}(W^{(1)}_{\Psi}))
\oplus J^{c+2}_{C}.
$$
and so, we obtain the desired $\mathbb{Z}$-module decomposition of
$J^{c+2}$ in the equation (27). From the equations (26) and (27),
we get $J \cap L^{c+2}_{\rm grad}(W^{(1)}_{\Psi}) \cong U^{c+2}$
as $\mathbb{Z}$-modules. From the analysis of cases 1, 2, 3 and 4,
we have the corresponding extended elements from the equations
$(C1)-(C5)$ and the Lie commutators which are not coming from
$(C1)-(C5)$ belong to $J \cap L^{c+2}_{\rm grad}(W^{(1)}_{\Psi})$.
Write ${\cal E}^{\rm corr}_{c+2}$ for the set of all
aforementioned Lie elements which belong to $J \cap L^{c+2}_{\rm
grad}(W^{(1)}_{\Psi})$. Since $J \cap L^{c+2}_{\rm
grad}(W^{(1)}_{\Psi}) \cong U^{c+2}$ as (torsion-free)
$\mathbb{Z}$-modules and ${\cal E}^{\rm corr}_{c+2}$ is
$\mathbb{Z}$-linearly independent, we have $J \cap L^{c+2}_{\rm
grad}(W^{(1)}_{\Psi})$ is $\mathbb{Z}$-spanned by ${\cal E}^{\rm
corr}_{c+2}$. By the form of the elements of ${\cal E}^{\rm
corr}_{c+2}$, it is easy to construct a set consisting of Lie
polynomials in $L^{c+2}_{\rm grad}(W^{(1)}_{\Psi})$, and let
$(L^{c+2}_{\rm grad}(W^{(1)}_{\Psi}))^{*}$ denote the
$\mathbb{Z}$-module spanned by this set, such that
$$
L^{c+2}_{\rm grad}(W^{(1)}_{\Psi}) = (J \cap L^{c+2}_{\rm
grad}(W^{(1)}_{\Psi})) \oplus (L^{c+2}_{\rm
grad}(W^{(1)}_{\Psi}))^{*} \eqno(28)
$$
Since
$$
L^{c+2} = L^{c+2}_{\rm grad}(W^{(1)}_{\Psi}) \oplus J^{c+2}_{C},
$$
we obtain from the equations (28) and (27),
$$
L^{c+2} = (L^{c+2}_{\rm grad}(W^{(1)}_{\Psi}))^{*} \oplus J^{c+2}.
$$

\begin{corollary}
$L/J$ is torsion-free as $\Z$-module. Furthermore, $J$ is a free Lie algebra. 
\end{corollary}

\pf By the above theorem, for $c \geq 2$, $L^{c+2}/J^{c+2}$ is
torsion free $\mathbb{Z}$-module. Since $J = \bigoplus_{c \geq
0}J^{c+2}$ and $L/J \cong \bigoplus_{c \geq 0} (L^{c+2}/J^{c+2})$
as $\mathbb{Z}$-modules, we have $L/J$ is torsion-free. By a result of Witt (see \cite[Theorem 2.4.2.5]{baht}), we have $J$ is a free Lie algebra.

\section{The Lie algebra of $M_3$}

In this section, we deduce a presentation of $\L(M_{3})$. In
particular, we prove in Theorem \ref{3} that $L/J \cong \L(M_{3})$
as Lie algebras. Recall that $F$ denotes a free group of rank $6$
with a free generating set $\{a_{1}, \ldots, a_{6}\}$ with
ordering $a_{1} < a_{2} < \cdots < a_{6}$. Furthermore, ${\cal
R}_{\cal V} = \{r_{1}, \ldots, r_{9}\}$ and $M_{3} = F/N$ where
$N$ is generated as a group by the set $\{r^{g} = g^{-1}rg: r \in
{\cal R}_{\cal V}, g \in F\}$. Write $\widetilde{{\cal R}}_{\cal
V} = \{(r^{\pm 1},g): r \in {\cal R}_{\cal V}, g \in F \setminus
\{1\}\}$. Since $N \gamma_{3}(F)/\gamma_{3}(F)$ is generated by
the set $\{r\gamma_{3}(F): r \in {\cal R}_{\cal V}\}$,
$N\gamma_{3}(F)/\gamma_{3}(F) = J^{2}$ and since $L^{2} =
(L^{2})^{*} \oplus J^{2}$ as $\mathbb{Z}$-modules, where
$(L^{2})^{*} = W^{(1)}_{2,\Psi}$ (by Lemma \ref{le10a} (I)), we
have the set $\{r\gamma_{3}(F): r \in {\cal R}_{\cal V}\}$ is
$\mathbb{Z}$-linearly independent. Hence ${\cal R}_{\cal V} \cap
\widetilde{{\cal R}_{\cal V}} = \emptyset$ and so, $N$ is
generated by the disjoint union ${\cal R}_{\cal V} \cup
\widetilde{{\cal R}_{\cal V}}$.  Write $g = a^{\pm 1}_{i_{1}}
\cdots a^{\pm 1}_{i_{\mu}}$ as a reduced word in $F$. Since
$(a,bc) = (a,c)(a,b)(a,b,c)$, we may write $(r^{\pm 1},g)$ as a
product of group commutators of the form $(r^{\pm 1},a^{\pm
1}_{j_{1}}, \ldots, a^{\pm 1}_{j_{s}})$ with $s \geq 1$ and
$j_{1}, \ldots, j_{s} \in \{1, \ldots, 6\}$. Writing
$$
{\cal R}^{\prime}_{\cal V} = \{(r^{\pm 1},a^{\pm 1}_{j_{1}},
\ldots, a^{\pm 1}_{j_{s}}), r \in {\cal R}_{\cal V}, s \geq 1,
j_{1}, \ldots, j_{s} \in \{1, \ldots, 6\}\},
$$
we have $N$ is generated by the disjoint union ${\cal R}_{\cal V}
\cup {\cal R}^{\prime}_{\cal V}$. Since $N$ is a non-trivial
normal subgroup of $F$ and the index of $N$ in $F$ is not finite
(since $N \subseteq F^{\prime}$), we have by a result of
Nielsen-Schreier that $N$ is a free group of infinite rank (see,
also, \cite[Proposition 3.12]{lysc}).

For a positive integer $d$, let $N_{d} = N \cap \gamma_{d}(F)$.
Note that for $d \leq 2$, we have $N_{d} = N$. Also, for $d \geq
2$, $N_{d+1} = N_{d} \cap \gamma_{d+1}(F)$. Since $ (g_{1},
\ldots, g_{\kappa})^{f} = (g_{1}^{f}, \ldots, g_{\kappa}^{f})$ for
all $f \in F$, and $N$ is normal, we obtain $\{N_{d}\}_{d \geq 2}$
is a normal (descending) series of $N$. Clearly each $N_{d}$ is
normal in $F$. Since $(N_{\kappa}, N_{\ell}) \subseteq N_{\kappa +
\ell}$ for all $\kappa, \ell \geq 2$, we have $\{N_{d}\}_{d \geq
2}$ is a central series of $N$. Define ${\cal I}_{d}(N) = N_{d}
\gamma_{d+1}(F)/\gamma_{d+1}(F)$. It is easily verified that
${\cal I}_{d}(N) \cong N_{d}/N_{d+1}$ as $\mathbb{Z}$-modules. By
our definitions, identifications and the above discussion, ${\cal
I}_{2}(N) = N\gamma_{3}(F)/\gamma_{3}(F) = J^{2}$. Let ${\cal
B}^{*}_{2}$, ${\cal B}_{2,J}$ be the natural $\mathbb{Z}$-bases
for $(L^{2})^{*}$, $J^{2}$, respectively. By Lemma \ref{le10a}
(I), we have the union ${\cal B}^{*}_{2} \cup {\cal B}_{2,J}$ is
disjoint, and it is a $\mathbb{Z}$-basis of $L^{2}$. Note that we
may consider ${\cal B}^{*}_{2} = {\cal V}^{*}$ and ${\cal B}_{2,J}
= {\cal V}$. Since ${\cal R}_{\cal V}\gamma_{3}(F)/\gamma_{3}(F)
\subseteq {\cal I}_{2}(N)$ and ${\cal V}$ is a $\mathbb{Z}$-basis
of $J^{2}$, we have ${\cal R}_{\cal V}\gamma_{3}(F)/\gamma_{3}(F)$
is a $\mathbb{Z}$-linearly independent subset of ${\cal
I}_{2}(N)$. Thus
$$
r_{1}^{m_{1}} \cdots r_{9}^{m_{9}} \notin \gamma_{3}(F)\eqno(34)
$$
for any $m_{1}, \ldots, m_{9} \in {\mathbb{Z}}$ with $m_{1} +
\cdots + m_{9} \neq 0$. Write $N_{21}$ and $N_{22}$ for the
subgroups of $N$ generated by ${\cal R}_{\cal V}$ and $\{(r^{\pm
1},g_{1}, \ldots, g_{s}): r \in {\cal R}_{\cal V}, s \geq 1,
g_{1}, \ldots, g_{s} \in F \setminus \{1\}\}$, respectively. Since
$(ab,c) = (a,c)(a,c,b)(b,c)$ for all $a,b,c \in F$, we have
$N_{22}$ is normal in $N$. Since $N = N_{21} N_{22}$ and $N_{22}
\subseteq \gamma_{3}(F)$, we have, by the modular law,
$$
N_{3} = N \cap \gamma_{3}(F) = (N_{21} \cap \gamma_{3}(F)) N_{22}.
$$
By (34), we get $N_{21} \cap \gamma_{3}(F)$ is generated by
elements of length at least $4$ and so, $N_{21} \cap \gamma_{3}(F)
\subseteq N_{22}$. Thus $N_{3}$ is generated by the set $\{(r^{\pm
1},g_{1}, \ldots, g_{s}): r \in {\cal R}_{\cal V}, s \geq 1,
g_{1}, \ldots, g_{s} \in F \setminus \{1\}\}$. By our definitions
and identifications, $ {\cal I}_{3}(N) =
N_{3}\gamma_{4}(F)/\gamma_{4}(F) = J^{3}$. Write $(L^{3})^{*} =
(W^{(1)}_{3,\Psi})^{*}$. By Lemma \ref{le10a} (II)
$$
J^{3} = (\bigoplus^{3}_{i=2}[L^{2}(V_{1}),V_{i}]) \oplus
(\bigoplus^{3}_{i=1 \atop i\neq 2}[L^{2}(V_{2}),V_{i}]) \oplus
(\bigoplus_{i=1}^{2}[L^{2}(V_{3}),V_{i}]) \oplus
[[V_{3},V_{2}]^{(2)},V_{1}] \oplus J^{3}_{C}.
$$
Let ${\cal B}^{*}_{3}$, ${\cal B}_{3,J}$ be the natural
$\mathbb{Z}$-bases for $(L^{3})^{*}$, $J^{3}$, respectively. By
Lemma \ref{le10a} (III), we have the union ${\cal B}^{*}_{3} \cup
{\cal B}_{3,J}$ is disjoint, and it is a $\mathbb{Z}$-basis of
$L^{3}$. Hence we choose a $\mathbb{Z}$-basis ${\cal B}_{3,J} =
\{b^{(3)}_{1}, \ldots, b^{(3)}_{m(3)}\}$, with $m(3) = {\rm
rank}(J^{3})$, of $J^{3}$ consisting of Lie commutators of length
$3$ where at least one element of $\cal V$ occurs in each Lie
commutator of the basis elements. So, each basis element of
$J^{3}$ is written as a $\mathbb{Z}$-linear combination of Lie
commutators of degree $3$ of the above form $[\ell_{1}, \ell_{2},
\ell_{3}]$. By replacing the elements of ${\cal V}$ occurring in
$b_{\kappa}^{(3)}$ by elements of ${\cal R}_{\cal V}$, using
$\chi$, we may view each $b_{\kappa}^{(3)}$ as an element of
$N_{3}$. We write ${\cal B}_{3,J,f}$ for the set ${\cal B}_{3,J}$
viewed its elements as elements of $N_{3}$. (For example, $[x_{2},
x_{1}, x_{1}]$ is replaced by $(r_{1},a_{1})$ or,
$[\psi_{2}(y_{15}), x_{2}]-[\psi_{2}(y_{10}), x_{4}]+
[\psi_{2}(y_{6}), x_{2}]$ is replaced by $(r_{5}, a_{2})(r_{4},
a_{4})^{-1}(r_{6}, a_{2})$.) Since ${\cal
B}_{3,J,f}\gamma_{4}(F)/\gamma_{4}(F) \subseteq {\cal I}_{3}(N) =
J^{3}$ and ${\cal B}_{3,J}$ is a $\mathbb{Z}$-basis of $J^{3}$, we
have ${\cal B}_{3,J,f}\gamma_{4}(F)/\gamma_{4}(F)$ is a
$\mathbb{Z}$-linearly independent subset of ${\cal I}_{3}(N)$.
Thus
$$
(b_{1}^{(3)})^{n_{1}} \cdots (b_{m(3)}^{(3)})^{n_{m(3)}} \notin
\gamma_{4}(F)\eqno(35)
$$
for any $n_{1}, \ldots, n_{m(3)} \in \mathbb{Z}$ with $n_{1} +
\cdots + n_{m(3)} \neq 0$. Note that the set ${\cal
B}_{3,J,f}\gamma_{4}(F)/\gamma_{4}(F)$ is a $\mathbb{Z}$-basis of
${\cal I}_{3}(N)$. By (35), ${\cal B}_{3,J,f} \cap \gamma_{4}(F) =
\emptyset$. Write $N_{31}$ and $N_{32}$ for the subgroups of
$N_{3}$ generated by the sets ${\cal B}_{3,J,f}$ and $\{(r^{\pm
1},g_{1}, \ldots, g_{s}): r \in {\cal R}_{\cal V}, s \geq 2,
g_{1}, \ldots, g_{s} \in F \setminus \{1\}\}$, respectively. Note
that $N_{32}$ is normal in $N_{3}$ and so, $N_{31}N_{32}$ is a
subgroup of $N_{3}$. We claim that $N_{3} = N_{31}N_{32}$. It is
enough to show that $N_{3} \subseteq N_{31}N_{32}$. Furthermore,
it is enough to show that $(r^{\pm 1},a_{i}) \in N_{31}N_{32}$ for
$r \in {\cal R}_{\cal V}$, $i = 1, \ldots, 6$. Since
$(r^{-1},a_{i}) = (r,a_{i})^{-1}((r,a_{i})^{-1},r^{-1}) \in N_{3}$
and $N_{3}$ is normal in $F$, it is enough to show that $(r,a_{i})
\in N_{31}N_{32}$ for $r \in {\cal R}_{\cal V}$, $i = 1, \ldots,
6$. For $i = 1, \ldots, 6$ and $r \in {\cal R}_{\cal V}$, there
are unique $n_{1i}, \ldots, n_{m(3)i} \in {\mathbb{Z}}$ such that
$$
(r,a_{i}) \gamma_{4}(F) = (b_{1}^{(3)})^{n_{1i}} \cdots
(b_{m(3)}^{(3)})^{n_{m(3)i}} \gamma_{4}(F).
$$
In fact, working the elements $b_{1}^{(3)}, \ldots,
b_{m(3)}^{(3)}$, appearing in the above product, in $N_{31}N_{32}
\subseteq \gamma_{3}(F)$ (Since the elements $b_{1}^{(3)}, \ldots,
b_{m(3)}^{(3)}$ are regarded in $N_{3}$.) and using the identity
$ab = ba(a,b)$, we form the product $(b_{1}^{(3)})^{n_{1i}} \cdots
(b_{m(3)}^{(3)})^{n_{m(3)i}}$ in such a way
$$
(r,a_{i}) = (b_{1}^{(3)})^{n_{1i}} \cdots
(b_{m(3)}^{(3)})^{n_{m(3)i}} v
$$
with $v \in \gamma_{4}(F)$. Following the procedure, it is clearly
enough that $v \in N_{31}N_{32}$. Thus $(r,a_{i}) \in
N_{31}N_{32}$ for all $r \in {\cal R}_{\cal V}$ and $i \in \{1,
\ldots, 6\}$. Therefore $N_{3} = N_{31}N_{32}$. Since $N_{3} =
N_{31} N_{32}$ and $N_{32} \subseteq \gamma_{4}(F)$, we have, by
the modular law,
$$
N_{4} = N_{3} \cap \gamma_{4}(F) = N_{31}N_{32} \cap \gamma_{4}(F)
= (N_{31} \cap \gamma_{4}(F)) N_{32}.
$$
By (35), we get $N_{31} \cap \gamma_{4}(F)$ is generated by
elements of length at least $6$. Having in mind the way in which
$b_{1}^{(3)}, \ldots, b_{m(3)}^{(3)}$ are regarded elements in
$N_{3}$ and the definition of $N_{32}$, we have $N_{31} \cap
\gamma_{4}(F) \subseteq N_{32}$. Thus $N_{4}$ is generated by the
set $\{(r^{\pm 1}, g_{1}, \ldots, g_{s}): r \in {\cal R}, s \geq
2, g_{1}, \ldots, g_{s} \in F \setminus \{1\}\}$ and so, ${\cal
I}_{4}(N) = J^{4}$.

\begin{proposition}\label{1}
For a positive integer $c$, $N_{c+2}$ is generated by the set
$\{(r^{\pm 1}, g_{1},$ $\ldots, g_{s}): r \in {\cal R}, s \geq c,
g_{1},\ldots, g_{s} \in F \setminus \{1\}\}$. Furthermore, ${\cal
I}_{c+2}(N) = J^{c+2}$ for all $c \geq 1$.
\end{proposition}

\pf We proceed by induction on $c$, with $c \geq 1$. We have
already shown our claim for $c = 1,2$. Thus we assume that $
N_{c+2}$, with $c \geq 2$, is generated by the set $\{(r^{\pm 1},
g_{1}, \ldots, g_{s}): r \in {\cal R}, s \geq c, g_{1}, \ldots,
g_{s} \in F \setminus \{1\}\}$ and so, by our definitions and
identifications, ${\cal I}_{c+2}(N) = J^{c+2}$. By the proof of
Theorem \ref{th1}, there exists a $\mathbb{Z}$-module
$(L^{c+2})^{*}$, say, such that $L^{c+2} = (L^{c+2})^{*} \oplus
J^{c+2}$ as $\mathbb{Z}$-modules. By the equation (3) and Theorem
\ref{th1}, we have, for $c \geq 2$,
$$
J^{c+2} = U^{c+2}_{1} \oplus U^{c+2}_{2} \oplus U^{c+2}_{3} \oplus
U^{c+2}_{\Psi} \oplus L^{c+2}(V_{1}) \oplus L^{c+2}(V_{2}) \oplus
L^{c+2}(V_{3}) \oplus L^{c+2}_{\rm grad}(W^{(2)}_{\Psi}) \oplus
L^{c+2}_{\rm grad}(\widetilde{W}_{\Psi,J}).
$$
Let ${\cal B}^{*}_{c+2}, {\cal B}_{c+2,J}$ be $\mathbb{Z}$-bases
of $(L^{c+2})^{*}, J^{c+2}$, respectively. (By using the proof of
Theorem \ref{th1}, and the chosen $\mathbb{Z}$-bases of
$L^{c+2}(V_{1}), L^{c+2}(V_{2}), L^{c+2}_{\rm
grad}(W^{(2)}_{\Psi})$ and $L^{c+2}_{\rm
grad}(\widetilde{W}_{\Psi,J})$, respectively.) Thus we have the
union ${\cal B}^{*}_{c+2} \cup {\cal B}_{c+2,J}$ is disjoint, and
it is a $\mathbb{Z}$-basis of $L^{c+2}$. Hence we choose a
$\mathbb{Z}$-basis ${\cal B}_{c+2,J} = \{b^{(c+2)}_{1}, \ldots,
b^{(c+2)}_{m(c+2)}\}$, with $m(c+2) = {\rm rank}(J^{c+2})$,  of
$J^{c+2}$ consisting of Lie commutators of length $c+2$ where at
least one element of $\cal V$ occurs in each Lie commutator of the
basis elements. So, each basis element of $J^{c+2}$ is written as
a $\mathbb{Z}$-linear combination of Lie commutators of length
$c+2$ of the aforementioned form $[\ell_{1}, \ldots, \ell_{c+2}]$.
By replacing (as in the case $c=2$) the elements of ${\cal V}$
occurring in $b_{\kappa}^{(c+2)}$ by the corresponding elements of
${\cal R}_{\cal V}$, we may view each $b_{\kappa}^{(c+2)}$ as an
element of $N_{c+2}$. We write ${\cal B}_{c+2,J,f}$ for the set
${\cal B}_{c+2,J}$ viewed its elements as elements of $N_{c+2}$.
Since ${\cal B}_{c+2,J,f}\gamma_{c+3}(F)/\gamma_{c+3}(F) \subseteq
{\cal I}_{c+2}(N) = J^{c+2}$ and ${\cal B}_{c+2,J}$ is a
$\mathbb{Z}$-basis of $J^{c+2}$, we have ${\cal
B}_{c+2,J,f}\gamma_{c+2}(F)/\gamma_{c+2}(F)$ is a
$\mathbb{Z}$-linearly independent subset of ${\cal I}_{c+2}(N)$.
Thus
$$
(b_{1}^{(c+2)})^{n_{1}} \cdots (b_{m(c+2)}^{(c+2)})^{n_{m(c+2)}}
\notin \gamma_{c+3}(F)\eqno(36)
$$
for any $n_{1}, \ldots, n_{m(c+2)} \in \mathbb{Z}$ with $n_{1} +
\cdots + n_{m(c+2)} \neq 0$. The set ${\cal
B}_{c+2,J,f}\gamma_{c+3}(F)/\gamma_{c+3}(F)$ is a
$\mathbb{Z}$-basis of ${\cal I}_{c+2}(N)$. By (36), ${\cal
B}_{c+2,f} \cap \gamma_{c+3}(F) = \emptyset$. Write $N_{(c+2)1}$
and $N_{(c+2)2}$ for the subgroups of $N_{c+2}$ generated by the
sets ${\cal B}_{c+2,J,f}$ and $\{(r^{\pm 1},g_{1}, \ldots, g_{s}):
r \in {\cal R}_{\cal V}, s \geq c+1, g_{1}, \ldots, g_{s} \in F
\setminus \{1\}\}$, respectively. Note that $N_{(c+2)2}$ is normal
in $N_{c+2}$. We claim that $N_{c+2} = N_{(c+2)1}N_{(c+2)2}$. It
is enough to show that $N_{c+2} \subseteq N_{(c+2)1}N_{(c+2)2}$.
Furthermore, it is enough to show that $(r^{\pm 1},a_{i_{1}},
\ldots,a_{i_{c}}) \in N_{(c+2)1}N_{(c+2)2}$ for $r \in {\cal
R}_{\cal V}$, $i = 1, \ldots, 6$. Since
$$
(a^{-1},b,c) = ((a,b)^{-1},c)(((a,b)^{-1},c),((a,b)^{-1},a^{-1}))
(((a,b)^{-1},a^{-1}),c)
$$
and $N_{c+2}$ is normal in $F$, it is enough to show that
$(r,a_{i_{1}}, \ldots,a_{i_{c}}) \in N_{(c+2)1}N_{(c+2)2}$ for $r
\in {\cal R}_{\cal V}$, $i = 1, \ldots, 6$. For $i = 1, \ldots, 6$
and $r \in {\cal R}_{\cal V}$, there are unique $n_{1i_{1}\cdots
i_{c}}, \ldots, n_{m(c+2)i_{1} \cdots i_{c}} \in {\mathbb{Z}}$
such that
$$
(r,a_{i_{1}}, \ldots, a_{i_{c}}) \gamma_{c+3}(F) =
(b_{1}^{(c+2)})^{n_{1i_{1} \cdots i_{c}}} \cdots
(b_{m(c+2)}^{(c+2)})^{n_{m(c+2)i_{1} \cdots i_{c}}}
\gamma_{c+3}(F).
$$
As in the case $c=2$, working the elements $b_{1}^{(c+2)}, \ldots,
b_{m(c+2)}^{(c+2)}$ in $N_{(c+2)1}N_{(c+2)2} \subseteq
\gamma_{c+2}(F)$ (Since the elements $b_{1}^{(c+2)}, \ldots,
b_{m(c+2)}^{(c+2)}$ are regarded in $N_{c+2}$.) and using the
identity $ab = ba(a,b)$, we form the product
$(b_{1}^{(c+2)})^{n_{1i_{1} \cdots i_{c}}} \cdots
(b_{m(c+2)}^{(c+2)})^{n_{m(c+2)i_{1} \cdots i_{c}}}$ in such a way
$$
(r,a_{i_{1}}, \ldots, a_{i_{c}}) = (b_{1}^{(c+2)})^{n_{1i_{1}
\cdots i_{c}}} \cdots (b_{m(c+2)}^{(c+2)})^{n_{m(c+2)i_{1} \cdots
i_{c}}} v
$$
with $v \in \gamma_{c+3}(F)$. Following the procedure, it is
clearly enough that $v \in N_{(c+2)1}N_{(c+2)2}$. Thus
$(r,a_{i_{1}}, \ldots, a_{i_{c}}) \in N_{(c+2)1}N_{(c+2)2}$ for
all $r \in {\cal R}_{\cal V}$ and $i \in \{1, \ldots, 6\}$.
Therefore $N_{c+2} = N_{(c+2)1}N_{(c+2)2}$. Since $N_{c+2} =
N_{(c+2)1} N_{(c+2)2}$ and $N_{(c+2)2} \subseteq \gamma_{c+3}(F)$,
we have, by the modular law,
$$
N_{c+3} = N_{c+2} \cap \gamma_{c+3}(F) = (N_{(c+2)1}N_{(c+2)2})
\cap \gamma_{c+3}(F) = (N_{(c+2)1} \cap \gamma_{c+3}(F))
N_{(c+2)2}.
$$
By (36), we get $N_{(c+2)1} \cap \gamma_{c+3}(F)$ is generated by
elements of length at least $2(c+1)$. Having in mind the way in
which $b_{1}^{(c+2)}, \ldots, b_{m(c+2)}^{(c+2)}$ are regarded
elements in $N_{c+3}$ and the definition of $N_{(c+2)2}$, we have
$N_{(c+2)1} \cap \gamma_{c+3}(F) \subseteq N_{(c+2)2}$. Thus
$N_{c+3}$ is generated by the set $\{(r^{\pm 1}, g_{1}, \ldots,
g_{s}): r \in {\cal R}, s \geq c+1, g_{1}, \ldots, g_{c+1} \in F
\setminus \{1\}\}$ and so, ${\cal I}_{c+3}(N) = J^{c+3}$. \qed

\vskip .120 in

Since $F$ is residually nilpotent, we have $ \bigcap_{d \geq
2}N_{d} = \{1\}$. Since $N \subseteq F^{\prime}$, we get ${\cal
I}_{1}(N) = 0$. Since ${\cal I}_{d}(N) \cong N_{d}/N_{d+1}$ as
$\mathbb{Z}$-modules for all $d \geq 2$, we have from Proposition
\ref{1}, $N_{d} \neq N_{d+1}$ for all $d \geq 2$. Define
$$
{\cal I}(N) = \bigoplus_{d \geq 2} N_{d}
\gamma_{d+1}(F)/\gamma_{d+1}(F) = \bigoplus_{d \geq 2}{\cal
I}_{d}(N).
$$
Since $N$ is a normal subgroup of $F$, we have ${\cal I}(N)$ is an
ideal of $L$ (see \cite{laz}).

\begin{corollary}
${\cal I}(N) = J$.
\end{corollary}

\pf Since $J = \bigoplus_{d \geq 2}J^{d}$ and ${\cal I}_{2}(N) =
J^{2}$, we have from Proposition \ref{1} that ${\cal I}(N) = J$.
\qed

\vskip .120 in

Our next result gives a presentation of $\L(M_{3})$.

\begin{theorem}\label{3} $L/J \cong \L(M_{3})$ as Lie
algebras. Consequently, $\L(M_3)$ is a torsion-free Lie algebra.
\end{theorem}

\emph{Proof.} Recall that
$$
\L(M_{3}) = \bigoplus_{c \geq 1}
\gamma_{c}(M_{3})/\gamma_{c+1}(M_{3}).
$$
Since $M_{3}/M^{\prime}_{3} \cong F/N F^{\prime} = F/F^{\prime}$,
we have $\L(M_{3})$ is generated as a Lie algebra by the set
$\{\alpha_{i}: i = 1, \ldots, 6\}$ with $\alpha_{i} =
a_{i}M^{\prime}_{3}$. Since $L$ is a free Lie algebra of rank $6$
with a free generating set $\{x_{1}, \ldots, x_{6}\}$, the map
$\psi$ from $L$ into $\L(M_{3})$ satisfying the conditions
$\psi(x_{i}) = \alpha_{i}$, $i = 1, \ldots, 6$, extends uniquely
to a Lie algebra homomorphism. Since $\L(M_{3})$ is generated as a
Lie algebra by the set $\{\alpha_{i}: i = 1, \ldots, 6\}$, we have
$\psi$ is onto. Hence $L/{\rm Ker}\psi \cong \L(M_{3})$ as Lie
algebras. By definition, $J \subseteq {\rm Ker} \psi$, and so
$\psi$ induces a Lie algebra epimorphism $\overline{\psi}$ from
$L/J$ onto $\L(M_{3})$. In particular, $\overline{\psi}(x_{i} + J)
= \alpha_{i}$, $i = 1, \ldots, 6$. Note that $\overline{\psi}$
induces $\overline{\psi}_{c}$, say, a $\mathbb{Z}$-linear mapping
from $(L^{c} + J)/J$ onto $\gamma_{c}(M_{3})/\gamma_{c+1}(M_{3})$.
For $c \geq 2$,
$$
\gamma_{c}(M_{3})/\gamma_{c+1}(M_{3}) \cong
\gamma_{c}(F)\gamma_{c+1}(F)N/\gamma_{c+1}(F)N \cong
\gamma_{c}(F)/(\gamma_{c}(F) \cap \gamma_{c+1}(F)N).
$$
Since $\gamma_{c+1}(F) \subseteq \gamma_{c}(F)$, we have by the
modular law,
$$
\gamma_{c}(F)/(\gamma_{c}(F) \cap \gamma_{c+1}(F)N) =
\gamma_{c}(F)/\gamma_{c+1}(F)N_{c}.
$$
But, by Proposition \ref{1},
$$
\gamma_{c}(F)/\gamma_{c+1}(F)N_{c} \cong
(\gamma_{c}(F)/\gamma_{c+1}(F))/{\cal I}_{c}(N) \cong L^{c}/J^{c}.
$$
Therefore
$$
\gamma_{c}(M_{3})/\gamma_{c+1}(M_{3}) \cong L^{c}/J^{c} \cong
(L^{c})^{*}.
$$
Hence ${\rm rank}(\gamma_{c}(M_{3})/\gamma_{c+1}(M_{3})) = {\rm
rank}(L^{c})^{*}$. Since $J = \bigoplus_{m \geq 2}J^{m}$, we have
$(L^{c} + J)/J \cong L^{c}/(L^{c} \cap J) = L^{c}/J^{c} \cong
(L^{c})^{*}$ (by Theorem \ref{th1}), we obtain ${\rm
Ker}\overline{\psi}_{c}$ is torsion-free. Since ${\rm
rank}(\gamma_{c}(M_{3})/\gamma_{c+1}(M_{3})) = {\rm
rank}(L^{c})^{*}$, we have ${\rm Ker}\overline{\psi}_{c} = \{1\}$
and so, $\overline{\psi}_{c}$ is isomorphism. Since
$\overline{\psi}$ is epimorphism and each $\overline{\psi}_{c}$ is
isomorphism, we have $\overline{\psi}$ is isomorphism. Hence $L/J
\cong \L(M_{3})$ as Lie algebras. \qed

\begin{corollary}
$\L(M_3)$ admits the presentation in \cite{cpvw} described by the
original presentation of $M_3$.
\end{corollary}

\bigskip

\noindent V. Metaftsis, Department of Mathematics, University of
the Aegean, Karlovassi, 832 00 Samos, Greece. {\it e-mail:}
vmet@aegean.gr

\bigskip

\noindent A.I. Papistas, Department of Mathematics, Aristotle
University of Thessaloniki, 541 24 Thessaloniki, Greece. {\it
e-mail:} apapist@math.auth.gr


\begin{thebibliography}{9}

\small

\bibitem{andreadakis} S. Andreadakis, On the automorphisms of free groups and free nilpotent groups.
{\it Proc. London Math. Soc.} {\bf 15} (1965), 239--268.

\bibitem{baht} Yu.A. Bakhturin, Identical Relations in Lie Algebras, Nauka, Moscow, 1985 (in Russian).
English translation: VNU Science Press, Utrecht, 1987.

\bibitem{barmik} V.G. Bardakov, R. Mikhailov, On certain questions of the free group automorphisms theory.
{\it Comm. Algebra} {\bf 36} (2008), no. 4, 1489--1499.

\bibitem{bour} N. Bourbaki, Lie Groups and Lie Algebras Part I, Hermann, Paris, 1987 (Chapters 1-3).

\bibitem{bks1} R.M. Bryant, L.G. Kov\'{a}cs, Ralph St\"{o}hr, Lie powers of modules for groups of prime order,
{\it Proc. London Math. Soc.} {\bf 84} (2002), 343--374.

\bibitem{bks2} R.M. Bryant, L.G. Kov\'{a}cs, Ralph St\"{o}hr, Invariant bases for free Lie rings,
{\it Q. J. Math.} {\bf 53} (2002), no. 1, 1--17.

\bibitem{brmic} R.M. Bryant, I.C. Michos, Lie powers of free modules for certain groups of prime power order,
{\it J. Austral. Math. Soc.} {\bf 71} (2001), 149--158.

\bibitem{bsc} R.M. Bryant, M. Schocker, The decomposition of Lie powers,
{\it Proc. London Math. Soc.} {\bf 93} (2006), 175--196.

\bibitem{cpvw} F.R. Cohen, J. Pakianathan, V.V. Vershinin and J. Wu, Basis-conjugating automorphisms of a free group and
associated Lie algebras. {\it Geometry \& Topology Monographs}
{\bf 13} (2008) 147--168.

\bibitem{fr} M. Falk, R. Randell , The lower central series of a fiber-type arrangement, {\it Invent. Math.}
{\bf 82} (1985) 77--88.

\bibitem{jac} N. Jacobson, Lie algebras, Interscience, New York, 1962.

\bibitem{kohno} T. Kohno , S\'erie de Poincar\'e-Koszul associ\'ee aux groupes de tresses pures, {\it Invent.
Math.} {\bf 82} (1985) 57--75.

\bibitem{laz} M. Lazard, Sur les groupes nilpotents et les anneaux de Lie, \emph{Ann. Sci. Ecole Norm. Sup.} (3) \textbf{71} (1954) 101-190.

\bibitem{loth} M. Lothaire, Combinatorics on words, Encyclopedia
of Mathematics and its applications v. 17, Cambridge University
Press, 1997.

\bibitem{lysc} R. Lyndon, P. Schupp, Combinatorial Group Theory, Ergebnisse der Mathematik und ihner Grenzgebiete 89, Springer-Verlag, 1977.

\bibitem{magnus} W. Magnus, \"{U}ber $n$-dimensinale Gittertransformationen, \emph{Acta Math.}, 64 (1935), 353-367.

\bibitem{mccool} J. McCool, On basis-conjugating automorphisms of free groups, {\it Can. J. Math.}
{\bf 38} (1986), 1525--1529.
\bibitem{nielsen} J. Nielsen, Die Isomorphismengruppe der freien Gruppen, \emph{Math. Annalen} 91 (1924), 169-209.

\bibitem{reut} C. Reutenauer, Free Lie algebras, London
Mathematical Society Monographs, new series, no 7, Oxford
University Press, 1993.

\bibitem{shme} A.L. Shmel'kin, \emph{Free polynilpotent groups}, Izv. Akad. Nauk. SSSR, Ser Mat. 28 (1964) 91-122.
\end{thebibliography}
\end{document}